\documentclass[fixNSE]{nseJournal}
\usepackage{style}

\title{High-Order Finite Element Second Moment Methods for Linear Transport}
\addAuthor{\correspondingAuthor{Samuel Olivier}}{a}
\correspondingEmail{solivier@lanl.gov}
\addAuthor{Terry S. Haut}{b}

\addAffiliation{a}{Los Alamos National Laboratory, Los Alamos, NM}
\addAffiliation{b}{Lawrence Livermore National Laboratory, Livermore, CA}

\addKeyword{moment methods}
\addKeyword{high-order low-order}
\addKeyword{second moment method}
\addKeyword{preconditioned iterative solvers}

\begin{document}
\titlePage

\begin{abstract}
We present high-order, finite element-based Second Moment Methods (SMMs) for solving radiation transport problems in two spatial dimensions. 
We leverage the close connection between the Variable Eddington Factor (VEF) method and SMM to convert existing discretizations of the VEF moment system to discretizations of the SMM moment system. 
The moment discretizations are coupled to a high-order Discontinuous Galerkin discretization of the \Sn transport equations. 
We show that the resulting methods achieve high-order accuracy on high-order (curved) meshes, preserve the thick diffusion limit, and are effective on a challenging multi-material problem both in outer fixed-point iterations and in inner preconditioned iterative solver iterations for the discrete moment systems. 
We also present parallel scaling results and provide direct comparisons to the VEF algorithms the SMM algorithms were derived from. 
\end{abstract}

\section{Introduction}
The Second Moment Method (SMM), introduced by \citet{lewis_miller} in 1976, is an iterative, moment-based algorithm for solving the Boltzmann transport equation in applications such as nuclear reactor design, medical physics, and high energy density physics. 
The original concept was to iteratively couple the transport equation to a diffusion equation with linear, transport-dependent source terms that vanish when the solution is a linear function in angle \cite{AL}. 
The transport-dependent source terms ``correct'' the diffusion approximation such that, upon iterative convergence, the corrected diffusion equation can reproduce the transport solution. 
Using the corrected diffusion solution to compute slow-to-converge physics, such as scattering, provides iterative acceleration, leading to a robust algorithm where the cost scales independent of the mean free path. 

SMM can also be viewed as a member of a broader class of multi-scale methods, known as moment or high-order low-order (HOLO) methods, developed to solve high-dimensional kinetic equations in the context of multiphysics simulations (see \citet{CHACON201721} for a review of their 60 year history) where the kinetic equation is iteratively coupled to a small number of its moments through discrete closures. 
In the case of SMM, the corrected diffusion equation can be viewed as a two-moment moment system closed with additive closures. 
Such methods are automatically asymptotic preserving (i.e.~preserve the thick diffusion limit), are conservative even when a non-conservative discretization of the kinetic equation is used, and can be directly coupled to other physics components in the simulation, isolating the expensive, high-dimensional kinetic equation from the evolution of stiff multiphysics. 
Furthermore, the rapid and robust convergence of moment methods is maintained even when the kinetic and moment equations are discretized independently \cite{two-level-independent-warsa}. 
These so-called independent moment methods \cite{doi:10.1080/00411459308203810} allow the flexibility to choose the discretization of the moment system to meet the requirements of the overall algorithm such as computational efficiency and multiphysics compatibility. 

SMM is particularly closely related to another moment method known as the Variable Eddington Factor (VEF) \cite{mihalas} or  Quasidiffusion \cite{goldin} method. 
The SMM and VEF algorithms are differentiated only by their choice of closure: SMM uses additive, linear closures while VEF uses multiplicative, nonlinear closures. 
VEF's multiplicative closures affect the moment operator (i.e.~the left had side of the linear system $\mat{A}x = b$), resulting in a non-self-adjoint moment system that has been difficult to discretize and precondition effectively (see \citet{dgvef_olivier} and the references therein). 
SMM's use of additive closures results in a moment system where the closures appear as source terms only, leading to clear computational advantages over VEF. 
In particular, the self-adjoint structure of the diffusion approximation is preserved, allowing the use of simpler iterative schemes, such as conjugate gradient, to solve the discretized moment system. 
By contrast, VEF's generally non-symmetric discrete moment system requires the use of general purpose solvers, such as GMRES or BiCGStab, which typically require more computation and memory. 
SMM can also directly leverage existing discretization and preconditioning strategies developed for standard elliptic problems that would not be possible for the unusual structure of the VEF moment system. 
Furthermore, where the VEF closures are well defined only when the transport solution is strictly positive, the SMM closures are valid for any transport solution, including solutions that possess negativities arising in under-resolved calculations. 

In the process of performing a Fourier stability analysis of the VEF method, \citet{cefus} found that SMM can be viewed as a VEF method that has been linearized about a linearly anisotropic solution and that the SMM and VEF algorithms had similar iterative convergence properties. 
We stress that even though SMM can be viewed as a linearized VEF method, the SMM moment system is still simply an algebraic reformulation of the transport equation and is thus able to achieve the full transport physics even when the problem is not diffusive. 
The methods' close iterative convergence rates indicate that SMM has the potential to converge as quickly as VEF while using a moment system that is less expensive to invert at each iteration and more robust to negativities in the transport solution. 

In this paper, we leverage the close connection between SMM and VEF -- both through linearization and by algebraically reformulating the closures -- to convert the Discontinuous Galerkin (DG) and mixed finite element independent VEF methods presented in \citet{dgvef_olivier} and \citet{rtvef_olivier}, respectively, to the corresponding independent SMMs. 
The SMM moment discretizations are coupled to a high-order DG discretization of the Discrete Ordinates (\Sn) transport equations to solve problems from linear transport. 
Our motivation for this work is in the context of high energy density physics simulations where the tightly coupled simulation of hydrodynamics and thermal radiative transfer is required, the latter of which typically employs the \Sn angular model. 
We are interested in designing discretizations for the SMM moment system that can be scalably solved with existing preconditioning technology and are compatible with the high-order finite element, Lagrangian hydrodynamics framework being developed at the Lawrence Livermore National Laboratory (LLNL) \cite{blast} where the thermodynamic variables are represented using discontinuous finite elements. 
Both the DG and mixed finite element approaches satisfy the compatibility requirement by producing the SMM scalar flux solution in a DG finite element space. 
The mixed finite element discretization has the added benefit that the SMM diffusion operator exactly matches the mixed finite element techniques used for radiation diffusion calculations at LLNL \cite{pete,doi:10.1137/17M1132562}. 
In such case, the existing diffusion package could be elevated to a transport algorithm by simply supplying a modified, transport-dependent source term to the linear and nonlinear solvers already in place for diffusion, thereby easing the implementational burden associated with multiphysics coupling and maintaining disparate code bases for diffusion and transport. 
We note that mixed finite elements are also used for radiation diffusion-based reactor physics calculations \cite{mfem_diffusion,doi:10.13182/NSE97-A28593,doi:10.13182/NSE07-A2660} for which the mixed finite element SMM could also serve as a drop-in replacement. 
We also present a continuous finite element (CG) discretization for the SMM moment system as a CG-based VEF algorithm was found to perform as well as DG-based algorithms while solving for fewer moment unknowns and requiring simpler preconditioning strategies \cite{dgvef_olivier}. 

Since the left hand side of the SMM moment system is equivalent to radiation diffusion, consistent SMMs can be designed by using any of the diffusion discretizations developed for consistent Diffusion Synthetic Acceleration (DSA) \cite{A} methods such as \citet{WWM}, \citet{AM}, \citet{WR}, and \citet{ldrd_dsa}. 
Furthermore, many consistent DSA schemes can be scalably solved with preconditioned iterative solvers (e.g.~\cite{WR,ldrd_dsa}). 
The design of efficient, consistent SMM algorithms then only requires developing consistent discretizations for the SMM correction sources. 
Thus, in the case of SMM, the consistent approach is likely to be effective in achieving our research goals of multiphysics compatibility and scalable preconditioned iterative solvers for the discretized moment system. 
However, such a method would not enjoy the flexibilities discussed above and, in particular, would preclude the possibility of a SMM algorithm acting as a drop-in replacement for an existing diffusion package. 

Previous work on discrete SMMs is limited in the literature. 
\citet{STEHLE2014325} developed a domain decomposition method where SMM was used to couple transport and diffusion domains through interface conditions. 
The SMM moment system and correction sources were discretized to be algebraically consistent with the simple corner balance \cite{ucb_maginot} \Sn transport method. 
\citet{https://doi.org/10.48550/arxiv.2102.09054} investigated a multilevel-in-energy algorithm for neutron transport problems. 
The discretization of the moment system was designed to be consistent with a lowest-order DG discretization of the \Sn transport equations. 
This method can be viewed as a linearization of the DG VEF method from \citet{dima_dfem}. 
Neither of these methods attain higher than second-order accuracy or are compatible with curved meshes. 
Furthermore, independent SMMs have not yet been investigated. 

We proceed as follows. 
First, we introduce a general framework for the iterative moment algorithm in the context of the mono-energetic, fixed-source transport problem with isotropic scattering. 
We illustrate the iteration in a general context by abstracting away the choice of closure (i.e.~VEF vs.~SMM) as well as casting the moment system in first or second-order form. 
We also demonstrate that the SMM and VEF algorithms are related through linearization using the Gateaux derivative. 
Next, we provide background on the relevant finite element technology used to discretize the transport and SMM moment equations on curved meshes and discuss the discretization overlap that occurs in computing the SMM correction sources and the transport equation's scattering source. 
We derive DG, CG, Raviart Thomas mixed finite element, and hybridized Raviart Thomas mixed finite element SMM moment system discretizations through their close connection to the corresponding discrete VEF methods derived in \cite{dgvef_olivier,rtvef_olivier}.
We then present numerical results. 
We show that the SMMs converge the scalar flux with optimal accuracy on a manufactured transport problem and preserve the diffusion limit on an orthogonal and a curved mesh. 
We also compare the performance of the methods on a multi-material problem and show parallel scaling results. 
Finally, we give conclusions and recommendations for future work. 

\section{Moment Methods for Radiation Transport}
We take the mono-energetic, fixed-source transport problem with isotropic scattering given by: 
	\begin{subequations}
	\begin{equation} \label{eq:transport}
		\Omegahat\cdot\nabla\psi + \sigma_t \psi = \frac{\sigma_s}{4\pi}\int \psi \ud \Omega' + q \,, \quad \x \in \D \,, 
	\end{equation}
	\begin{equation} \label{eq:inflow}
		\psi(\x,\Omegahat) = \bar{\psi}(\x,\Omegahat) \,, \quad \x \in \partial \D \ \mathrm{and} \ \Omegahat\cdot\n < 0 \,,
	\end{equation}
	\end{subequations}
as our model problem. 
Here, $\x \in \R^{\dim}$ and $\Omegahat\in \mathbb{S}^2$ are the spatial and angular variables, respectively, $\psi(\x,\Omegahat)$ the angular flux, $\D\subset\R^{\dim}$ the spatial domain of the problem with $\partial\D$ its boundary, $\sigma_t(\x)$ and $\sigma_s(\x)$ the total and scattering macroscopic cross sections, respectively, $q(\x,\Omegahat)$ the fixed-source, and $\bar{\psi}(\x,\Omegahat)$ the inflow boundary function. 
In this section, we derive the moment system and its boundary conditions. 
Due to the streaming term, $\Omegahat\cdot\nabla\psi$, angular moments always produce more unknowns than equations. 
We define two types of closures that give rise to the VEF and SMM moment systems. 
We then describe the iterative schemes used to solve the coupled transport-moment systems. 
The section concludes with a discussion of the close connection between the VEF and SMM algorithms established by \citet{cefus}. 

\subsection{Derivation of Moment System}
The moment system is comprised of the zeroth and first angular moments of the transport equation: 
	\begin{subequations}
	\begin{equation}
		\nabla\cdot\vec{J} + \sigma_a \varphi = Q_0 \,,  
	\end{equation}
	\begin{equation}
		\nabla\cdot\P + \sigma_t \vec{J} = \vec{Q}_1 \,, 
	\end{equation}
	\end{subequations}
where $\sigma_a(\x) \equiv \sigma_t(\x) - \sigma_s(\x)$ is the absorption macroscopic cross section, $Q_0 = \int q \ud \Omega$ and $\vec{Q}_1 = \int \Omegahat\, q \ud \Omega$ the zeroth and first moments of the fixed-source, $q$, and $\varphi$, $\vec{J}$, and $\P$ the zeroth, first, and second angular moments of the angular flux, respectively. 
We refer to the first three moments of the angular flux as the scalar flux, current, and pressure, respectively. 

Boundary conditions are derived by manipulating partial currents. 
Let $J_n^\pm = \int_{\Omegahat\cdot\n \gtrless 0} \Omegahat\cdot\n \, \psi \ud \Omega$ denote the partial currents where $\n$ is the outward unit normal to the boundary of the domain. 
The net current can be expressed in terms of the partial currents and is manipulated to yield: 
	\begin{equation}
	\begin{aligned}
		\vec{J}\cdot\n &= J_n^+ + J_n^- \\ 
		&= 2J_n^- + (J_n^+ - J_n^-) \\
		&= 2J_n^- + \int |\Omegahat\cdot\n|\, \psi \ud \Omega \,. 
	\end{aligned}
	\end{equation}
Defining 
	\begin{equation}
		B(\psi) = \int |\Omegahat\cdot\n|\, \psi \ud \Omega\,, 
	\end{equation}
the boundary conditions for the moment system are: 
	\begin{equation}
		\vec{J}\cdot\n = B + 2\Jin \,, 
	\end{equation}
where $\Jin = \int_{\Omegahat\cdot\n<0} \Omegahat\cdot\n\, \bar{\psi} \ud \Omega$ is the incoming partial current computed using the inflow boundary function, $\bar{\psi}$. 

Combined, the moment system is given by: 
	\begin{subequations} \label{eq:moment_unclosed}
	\begin{equation}
		\nabla\cdot\vec{J} + \sigma_a \varphi = Q_0 \,,  \quad \x \in \D\,,
	\end{equation}
	\begin{equation}
		\nabla\cdot\P + \sigma_t \vec{J} = \vec{Q}_1 \,, \quad \x \in \D\,,
	\end{equation}
	\begin{equation}
		\vec{J}\cdot\n = B + 2\Jin \,, \quad \x \in \partial \D \,. 
	\end{equation}
	\end{subequations}
In three dimensions, there are 10 unknowns corresponding to the scalar flux, the three components of the current, and the six unique components of the symmetric pressure tensor but only four equations arising from the scalar zeroth moment and vector first moment equations. 
In addition, we have only one equation on the boundary but two unknowns corresponding to the normal component of the current and the boundary functional, $B$. 
Thus, we need to provide additional symmetric tensor and scalar-valued equations on the interior and boundary of the domain, respectively, to form a solvable system of equations. 
These additional equations are referred to as closures. 
For the moment methods considered in this document, the closures are exact such that the moment system is an equivalent reformulation of the transport equation and trivial in that the transport solution must already be known to define the closures. 
\begin{rem} \label{rem:independent_phis}
When the discretized moment and transport equations are not algebraically consistent (e.g.~in the context of an independent moment method), the moments of the angular flux and the solution of the moment system will not be equivalent; they will differ on the order of the spatial and angular discretization errors associated with the discrete transport and moment systems. 
To notationally separate the two scalar flux solutions, we use $\varphi$ (varphi) to denote the moment system's scalar flux and $\phi = \int \psi \ud \Omega$ (phi) for the zeroth moment of the angular flux.
Since the two solutions differ on the order of the discretization error, we can approximate $\varphi/\phi \approx 1$ and $\varphi - \phi \approx 0$ without altering the overall transport algorithm's order of accuracy in space or angle. 
\end{rem}

\subsection{Variable Eddington Factor and Second Moment Method Closures}
We consider two types of closures for the pressure and boundary functional which give rise to the VEF and SMM algorithms. 
For VEF, the pressure and boundary functional are closed by multiplying and dividing by the scalar flux: 
	\begin{equation}
		\P = \E\varphi \,,
	\end{equation}
	\begin{equation}
		B = E_b \varphi \,,
	\end{equation}
where 
	\begin{subequations}
	\begin{equation}
		\E = \frac{\int \Omegahat\otimes\Omegahat\, \psi \ud \Omega}{\int \psi \ud \Omega} \,,
	\end{equation}
	\begin{equation}
		E_b = \frac{\int |\Omegahat\cdot\n|\, \psi \ud \Omega}{\int \psi \ud \Omega} \,,
	\end{equation}
	\end{subequations}
are the Eddington tensor and boundary factor, respectively. 
Since, for the continuous equations, $\varphi / \int \psi \ud\Omega = 1$, these closures are simply algebraic manipulations of the pressure and boundary functional. 
When the transport and moment systems are discretized independently, the closure procedure induces errors on the order of the discretization error that can safely be ignored without altering the discrete transport algorithm's spatial or angular order of accuracy (see Remark \ref{rem:independent_phis}). 
Note that the VEF closures are bounded, nonlinear functions of the transport solution due to their normalization by the scalar flux. 
Furthermore, this normalization means the VEF closures do not depend on the magnitude of the angular flux, only its angular shape. 
If $\psi = \frac{1}{4\pi}(f + \Omegahat\cdot\vec{g})$ is a linearly anisotropic function for some spatially dependent functions $f(\x)$ and $\vec{g}(\x)$, then $\E = 1/3\I$ and $E_b = 1/2$. 

On the other hand, SMM uses additive closures of the form: 
	\begin{subequations}
	\begin{equation} \label{eq:smm_Pclosure}
		\P = \T + \frac{1}{3}\I \varphi \,,  
	\end{equation}
	\begin{equation} \label{eq:smm_Bclosure}
		B = \beta + \frac{1}{2}\varphi \,,
	\end{equation}
	\end{subequations}
where 
	\begin{subequations}
	\begin{equation} \label{eq:correction_tensor}
		\T = \int \Omegahat\otimes\Omegahat\, \psi \ud \Omega - \frac{1}{3}\I \int \psi \ud \Omega \,, 
	\end{equation}
	\begin{equation} \label{eq:bdr_correction_factor}
		\beta = \int |\Omegahat\cdot\n|\,\psi \ud \Omega - \frac{1}{2} \int \psi \ud \Omega \,,
	\end{equation}
	\end{subequations}
are the SMM correction tensor and boundary factor, respectively. 
Here, $\I$ denotes the $\dim \times \dim$ identity tensor. 
In the same way that the VEF closures are derived by multiplying and dividing by the scalar flux, the SMM closures are formulated by adding and subtracting the scalar flux. 
Analogously to VEF, this closure procedure is simply an algebraic reformulation of the pressure and boundary functional for the continuous equations. 
When discretized independently, $\varphi - \int\psi \ud \Omega$ is on the order of the discretization error and, by the same arguments used for VEF, this numerical discrepancy can be ignored. 
Due to their additive formulation, the SMM closures are linear functions of the transport solution that depend on both the magnitude and angular shape of the transport solution. 
The SMM closures $\T$ and $\beta$ vanish when the transport solution is linearly anisotropic. 
The factors of $1/3$ and $1/2$ used in the correction tensor and boundary factor, respectively, cause the SMM moment system to match the radiation diffusion equation with Marshak boundary conditions when the transport solution is linearly anisotropic. 

Using these closures in Eqs.~\ref{eq:moment_unclosed}, the VEF and SMM moment systems are given by: 
	\begin{subequations} \label{eq:vef_first}
	\begin{equation}
		\nabla\cdot\vec{J} + \sigma_a \varphi = Q_0 \,, \quad \x \in \D \,, 
	\end{equation}
	\begin{equation}
		\nabla\cdot\paren{\E\varphi} + \sigma_t \vec{J} = \vec{Q}_1 \,, \quad \x \in \D \,,
	\end{equation}
	\begin{equation} \label{eq:vef_bc}
		\vec{J}\cdot\n = E_b \varphi + 2\Jin \,, \quad \x \in \partial \D \,, 
	\end{equation}
	\end{subequations}
and 
	\begin{subequations} \label{eq:smm_first}
	\begin{equation}
		\nabla\cdot\vec{J} + \sigma_a \varphi = Q_0 \,, \quad \x \in \D \,, 
	\end{equation}
	\begin{equation}
		\frac{1}{3}\nabla\varphi + \sigma_t \vec{J} = \vec{Q}_1 - \nabla\cdot \T \,, \quad \x \in \D \,,
	\end{equation}
	\begin{equation} \label{eq:smm_bc}
		\vec{J}\cdot\n = \frac{1}{2} \varphi + 2\Jin + \beta \,, \quad \x \in \partial \D \,, 
	\end{equation}
	\end{subequations}
respectively. 
By eliminating the first moment equation, the moment system can be recast as a second-order partial differential equation. 
For VEF, the second-order form is: 
	\begin{equation} \label{eq:vef_second}
		-\nabla\cdot \frac{1}{\sigma_t}\nabla\cdot\paren{\E\varphi} + \sigma_a \varphi = Q_0 - \nabla\cdot \frac{\vec{Q}_1}{\sigma_t} \,,
	\end{equation}
while for SMM we have: 
	\begin{equation} \label{eq:smm_second}
		-\nabla\cdot \frac{1}{3\sigma_t} \nabla\varphi + \sigma_a \varphi = Q_0 - \nabla\cdot\frac{\vec{Q}_1}{\sigma_t} + \nabla\cdot \frac{1}{\sigma_t}\nabla\cdot\T \,. 
	\end{equation}
Here, boundary conditions are specified by Eqs.~\ref{eq:vef_bc} and \ref{eq:smm_bc} for the VEF and SMM second-order forms, respectively. 
In both cases, the moment system is equivalent to radiation diffusion with Marshak boundary conditions when the transport solution is linearly anisotropic. 
Note that the left hand side of the VEF operator is non-self-adjoint due to the presence of the Eddington tensor inside the divergence while the SMM left hand side is simply the self-adjoint radiation diffusion approximation. 
Observe that the SMM correction source, $\nabla\cdot \frac{1}{\sigma_t}\nabla\cdot\T$, shares the same mathematical structure as the differential term in the VEF operator, $\nabla\cdot\frac{1}{\sigma_t}\nabla\cdot\paren{\E\varphi}$, and thus will have analogous discretization challenges. 
In this document, we leverage recently developed discretization techniques for the VEF moment system to derive discretizations for the SMM moment system. 

\subsection{The Moment Method Algorithm}
Moment methods solve the coupled transport-moment system simultaneously. 
The transport equation is used to provide the VEF or SMM closures while the moment system is used to compute the moment-dependent physics. 
In the present discussion, the moment system's scalar flux solution is used to compute the isotropic scattering source. 
In this way, the coupling of the angular phase space induced by integrating over angle is avoided, allowing use of the transport sweep to efficiently invert the transport equation's streaming and collision operator. 
Simple iterative schemes often converge rapidly and robustly due to the closures' weak dependence on the transport solution. 

We first introduce notation that abstracts away the choice of the closures and casting the moment system in first or second-order form. 
Let $\mathcal{M}(\psi,\mat{X}) = 0$ denote one of the moment systems derived in the previous subsection with $\mat{X}$ the moment system's unknowns. 
For example, $\mathcal{M}(\psi,\mat{X})$ could represent the VEF moment system in first-order form given by Eqs.~\ref{eq:vef_first} where $\mat{X}$ would include both the scalar flux and current. 
In the case of the second-order form, we would set $\mat{X} = \varphi$ since the scalar flux is the only unknown. 
For VEF, $\mathcal{M}(\psi,\mat{X})$ is nonlinear in $\psi$ and linear in the moments, $\mat{X}$. 
For SMM, $\mathcal{M}(\psi,\mat{X})$ is linear in both arguments. 

The moment algorithm solves the coupled system given by: 
	\begin{subequations} \label{eq:moment_coupled}
	\begin{equation}
		\Omegahat\cdot\nabla\psi + \sigma_t\psi = \frac{\sigma_s}{4\pi}\varphi + q \,,
	\end{equation}
	\begin{equation}
		\mathcal{M}(\psi,\mat{X}) = 0 \,,
	\end{equation}
	\end{subequations}
where transport boundary conditions are specified in Eq.~\ref{eq:inflow}. 
The moment system's boundary conditions are given by Eq.~\ref{eq:vef_bc} for a VEF method and Eq.~\ref{eq:smm_bc} for SMM. 
Here, the moment system is coupled to the transport equation through the closures and the transport equation's scattering source is coupled to the moment system through the moment system's scalar flux. 
We have increased the complexity of the problem by adding the moment system's unknowns. 
In the case of VEF, the coupled system in Eq.~\ref{eq:moment_coupled} is also nonlinear due to the use of nonlinear closures. 
However, solving the coupled system is still advantageous due to the ability to use the transport sweep and the rapid convergence of the closures.  

Let 
	\begin{equation} 
		\mat{L}\psi = \Omegahat\cdot\nabla\psi + \sigma_t\psi 
	\end{equation}
represent the streaming and collision operator. 
The coupled transport-moment system can then be rewritten 
	\begin{subequations}
	\begin{equation}
		\mat{L}\psi = \frac{\sigma_s}{4\pi}\varphi + q \,, 
	\end{equation}
	\begin{equation}
		\mathcal{M}(\psi,\mat{X}) = 0 \,. 
	\end{equation}
	\end{subequations}
By linearly eliminating the angular flux, the coupled system is equivalent to: 
	\begin{equation} \label{eq:lo_only}
		\mathcal{M}\!\paren{\mat{L}^{-1}\!\paren{\frac{\sigma_s}{4\pi}\varphi +q},\mat{X}} = 0 \,. 
	\end{equation}
Observe that Eq.~\ref{eq:lo_only} is now a function of the moment solution only. 
That is, we can define 
	\begin{equation} 
		\F(\mat{X}) = \mathcal{M}\!\paren{\mat{L}^{-1}\!\paren{\frac{\sigma_s}{4\pi}\varphi +q},\mat{X}} 
	\end{equation}
and equivalently solve $\F(\mat{X}) = 0$. 
In this reduced problem, the angular flux appears only as an auxiliary variable used to compute the residual $\mat{F}(\mat{X})$ and we say that the angular flux is determined by the moment system. 
This reduced formulation $\mat{F}(\mat{X}) = 0$ has much lower dimension than the original coupled system given in Eq.~\ref{eq:moment_coupled} but has the same solution. 
Due to this, advanced solvers for $\F(\mat{X})$ can be applied that would otherwise be impractical for Eq.~\ref{eq:moment_coupled} due to the storage and computational costs associated with the high-dimensionality of the angular flux. 

We now leverage the structure of the VEF and SMM moment systems to further simplify the above algorithm. 
Let 
	\begin{equation}
		\mat{V}(\psi)\mat{X} = f \,, 
	\end{equation}
represent the VEF moment system such that $\mathcal{M}(\psi,\mat{X}) = \mat{V}(\psi)\mat{X} - f$. 
We then have that 
	\begin{equation}
		\F(\mat{X}) = \mat{V}\!\paren{\mat{L}^{-1}\!\paren{\frac{\sigma_s}{4\pi}\varphi + q}} \mat{X} - f = 0 \,. 
	\end{equation}
Operating by the inverse of the VEF moment system, the coupled transport-VEF system is equivalent to: 
	\begin{equation}
		\mat{X} = \mat{V}\!\paren{\mat{L}^{-1}\!\paren{\frac{\sigma_s}{4\pi}\varphi + q}}^{-1} f \,. 
	\end{equation}
For SMM, the moment system is of the form
	\begin{equation}
		\mat{D}\mat{X} = b(\psi) \,,
	\end{equation}
where $\mat{D}$ is a diffusion operator and $b(\psi)$ includes the moments of the fixed-source and the transport-dependent correction sources. 
The root-finding problem $\mat{F}(\mat{X}) = 0$ is then equivalent to 
	\begin{equation} 
		\mat{X} = \mat{D}^{-1}b\!\paren{\mat{L}^{-1}\!\paren{\frac{\sigma_s}{4\pi}\varphi + q}} \,. 
	\end{equation}
Thus, for both VEF and SMM, the solution of the coupled transport-moment system is the fixed-point: 
	\begin{equation} \label{moment:G}
		\mat{X} = \mat{G}(\mat{X}) \,,
	\end{equation}
where $\mat{G}(\mat{X})$ is given by 
	\begin{equation} \label{moment:vef_fp}
		\mat{G}(\mat{X}) = \mat{V}\!\paren{\mat{L}^{-1}\!\paren{\frac{\sigma_s}{4\pi}\varphi + q}}^{-1} f
	\end{equation}
for VEF and 
	\begin{equation} \label{moment:smm_fp}
		\mat{G}(\mat{X}) = \mat{D}^{-1}b\!\paren{\mat{L}^{-1}\!\paren{\frac{\sigma_s}{4\pi}\varphi + q}}
	\end{equation}
for SMM.
The fixed-point operator $\mat{G}$ is applied in two stages: 1) apply a transport sweep to invert the streaming and collision operator on a scattering source formed from the moment system's scalar flux and 2) solve the moment system using the closures computed with the angular flux from stage 1). 
The definitions of the fixed-point operator, $\mat{G}$, show the key differences between the VEF and SMM algorithms. 
VEF has a transport-dependent left hand side operator while the right hand side sources are fixed. 
On the other hand, SMM has transport-dependent sources but a fixed left hand side operator corresponding to radiation diffusion. 

\begin{figure}
\centering
\begin{subfigure}{.47\textwidth}
	\centering
	\includegraphics[width=\textwidth]{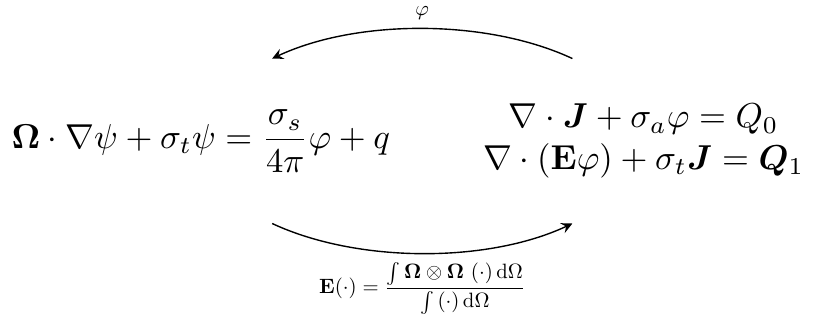}	
	\caption{}
\end{subfigure}
\hfill
\begin{subfigure}{.47\textwidth}
	\centering
	\includegraphics[width=\textwidth]{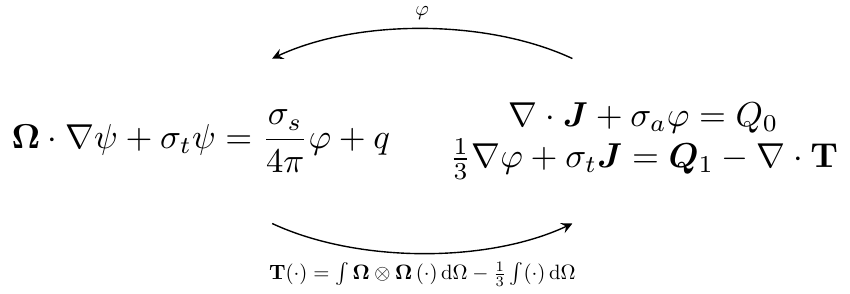}	
	\caption{}
\end{subfigure}
\caption{A depiction of the iterative coupling used in the (a) VEF and (b) SMM algorithms. The transport equation informs the moment system through the closures while the moment system drives the transport equation by computing the scattering physics. By lagging the scattering term, the transport equation can be efficiently inverted. Rapid convergence occurs because the closures are weak functions of the solution. We omit boundary conditions for clarity of presentation.}
\label{fig:moment_alg}
\end{figure}

The simplest algorithm to solve $\mat{X} = \mat{G}(\mat{X})$ is fixed-point iteration: 
	\begin{equation} \label{eq:fixed-point}
		\mat{X}^{k+1} = \mat{G}(\mat{X}^k) 
	\end{equation}
where $\mat{X}^0$ is an initial guess. 
This process is repeated until the difference between successive iterates is small enough. 
We will also use Anderson acceleration \cite{10.1145/321296.321305}, a more advanced algorithm for solving fixed-point problems that stores a set of $k$ previous evaluations of the fixed-point operator and chooses the next iterate by finding the optimal linear combination of the $k$ operator evaluations that minimizes the residual, $\mat{F}(\mat{X})$. 
For a storage cost of $k$ moment system-sized vectors, Anderson acceleration increases iterative efficiency and robustness. 
We stress that the fixed-point problem in Eq.~\ref{eq:fixed-point} corresponds to the reduced formulation where the angular flux is determined by the moment system, avoiding the need to store $k$ angular flux-sized vectors. 
The coupling between the transport and moment equations for the VEF and SMM algorithms is depicted in Fig.~\ref{fig:moment_alg}. 
Convergence of the fixed-point iteration is expected to be rapid since the closures are weak functions of the transport solution \cite{goldin}. 

\begin{rem} \label{rem:smm_advantage}
The changing left hand side operator in the VEF algorithm means that costs associated with forming a new sparse matrix, such as finite element assembly and solver setup (e.g.~LU factorization or the AMG setup phase), must be incurred at each iteration. 
Since SMM's left hand side is iteratively fixed, setup costs need only be incurred once, amortizing them over the entire iteration. 
Furthermore, less computational work is needed to assemble a right hand side vector than is needed to assemble a left hand side sparse matrix. 
Finally, SMM only requires solving the self-adjoint diffusion operator at each iteration whereas VEF requires the solution of a generally non-symmetric operator. 
In this way, simpler iterative solvers can be applied to the SMM moment system that would not be convergent for the VEF moment system. 
For example, we will use conjugate gradient to solve three of the proposed SMM discretizations while the analogous VEF methods require the use of non-symmetric solvers such as BiCGStab or GMRES. 
Thus, it is expected that SMM's fixed-point operator will be less expensive to evaluate than VEF's. 
\end{rem}

\subsection{Connection via Linearization}
In the process of performing a Fourier stability analysis of the VEF algorithm, \citet{cefus} showed that SMM is equivalent to the VEF algorithm linearized about a linearly anisotropic solution. 
Here, we use this connection between VEF and SMM to provide an alternate method for deriving the SMM closures. 
We stress that the pressure and boundary functional are both \emph{linear} functions of the angular flux. 
Thus, the following linearizations of the pressure and boundary functional do not possess linearization error and are instead simply algebraic reformulations. 

We first derive the functional derivatives of the Eddington tensor and boundary factor. 
We use the Gateaux derivative \cite{vainberg64}, a generalization of the directional derivative that supports taking derivatives of functionals. 
The Gateaux derivative of a functional $f$ evaluated at some function $u$ in the direction of $v$ is defined as: 
	\begin{equation}
		D[f](u,v) = \lim_{\omega\rightarrow 0} \frac{f(u + \omega v) - f(u)}{\omega} \,. 
	\end{equation}
Assuming continuity of $f$, the Gateaux derivative is equivalent to 
	\begin{equation} \label{eq:gateaux}
		D[f](u,v) = \bracket{\pderiv{}{\omega}f(u+\omega v)}_{\omega = 0} \,. 
	\end{equation}
This second definition is preferred as it leads to simpler calculations by leveraging the familiar machinery of the partial derivative. 
For multi-variate arguments, such as when $\mat{u} = \vector{u_1 & \ldots & u_n}$ and $\mat{v} = \vector{v_1 & \ldots & v_n}$ for some $n>1$, the Gateaux derivative is applied in a manner analogous to the partial derivative where differentiation is applied to one variable at a time with the rest fixed. 
Using the Gateaux derivative, the first-order Taylor series of $f(u)$ expanded about $u_0$ is: 
	\begin{equation}
		f(u) \xrightarrow{\text{TSE}} f(u_0) + D[f](u_0\,,\, u - u_0) \,. 
	\end{equation}
In this way, the Gateaux derivative can be used to apply the action of the Jacobian in the linearization process. 

Here, we seek to expand the pressure and boundary functional about a linearly anisotropic solution using the Gateaux derivative. 
Let $\psi_0$ be such a linearly anisotropic approximation to the transport problem at hand. 
Using Eq.~\ref{eq:gateaux}, the Gateaux derivative of the Eddington tensor at $\psi_0$ in the direction $\psi$ is: 
	\begin{equation}
	\begin{aligned}
		D[\E](\psi_0,\psi) &= \bracket{\pderiv{}{\omega}\E(\psi_0 + \omega \psi)}_{\omega = 0} \\
		&= \bracket{\pderiv{}{\omega} \frac{\P_0 + \omega \P}{\phi_0 + \omega \phi}}_{\omega = 0} \\ 
		&= \frac{\P(\phi_0 + \omega \phi) - (\P_0 + \omega \P) \phi}{(\phi_0 + \omega \phi)^2}\biggr\rvert_{\omega = 0} \\
		&= \frac{\P \phi_0 - \P_0 \phi}{\phi_0^2} \\
		&= \frac{1}{\phi_0}\paren{\P - \E_0 \phi} \,,
	\end{aligned}
	\end{equation}
where we have used the quotient rule and have set $\P_0 = \int \Omegahat\otimes \Omegahat \, \psi_0 \ud \Omega$, $\phi_0 = \int \psi_0 \ud \Omega$, and $\E_0 = \P_0/\phi_0$. 
Since $\psi_0$ is linearly anisotropic, $\E_0 = 1/3 \I$. 
Thus, the Gateaux derivative of the Eddington tensor is equivalent to 
	\begin{equation} \label{eq:closure_connection}
		D[\E](\psi_0,\psi) = \frac{1}{\phi_0}\paren{\int \Omegahat\otimes \Omegahat\, \psi \ud \Omega - \frac{1}{3}\I \int \psi \ud \Omega} \equiv \frac{1}{\phi_0}\T(\psi) \,, 
	\end{equation}
where $\T(\psi)$ is the SMM correction tensor defined in Eq.~\ref{eq:correction_tensor}. 
By an analogous calculation, the Gateaux derivative of the Eddington boundary factor is: 
	\begin{equation}
		D[E_b](\psi_0,\psi) = \frac{1}{\phi_0}\paren{\int |\Omegahat\cdot\n|\,\psi \ud \Omega - E_{b0}\int \psi \ud \Omega}\,, 
	\end{equation}
where 
	\begin{equation} \label{eq:eb0}
		E_{b0} = \frac{\int |\Omegahat\cdot\n|\,\psi_0 \ud \Omega}{\int \psi_0 \ud \Omega} = \frac{\int |\Omegahat\cdot\n| \ud \Omega}{4\pi} = \frac{1}{2} \,, 
	\end{equation}
given that $\psi_0$ is linearly anisotropic. 
Thus, the Gateaux derivative of the Eddington boundary factor and the SMM boundary correction factor are related by: 
	\begin{equation} \label{eq:bdr_closure_connection}
		D[E_b](\psi_0,\psi) = \frac{1}{\phi_0}\beta(\psi) \,, 
	\end{equation}
where $\beta(\psi)$ is the SMM boundary factor defined in Eq.~\ref{eq:bdr_correction_factor}. 

We now linearize the pressure and boundary functional using the above functional derivatives of the VEF closures. 
Let $\y = \vector{\psi & \varphi}$ represent the solution of the coupled transport-moment system where, without loss of generality, we consider the moment system to be written in second-order form. 
Further, we set $\y_0 = \vector{\psi_0 & \varphi_0}$ to be the linearly anisotropic approximation of the coupled system. 
Here, we consider the moment solution, $\varphi$, to be an independent variable. 
Using the Gateaux derivative, we can write the first-order Taylor series expansion of the pressure as: 
	\begin{equation}
	\begin{aligned}
		\P(\y) &\xrightarrow{\text{TSE}} \bracket{\E\varphi}_{\y=\y_0} + D[\E\varphi](\y_0,\y - \y_0) \\
		&= \E_0 \varphi_0 + D[\E\varphi_0](\psi_0,\psi-\psi_0) + D[\E_0\varphi](\varphi_0,\varphi-\varphi_0) \\ 
		&= \E_0 \varphi_0 + \frac{\varphi_0}{\phi_0} \T(\psi - \psi_0) + \E_0 (\varphi - \varphi_0) \\
		&= \E_0 \varphi + \frac{\varphi_0}{\phi_0}\T(\psi - \psi_0) \\
		&= \E_0 \varphi + \T(\psi) \,, 
	\end{aligned}
	\end{equation}
where $\T(\psi - \psi_0) = \T(\psi) - \T(\psi_0) = \T(\psi)$ since the SMM correction tensor is linear and zero when the argument is linearly anisotropic, respectively. 
Further, we have simplified $\varphi_0/\phi_0 \approx 1$ by ignoring terms on the order of the discretization error (see Remark \ref{rem:independent_phis}). 
Using $\E_0 = 1/3\I$, observe that the extreme left and right hand sides of the linearized pressure exactly match the SMM closure given in Eq.~\ref{eq:smm_Pclosure}. 

Analogously, the first-order Taylor series expansion of the boundary functional is: 
	\begin{equation}
	\begin{aligned}
		B(\y) &\xrightarrow{\text{TSE}} \bracket{E_b \varphi}_{\y=\y_0} + D[E_b \varphi](\y_0, \y-\y_0) \\
		&= E_{b0} \varphi + \frac{\varphi_0}{\phi_0}\beta(\psi - \psi_0) \\
		&= E_{b0} \varphi + \beta(\psi) \,, 
	\end{aligned}
	\end{equation}
where we have applied the same arguments to simplify $\beta(\psi - \psi_0)$ and $\varphi_0/\phi_0$ as was used for the pressure. 
Recognizing that $E_{b0} = 1/2$, the linearized boundary functional is equivalent to the SMM closure in Eq.~\ref{eq:smm_Bclosure}. 
\begin{rem} \label{rem:disc_Eb}
While angular quadrature rules can exactly evaluate $E_0 = 1/3 \I$, it is not possible to exactly integrate $E_{b0}$ with numerical quadrature due to the the numerator's integrand, $|\Omegahat\cdot\n|\,\psi_0$, having a discontinuous derivative. 
Generally, $\beta$, $\Jin$, and $E_{b0}$ must all be computed in the same manner, either analytically or by angular quadrature, to ensure the discrete boundary conditions reduce to the Marshak case when the transport solution is linearly anisotropic. 
\end{rem}

In the VEF algorithm, the Eddington tensor and boundary factor are the only sources of nonlinearity. 
Thus, the linearized algorithm is equivalent to the original VEF algorithm with $\P = \E\varphi$ and $B = E_b \varphi$ replaced with $\T + 1/3\I\varphi$ and $\beta + 1/2 \varphi$, respectively. 
Such an algorithm is equivalent to SMM. 
Due to this, SMM can be viewed as a VEF algorithm that has been linearized about a linearly anisotropic solution.
\begin{rem} \label{rem:linearization}
Although SMM can be viewed as a VEF algorithm that has been linearized about a linearly anisotropic solution, we stress that the linearization is actually exact since the pressure and boundary functional are linear functions of the transport solution. 
In other words, SMM is still simply an algebraic reformulation of the transport equation and is thus accurate even in transport regimes. 
\end{rem}

\section{Finite Element Preliminaries}
In this section, we provide background on the finite element technology used to discretize the transport and SMM moment systems. 
We describe the machinery used to represent high-order meshes and the related transformations used to integrate the bilinear and linear forms that comprise the discretizations. 
Finally, we define the finite element spaces used to approximate the angular flux, scalar flux, and current in subsequent sections. 

\subsection{Description of the Mesh}
\begin{figure}
\centering
\begin{subfigure}{.49\textwidth}
	\centering
	\includegraphics[height=1.3in]{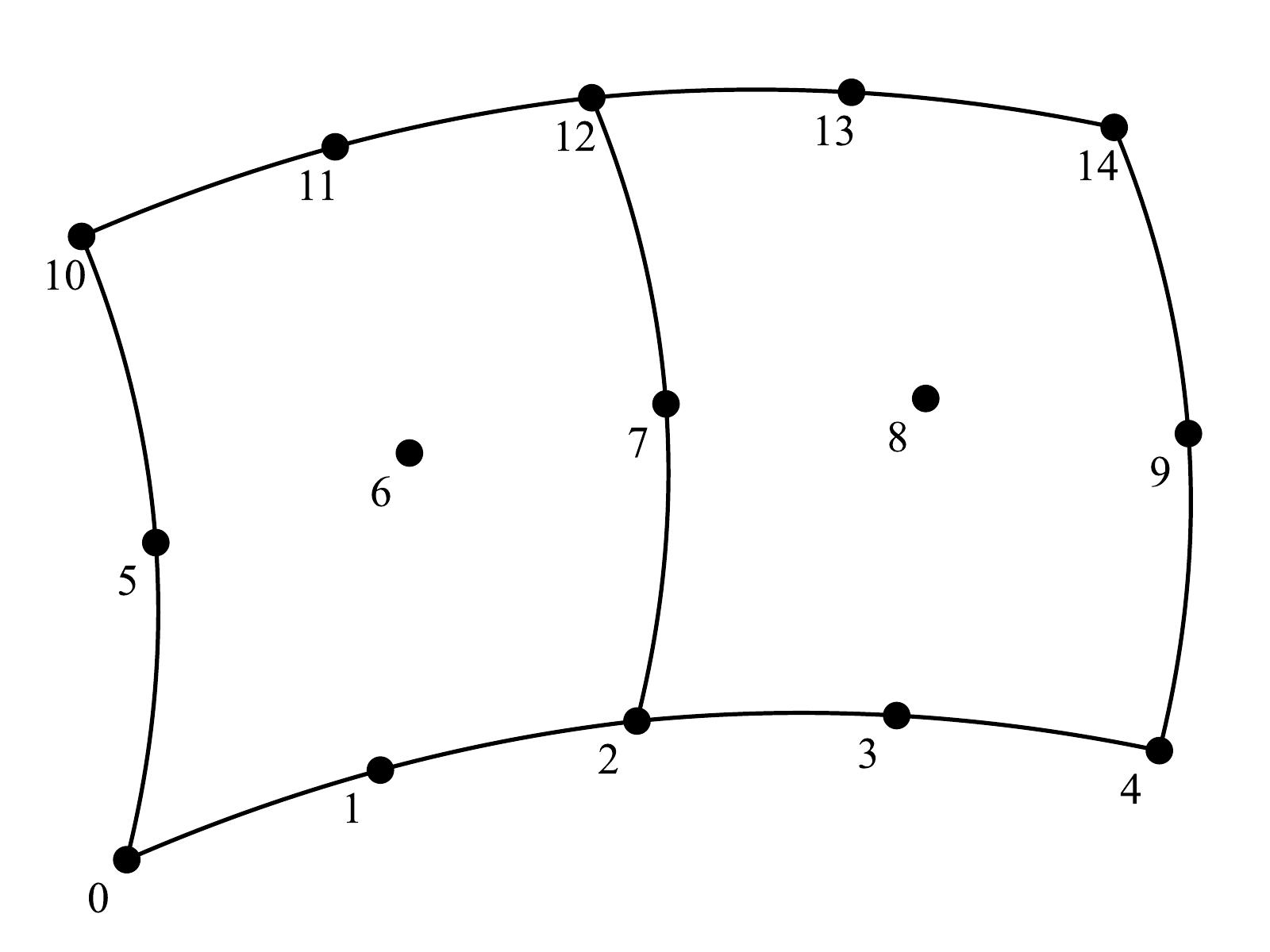}
	\caption{}
	\label{fig:curved_mesh}
\end{subfigure}
\begin{subfigure}{.49\textwidth}
	\centering
	\includegraphics[height=1.3in]{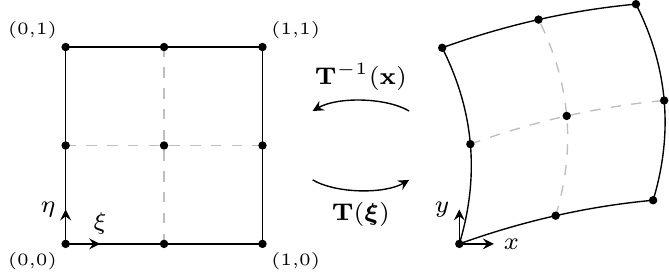}
	\caption{}
	\label{fig:eltrans}
\end{subfigure}
\caption{Depictions of (a) a quadratic, quadrilateral mesh where the mesh control points have been labeled and (b) the element transformation used to describe the left element of the mesh in (a). }
\end{figure}

Let $\D \subset \R^2$ be the domain of the problem tessellated into a collection of elements $\meshT$ such that 
	\begin{equation}
		\D = \bigcup_{K \in \meshT} K \,,	
	\end{equation}
where $K \in \meshT$ is a quadrilateral element in the mesh obtained as the image of the reference square $\hat{K} = [0,1]^2$ under an invertible, polynomial mapping $\T : \hat{K} \rightarrow K$. 
Here, the element transformation belongs to $[\Q_m(\hat{K})]^2$ where $\Q_m(\hat{K})$ is the space of polynomials of degree at most $m$ in each variable. 
We use a nodal basis to describe the element transformation. 
Let $\xi_i$ denote the $m+1$ Gauss-Lobatto points in the interval $[0,1]$. The $(m+1)^2$ points $\vec{\xi}_i$ on the unit square $\hat{K} = [0,1]^2$ are given by the two-fold Cartesian product of the one-dimensional points. 
Further, let $\ell_i$ denote the Lagrange interpolating polynomial that satisfies $\ell_i(\vec{\xi}_j) = \delta_{ij}$ where $\delta_{ij}$ is the Kronecker delta. 
The set of functions $\{\ell_i\}_{i=1}^{(m+1)^2}$ form a basis for $\Q_m(\hat{K})$. 
For each element $K \in \meshT$, the mapping is of the form: 
	\begin{equation}
		\x(\vec{\xi}) = \T(\vec{\xi}) = \sum_{i=1}^{(m+1)^2} \x_{K,i} \ell_i(\vec{\xi}) \,, 
	\end{equation}
where $\x \in K$, $\vec{\xi} \in \hat{K}$, and the $\{\x_{K,i}\}_{i=1}^{(m+1)^2}$ are the control points for element $K$. 
A quadratic quadrilateral mesh is depicted in Fig.~\ref{fig:curved_mesh} where, for example, the left element is described by the control points labeled $(0,1,2,5,6,7,10,11,12)$. 
Observe that the control points along the interface between the two elements are shared, ensuring that the mesh is continuous. 
The reference space to physical space transformation for the left element is depicted in Fig.~\ref{fig:eltrans}. 

We define the set of unique faces in the mesh as $\Gamma$ with $\Gamma_0 = \Gamma \setminus \partial \D$ the set of interior faces in the mesh and $\Gamma_b = \Gamma \cap \partial \D$ the set of boundary faces so that $\Gamma = \Gamma_0 \cup \Gamma_b$. 
We will often work with piecewise functions defined element-wise on the mesh. 
On an interior face $\mathcal{F} \in \Gamma_0$ between two arbitrary elements $K_1$ and $K_2$, we use the convention that $\n$ is a unit vector perpendicular to the shared face $K_1 \cap K_2$ pointing from $K_1$ to $K_2$. 
For a piecewise polynomial $u$ described by $u_1$ for $\x \in K_1$ and $u_2$ for $\x \in K_2$, the average, $\avg{\cdot}$, and jump, $\jump{\cdot}$, of $u$ across the face $\mathcal{F} = K_1 \cap K_2$ are defined as 
	\begin{equation}
		\avg{u} = \frac{1}{2}(u_1 + u_2) \,, \quad \jump{u} = u_1 - u_2 \,, \quad \text{on} \ \mathcal{F} = K_1 \cap K_2 \in \Gamma_0 \,, 
	\end{equation}
with an analogous definition for vector-valued piecewise functions. 
On boundary faces, the jump and average are defined as 
	\begin{equation}
		\jump{u} = u \,, \quad \avg{u} = u \,, \quad \text{on} \ \mathcal{F} \in \Gamma_b \,, 
	\end{equation}
and likewise for vectors on the boundary. 
We note that a piecewise-continuous function $u$ satisfies $\jump{u} = 0$ and $\avg{u} = u$ on each interior face. 

Finally, we define the ``broken'' gradient, denoted by $\nablah$, obtained by applying the gradient locally on each element. That is, 
	\begin{equation}
		(\nablah u)|_K = \nabla(u|_K) \,, \quad \forall K \in \meshT \,. 
	\end{equation}
The broken gradient is important for piecewise-discontinuous functions where the global gradient is not well defined due to the discontinuities present at interior mesh interfaces. 

\subsection{Integration Transformations}
The mesh transformations $\T$ are used to facilitate numerical integration on arbitrary elements. 
Let $\vec{\xi} = \vector{\xi & \eta} \in \hat{K}$ denote the reference coordinates and $\x = \vector{x & y} = \T(\vec{\xi})$ the physical coordinates.  
The Jacobian of the transformation is 
	\begin{equation}
	 	\F = \pderiv{\T}{\vec{\xi}} = \begin{bmatrix} 
	 		\pderiv{x}{\xi} & \pderiv{x}{\eta} \\ 
	 		\pderiv{y}{\xi} & \pderiv{y}{\eta} 
	 	\end{bmatrix} \,, 
	\end{equation} 
with $J = |\F|$ its determinant. 
The partial derivatives of the mesh transformation are computed using the basis expansion of the mesh transformations. 
In other words, 
	\begin{equation}
		\F = \sum_{i=1}^{(m+1)^2} \x_{K,i} \otimes \hnabla \ell_i = \sum_{i=1}^{(m+1)^2} \begin{bmatrix} 
			x_{K,i} \pderiv{\ell_i}{\xi} & x_{K,i} \pderiv{\ell_i}{\eta} \\
			y_{K,i} \pderiv{\ell_i}{\xi} & y_{K,i} \pderiv{\ell_i}{\eta} 
		\end{bmatrix} \,,
	\end{equation}
where $\x_{K,i} = \vector{x_{K,i} & y_{K,i}}$ and $\hnabla$ denotes the gradient with respect to the reference coordinates, $\vec{\xi}$. 
In this document, integration over the domain is implicitly computed in reference space using the following sum: 
	\begin{equation} \label{eq:int_sum}
		\int_\D \paren{\cdot} \ud \x = \sum_{K\in\meshT} \int_K \paren{\cdot} \ud \x = \sum_{K\in\meshT} \int_{\hat{K}} \paren{\cdot}\,J\!\ud \vec{\xi} \,. 
	\end{equation}
This provides a systematic way to integrate over arbitrary domains composed of arbitrarily shaped elements as well as the use of numerical quadrature rules defined on the reference element $\hat{K}$. 

We now discuss the transformations used to represent the integrand of Eq.~\ref{eq:int_sum} in reference space. 
For a scalar function $u: \D \rightarrow \R$, denote by $\hat{u} : \hat{K} \rightarrow \R$ its representation in reference space. The functions $u$ and $\hat{u}$ are related by 
	\begin{equation}
		u(\x) = \hat{u}(\T^{-1}(\x)) \,,
	\end{equation}
where $\T^{-1} : K \rightarrow \hat{K}$ denotes the inverse mesh transformation. 
Integration over the physical element is then equivalent to 
	\begin{equation}
		\int_K u \ud \x = \int_{\hat{K}} \hat{u}\, J\!\ud\vec{\xi} \,. 
	\end{equation}
Using the chain rule, the gradient transforms as 
	\begin{equation} \label{eq:scalar_trans_grad}
		\nabla \hat{u} = \begin{bmatrix} 
			\pderiv{\hat{u}}{\xi}\pderiv{\xi}{x} + \pderiv{\hat{u}}{\eta}\pderiv{\eta}{x} \\
			\pderiv{\hat{u}}{\xi}\pderiv{\xi}{y} + \pderiv{\hat{u}}{\eta}\pderiv{\eta}{y} 
		\end{bmatrix}
		= \mat{F}^{-T} \hnabla \hat{u} \,. 
	\end{equation}
In this way, the gradient in physical space can be computed using the Jacobian of the mesh transformation and the gradient in reference space. 

For the vector-valued functions used in the Raviart Thomas space, the contravariant Piola transformation: 
	\begin{equation} \label{eq:piola}
		\vec{v} = \frac{1}{J}\F \hvec{v} \circ \T^{-1} 
	\end{equation}
is used \cite{ciarlet_elasticity,piola_cisc}. 
Here, $\vec{v} : \D \rightarrow \R^2$ and $\hvec{v} : \hat{K} \rightarrow \R^2$ are the representations of a vector in physical and reference space, respectively. 
Members of the Raviart Thomas space have a continuous normal component, a requirement of any $H(\div;\D)$-conforming space \cite{quateroni}.  
The contravariant Piola transformation represents vectors on the so-called tangent space spanned by the columns of the Jacobian matrix, $\F$. 
Using this basis, the components of a contravariant vector represent the normal and tangential components of the vector along a face. 
This facilitates the sharing of degrees of freedom that represent the normal component across interior mesh interfaces, enforcing the normal continuity required by $H(\div;\D)$. 

For a contravariant vector, 
	\begin{equation} \label{eq:piola_ibp1}
		\int_K \nabla u \cdot \vec{v} \ud \x = \int_{\hat{K}} \mat{F}^{-T}\hnabla \hat{u} \cdot \frac{1}{J}\mat{F} \hat{v}\, J\!\ud\vec{\xi} = \int_{\hat{K}} \hnabla \hat{u} \cdot \hvec{v} \ud\vec{\xi} \,. 
	\end{equation}
The gradient transforms as 
	\begin{equation} \label{eq:piola_grad}
		\nabla\vec{v} = \nabla\paren{\frac{1}{J}\mat{F}\hat{v}\circ\mat{T}^{-1}} = \frac{1}{J}\mat{F}\!\paren{\hnabla\hvec{v} - \hat{\mat{B}}}\!\mat{F}^{-1} 
	\end{equation}
where 
	\begin{equation}
		\hat{\mat{B}} = \frac{1}{J}\hnabla\!\paren{J\mat{F}^{-1}}\!\mat{F}\hvec{v} 
	\end{equation}
represents the portion of the transformation of the gradient arising from the derivatives of the Piola transformation. 
This result is derived by direct computation in \citet[Appendix A]{rtvef_olivier} where implementation details are also provided. 
We note that $\hmat{B}$ is related to the Hessian of the mesh transformation and is thus zero when the mesh transformation is affine such that $\T(\vec{\xi}) = \mat{A} \vec{\xi} + \vec{b}$ for some $\mat{A}\in \R^{\dim\times\dim}$ and $\vec{b} \in \R^{\dim}$ independent of $\vec{\xi}$. 
In addition, $\tr(\hmat{B}) = 0$, a result known as the Piola identity \cite{ciarlet_elasticity}. 
Using the Piola identity, the linearity of the trace, and the invariance of the trace under similarity transformations, the divergence of a contravariant vector transforms as 
	\begin{equation} \label{eq:piola_div}
		\nabla\cdot\vec{v} = \tr\paren{\nabla\vec{v}} = \frac{1}{J}\tr\!\paren{\mat{F}\!\paren{\hnabla\hvec{v} - \hat{\mat{B}}}\!\mat{F}^{-1}} = 
		\frac{1}{J}\!\paren{\tr(\hnabla\hvec{v}) - \tr(\hmat{B})} =
		\frac{1}{J}\hnabla\cdot\hvec{v} \,. 
	\end{equation}
Thus, 
	\begin{equation} \label{eq:piola_ibp2}
		\int_K u\, \nabla\cdot\vec{v} \ud \x = \int_{\hat{K}} \hat{u}\, \hnabla\cdot\hvec{v} \ud\vec{\xi} \,. 
	\end{equation}
Combining the results from Eqs.~\ref{eq:piola_ibp1} and \ref{eq:piola_ibp2} yields: 
	\begin{equation}
		\int_{\partial K} u\,\vec{v}\cdot\n \ud s = \int_{\partial\hat{K}} \hat{u}\, \hvec{v}\cdot\hat{\n} \ud \hat{s} \,, 
	\end{equation}
where $\hat{\n}$ is the normal vector in reference space corresponding to the physical space normal $\n$. 
In other words, the contravariant Piola transformation preserves the normal component. 

In this document, integration is implicitly computed using numerical quadrature on the reference element. 
Integration over surfaces is performed over the one-dimensional reference element using the transformed element of length. 

\subsection{Finite Element Spaces}
Here, we define the finite element spaces used to approximate the transport and SMM moment systems. 
These finite element spaces are defined on the mesh $\meshT$ or the interior skeleton of the mesh $\Gamma_0$ and consist of an element-local function space and a set of inter-element matching conditions. 
The combination of a smooth element-local space and suitable matching conditions allows these finite-dimensional spaces to be subspaces of Sobolev spaces such as $L^2(\D)$, $H^1(\D)$, and $H(\div;\D)$. 

\subsubsection{Discontinuous Galerkin}
Let $\Qbb{p}(K)$ denote the space of polynomials of at most degree $p$ mapped from the reference element. 
That is, 
	\begin{equation}
		\Qbb{p}(K) = \{ u = \hat{u} \circ \T^{-1} : \hat{u} \in \Q_p(\hat{K}) \} \,, 
	\end{equation}
where $\hat{u}$ indicates a function defined on the reference element. 
The delineation between $\Q$ and $\mathbb{Q}$ is required when non-affine mesh transformations are used. 
In such case, $u = \hat{u} \circ \T^{-1} \notin \Q_p(K)$ even if $\hat{u} \in \Q_p(\hat{K})$. 
In other words, the solution can be non-polynomial due to the composition with the inverse mesh transformation.  

\begin{figure}
\centering
\includegraphics[width=.35\textwidth]{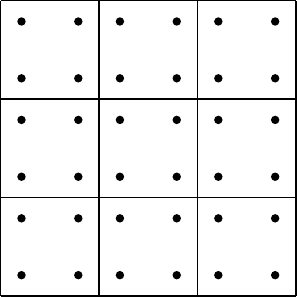}
\caption{A depiction of the distribution of degrees of freedom in the linear DG space. The nodal basis on each element is built from the Legendre nodes to illustrate that degrees of freedom are not shared between elements.}
\label{fig:dgfes}
\end{figure}

The Discontinuous Galerkin (DG) space is a finite-dimensional subspace of the $L^2(\D)$ space defined by 
	\begin{equation}
		L^2(\D) = \{ u : \int u^2 \ud \x < \infty \} \,,
	\end{equation}
the space of square-integrable functions. 
We denote the degree-$p$ DG space with 
	\begin{equation}
		Y_p = \{ u \in L^2(\D) : u|_K \in \Qbb{p}(K) \,, \quad \forall K \in \meshT \}\,,
	\end{equation}
so that each function $u \in Y_p$ is a piecewise polynomial mapped from the reference element with no continuity requirements enforced between elements. 
Figure \ref{fig:dgfes} shows the distribution of degrees of freedom in a $3\times 3$ mesh for the space $Y_1$. 
The Gauss-Legendre interpolation points are used to define the nodal basis for the local polynomial space to highlight that degrees of freedom are not shared between elements. 

\subsubsection{Continuous Finite Element}
\begin{figure}
\centering
\includegraphics[width=.35\textwidth]{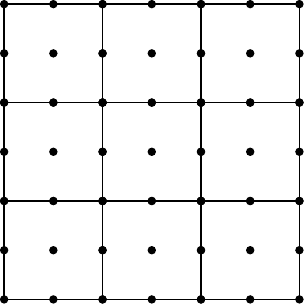}
\caption{A depiction of the distribution of degrees of freedom for the quadratic continuous finite element space. Continuity of members of the finite element space is enforced by sharing degrees of freedom across neighboring elements.}
\label{fig:cgfes}
\end{figure}

The continuous finite element (CG) space is a finite-dimensional subspace of the $H^1(\D)$ Sobolev space defined as: 
	\begin{equation}
		H^1(\D) = \{ u \in L^2(\D) : \int \nabla u \cdot \nabla u \ud \x < \infty \} \,,
	\end{equation}
the space of square-integrable functions with square-integrable gradient. 
It can be shown that, if a piecewise function $u \in C^0(\D)$ -- the space of continuous functions with zero continuous derivatives -- and satisfies $u|_K \in H^1(K)$ for each $K \in \meshT$ , then $u \in H^1(\D)$ (see \citet[\S3.2.1]{quateroni}). 
In other words, a finite-dimensional subspace of $H^1(\D)$ must be continuous across interior mesh interfaces and have a locally smooth function space on each element. 
Thus, we take the degree-$p$ continuous finite element space to be: 
	\begin{equation}
		V_p = \{ u \in C^0(\D) : u|_K \in \Qbb{p}(K) \,, \quad \forall K \in \meshT \} \,. 
	\end{equation}
Here, $V_p \subset H^1(\D)$ since $u \in V_p$ is continuous and $u|_K \in \Qbb{p}(K) \subset H^1(K)$ for each $K$. 
We use a nodal basis for $\Q_p(\hat{K})$ that employs points on the boundary of the element such as the Gauss-Lobatto points. 
This allows the sharing of degrees of freedom between elements that strongly enforces the continuity requirements of the discrete space. 
The distribution of unknowns for the space $V_2$ is shown in Fig.~\ref{fig:cgfes}. 
Observe that degrees of freedom are shared between neighboring elements along the interior interfaces between elements. 

\subsubsection{Raviart Thomas, Broken Raviart Thomas, and the Raviart Thomas Trace Space}
The Raviart Thomas (RT) space is a finite-dimensional subspace of the $H(\div;\D)$ space where 
	\begin{equation}
		H(\div;\D) = \{ \vec{v} \in [L^2(\D)]^2 : \nabla\cdot \vec{v} \in L^2(\D) \} \,,
	\end{equation}
is the space of square-integrable, vector-valued functions with square-integrable divergence \cite{Raviart1977,raviart_thomas}. 
It can be shown that $\vec{v} \in H(\div;\D)$ when $\vec{v}$ satisfies: 1) $\vec{v}$ is locally smooth on each element such that $\vec{v}|_K \in [H^1(K)]^2$ for each $K \in \meshT$ and 2) $\vec{v}\cdot\n$ is continuous across each interior face. 
It is also desirable that the divergence of a vector in the degree-$p$ Raviart Thomas space belong to the degree-$p$ DG space, $Y_p$. 
Let $\Q_{m,n}(\hat{K})$ denote the space of polynomials of at most $m$ and $n$ in the first and second variables, respectively, such that $\Q_{p,p}(\hat{K}) = \Q_p(\hat{K})$. 
The Raviart Thomas space uses the local polynomial space $\Q_{p+1,p}(\hat{K})\times\Q_{p,p+1}(\hat{K})$ on each element. 
In this way, the divergence of any vector $\vec{v} \in \Q_{p+1,p}(\hat{K})\times\Q_{p,p+1}(\hat{K})$ belongs to $\Q_p(\hat{K})$, the polynomial space used by $Y_p$ on each element. 
Finally, the contravariant Piola transformation, defined in Eq.~\ref{eq:piola}, is used to facilitate the sharing of the degrees of freedom that represent the normal component across an interior face in the mesh. 

\begin{figure}
\centering
\begin{subfigure}{.35\textwidth}
\centering
\includegraphics[width=\textwidth]{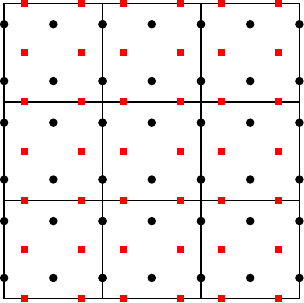}
\caption{}
\label{fig:rtfes}
\end{subfigure}
\hspace{1cm}
\begin{subfigure}{.35\textwidth}
\centering
\includegraphics[width=\textwidth]{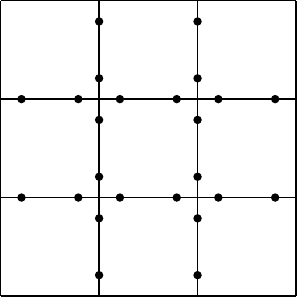}
\caption{}
\label{fig:ifes}
\end{subfigure}
\caption{(a) The distribution of degrees of freedom corresponding to the first degree Raviart Thomas space. Continuity of the normal component is enforced by sharing the degrees of freedom corresponding to the normal component along the interior face between neighboring elements. The circles and squares denote degrees of freedom representing the $x$ and $y$ components, respectively. (b) The distribution of degrees of freedom corresponding to $\Lambda_1$, the space defined as the interior trace of the first degree Raviart Thomas space, on a $3\times 3$ mesh. Observe that the interpolation for the normal component in (a) exactly matches the interpolation used in (b) on each face. }
\end{figure}

The degree-$p$ Raviart Thomas space is then 
	\begin{equation}
		\RT_p = \{ \vec{v} \in [L^2(\D)] : \vec{v}|_K \in \mathbb{D}_p(K) \,, \ \forall K \in \meshT \ \text{and}\ \jump{\vec{v}\cdot\n} = 0 \,, \ \forall \mathcal{F} \in \Gamma_0\} \,,
	\end{equation}
where 
	\begin{equation}
		\mathbb{D}_p(K) = \{ \vec{v} = \frac{1}{J}\F\hvec{v}\circ\T^{-1} : \hvec{v} \in \Q_{p+1,p}(\hat{K}) \times \Q_{p,p+1}(\hat{K}) \}
	\end{equation}
denotes the required Raviart Thomas element-local polynomial space composed with the contravariant Piola transformation. 
The constraint $\jump{\vec{v}\cdot\n} = 0$ for each interior face holds only when $\vec{v} \in \RT_p$ has a continuous normal component. 
Thus, since $\mathbb{D}_p(K) \subset [H^1(K)]^2$ and each element of $\RT_p$ has a continuous normal component, $\RT_p \subset H(\div;\D)$. 
Figure \ref{fig:rtfes} depicts the degrees of freedom in $\RT_1$ on a $3\times 3$ mesh. 
The black circles and red squares denote degrees of freedom representing the $x$ and $y$ components, respectively. 
Observe that on each element the $x$ component is quadratic in the $x$ direction and linear in the $y$ direction with the opposite holding for the $y$ component. 
Furthermore, the degrees of freedom representing the normal component are shared between elements, ensuring that $\jump{\vec{v}\cdot\n} = 0$ for each interior face and $\vec{v} \in \RT_p$. 

We also define two related spaces that arise in hybridizing a Raviart Thomas discretization. 
The first is the broken Raviart Thomas space, $\hRT_p$, defined as the degree-$p$ Raviart Thomas space without the requirement of continuity in the normal component. 
That is, 
	\begin{equation}
		\hRT_p = \{ \vec{v} \in [L^2(\D)] : \vec{v}|_K \in \mathbb{D}_p(K) \,, \quad \forall K \in \meshT \} \,. 
	\end{equation}
Note that $\RT_p \subset \hRT_p$ and that $\vec{v} \in \hRT_p$ belongs to $\RT_p$ if and only if $\jump{\vec{v}\cdot\n} = 0$ on all interior mesh faces. 
For the broken space, the degrees of freedom are distributed as in Fig.~\ref{fig:rtfes}, however, degrees of freedom are not shared between elements. 
Instead, these degrees of freedom are duplicated so that the normal component is doubly defined on interior interfaces. 

Second, we define the interior trace of the Raviart Thomas space. 
This space represents the normal component of $\RT_p$ along the interior faces in the mesh. 
Let $\mathcal{P}_p(\hat{\mathcal{F}})$ be the space of univariate polynomials with degree at most $p$ defined over the reference line $\hat{\mathcal{F}} = [0,1]$ and 
	\begin{equation}
		\mathbb{P}(\mathcal{F}) = \{ u = \hat{u} \circ \T^{-1} : \hat{u} \in \mathcal{P}_p(\hat{\mathcal{F}})\} \,,
	\end{equation}
the associated polynomials mapped to physical space.
The trace space is then 
	\begin{equation}
		\Lambda_p = \{ \mu \in L^2(\Gamma_0) : \mu|_\mathcal{F} \in \mathbb{P}_p(\mathcal{F}) \,, \quad \forall \mathcal{F} \in \Gamma_0 \} \,. 
	\end{equation}
The interior trace space for $\RT_1$, $\Lambda_1$, is depicted in Fig.~\ref{fig:ifes}. 
Note that these degrees of freedom are exactly the degrees of freedom corresponding to the normal component of $\RT_1$ on the interior faces of the mesh. 

\section{Transport Discretization}
We discretize the \Sn transport equations with a high-order DG spatial discretization compatible with curved meshes such as from \citet{woods_thesis} or \citet{graph_sweeps}. 
In \Sn, the transport equation is collocated at discrete angles, $\Omegahat_d$, and integration over the unit sphere is approximated using a suitable angular quadrature rule $\{\Omegahat_d,w_d\}_{d=1}^{N_\Omega}$. 
We use $\psi_d$ to denote the angular flux in discrete direction $\Omegahat_d$ such that $\psi_d(\x) = \psi(\x,\Omegahat_d)$. 
The DG discretization approximates each $\psi_d$ using the degree-$p$ DG finite element space, $Y_p$. 
Through finite element interpolation, $\psi_d(\x)$ can be evaluated at any point in the mesh. 
The SMM closures are computed using the \Sn angular quadrature and finite element interpolation with: 
	\begin{subequations}
	\begin{equation} \label{eq:correction_tensor_disc}
		\T(\x) = \sum_{d=1}^{N_\Omega} w_d\,\Omegahat_d\otimes\Omegahat_d\,\psi_d(\x) - \frac{1}{3}\I \sum_{d=1}^{N_\Omega} w_d\psi_d(\x) \,,  
	\end{equation}
	\begin{equation} \label{eq:bdr_correction_factor_disc}
		\beta(\x) = \sum_{d=1}^{N_\Omega} w_d\,|\Omegahat_d\cdot\n|\,\psi_d(\x) - E_{b0} \sum_{d=1}^{N_\Omega} w_d\psi_d(\x) \,. 
	\end{equation}
	\end{subequations}
In light of Remark \ref{rem:disc_Eb}, we elect to use $E_{b0}$, as defined in Eq.~\ref{eq:eb0}, in place of the factor of $1/2$ that scales the zeroth moment term in the boundary correction. 
Here, $E_{b0}$ is approximated with \Sn quadrature according to: 
	\begin{equation} \label{eq:eb0_disc}
		E_{b0} = \frac{\sum_{d=1}^{N_\Omega} w_d \, |\Omegahat_d\cdot\n|}{\sum_{d=1}^{N_\Omega} w_d }  = \frac{\sum_{d=1}^{N_\Omega} w_d \, |\Omegahat_d\cdot\n|}{4\pi}\,. 
	\end{equation}
With this discrete definition, the boundary conditions for the continuous system from Eq.~\ref{eq:smm_bc} are now modified to: 
	\begin{equation} \label{eq:smm_bc_disc}
		\vec{J}\cdot\n = E_{b0} \varphi + 2\Jin + \beta \,, 
	\end{equation}
where $\beta$ is defined using the discrete definition in Eq.~\ref{eq:bdr_correction_factor_disc}. 
We also approximate the inflow partial current, $\Jin$, using \Sn quadrature with $\Jin(\x) = \sum_{\Omegahat_d\cdot\n<0} w_d\,\Omegahat_d\,\bar{\psi}(\x,\Omegahat_d)$. 
In the limit as $N_\Omega \rightarrow \infty$, $E_{b0} \rightarrow 1/2$ and the discrete $\Jin \rightarrow \int_{\Omegahat\cdot\n<0} \Omegahat\cdot\n\, \bar{\psi} \ud \Omega$ so that the continuous boundary condition is recovered. 
This discrete prescription ensures that the boundary correction factor will vanish when the discrete transport solution is a linear function in angle and was found to be a requirement to obtain the correct moment solution and preserve the optimal convergence rate of the moment algorithm in the thick diffusion limit. 
Note that we use the same notation to denote the analytic and discrete SMM closures for notational simplicity. 
For the remainder of this document, only the discrete definitions of the closures given by Eqs.~\ref{eq:correction_tensor_disc} and \ref{eq:bdr_correction_factor_disc} will be used. 

Where needed, we compute the broken divergence of the SMM correction tensor using the finite element interpolation operator for $\T$. 
In other words, 
	\begin{equation} \label{eq:broken_Tdiv}
		\nablah\cdot\T = \nablah\cdot\P - \frac{1}{3} \nablah \phi \,,
	\end{equation}
where $\P$ and $\phi$ are the second and zeroth moments of the discrete angular flux, respectively, and the broken derivatives are applied to the finite element interpolation operators used to represent the angular flux on each element. 

\begin{rem} \label{rem:exact_integration}
For VEF, the discrete closures are locally represented in space as ratios of polynomials and are thus rational functions (see \citet[\S4]{dgvef_olivier}). 
The bilinear and linear forms that comprise the VEF moment system discretization will then contain numerical integration error since rational functions cannot be integrated exactly with numerical quadrature. 
For SMM, the linear closures simply inherit the finite element representation used for the angular flux. 
In this way, the SMM moment system can be integrated exactly with numerical quadrature. 
\end{rem}
\begin{rem} \label{rem:nofixup}
The VEF closures are normalized and thus are only well defined for strictly positive angular fluxes, requiring the use of a negative flux fixup for robustness. 
The SMM closures are not normalized and are thus well defined even when the transport solution contains zeros or non-physical negativities arising from numerical error.  
This could be particularly beneficial in shielding or other deep-penetration problems where the solution is significantly attenuated. 
Here, we leverage this behavior to build an algorithm that remains robust to numerically induced negativities both with and without the use of a negative flux fixup. 
\end{rem}

The SMM scalar flux is coupled to the transport equation through the scattering source. 
We use a mixed-space scattering mass matrix to support the use of differing finite element spaces for the SMM scalar flux and the angular flux. 
In other words, the scattering mass matrix is assembled using test functions from the space used to approximate the angular flux and trial functions from the space used to approximate the SMM scalar flux. 

\section{Discrete Second Moment Methods} \label{sec:smm}
In this section, we leverage the close connection between VEF and SMM to systematically convert the existing VEF methods derived in \citet{dgvef_olivier} and \citet{rtvef_olivier} to discretizations of the SMM moment system. 
In particular, we derive analogs of the interior penalty discontinuous Galerkin, continuous finite element, Raviart Thomas mixed finite element, and hybridized Raviart Thomas mixed finite element methods. 
The interior penalty method is derived by linearizing the interior penalty VEF discretization about a linearly anisotropic solution. 
The Raviart Thomas mixed finite element method is derived by algebraically manipulating the Raviart Thomas VEF discretization. 
The continuous finite element and hybridized Raviart Thomas methods are found by leveraging their close connection to the interior penalty and unhybridized Raviart Thomas methods, respectively. 
We have elected to derive the interior penalty and Raviart Thomas methods via linearization and algebraic manipulation, respectively, to give examples of these two equivalent strategies. 
However, either method can be used to convert any discretization of the VEF moment system to a discretization of the SMM moment system. 

\subsection{Discontinuous Galerkin and Continuous Finite Element}
Here, we linearize an existing DG VEF method about a linearly anisotropic solution to derive an analogous SMM moment system discretization. 
Note that since the transport equation is linear, the linearization process does not alter the transport equation, and thus we only need to linearize the VEF moment system to derive the SMM moment system. 
Let $\mat{V}(\psi,\varphi) = f$ represent the DG VEF moment system. 
Here, $\mat{V}$ is nonlinear in $\psi$, due to the nonlinear Eddington tensor and boundary factor, but linear in $\varphi$. 
Further, we set $\y_0 = \vector{\psi_0 & \varphi_0}$ to be the linearly anisotropic approximation to the coupled transport-VEF system. 
The linearized operator is then: 
	\begin{equation}
	\begin{aligned}
		0 = \mat{V}(\psi,\varphi) &\xrightarrow{\text{TSE}} \mat{V}(\psi_0,\varphi_0) + D[\mat{V}](\y_0, \y- \y_0) \\
		&= \mat{V}(\psi_0,\varphi_0) + D[\mat{V}\rvert_{\varphi = \varphi_0}](\psi_0,\psi - \psi_0) + D[\mat{V}\rvert_{\psi = \psi_0}](\varphi_0, \varphi - \varphi_0) \\
		&= \mat{V}(\psi_0,\varphi) + D[\mat{V}\rvert_{\varphi=\varphi_0}](\psi_0, \psi) \,,
	\end{aligned}
	\end{equation}
where $D[\cdot](\cdot,\cdot)$ denotes the Gateaux derivative defined in Eq.~\ref{eq:gateaux}.
We have simplified $\mat{V}(\psi_0,\varphi_0) + D[\mat{V}\rvert_{\psi = \psi_0}](\varphi_0, \varphi - \varphi_0) = \mat{V}(\psi_0,\varphi)$ using the linearity of $\mat{V}$ in the second argument.
Furthermore, the directional derivative of the Eddington tensor satisfies $D[\E](\psi_0,\psi-\psi_0) = D[\E](\psi_0,\psi)$ with an analogous result for the directional derivative of the Eddington boundary factor. 
Since $\psi_0$ is linearly anisotropic, $\mat{V}(\psi_0,\varphi)$ represents a radiation diffusion system. 
The correction sources are then given by $D[\mat{V}\rvert_{\varphi=\varphi_0}](\psi_0,\psi)$: the directional derivative of the VEF operator evaluated at $\psi_0$ in the direction $\psi$ where $\varphi$ is fixed at the linearly anisotropic moment solution, $\varphi_0$. 
We proceed by evaluating the interior penalty VEF method from \citet[\S 5.2.1]{dgvef_olivier} at a linearly anisotropic transport solution, $\psi_0$, to define the diffusion discretization and then derive the correction source by applying the Gateaux derivative to each term in the VEF discretization. 

From \citet[Eq.~59]{dgvef_olivier}, the interior penalty VEF discretization is: find $\varphi \in Y_p$ such that 
	\begin{multline} \label{eq:ipvef}
		\int_{\Gamma_b} E_b\, u \varphi \ud s + \int_{\Gamma_0} \kappa \jump{u} \jump{\varphi} \ud s - \int_{\Gamma_0} \jump{u} \avg{\frac{1}{\sigma_t}\nablah\cdot\paren{\E\varphi}\cdot\n} \ud s - \int_{\Gamma_0} \avg{\frac{\nablah u}{\sigma_t}} \cdot \jump{\E\varphi\n} \ud s \\
		+ \int \nablah u \cdot \frac{1}{\sigma_t}\nablah\cdot\paren{\E\varphi} \ud \x + \int \sigma_a\, u \varphi \ud \x \\ 
		= \int u\, Q_0 \ud \x + \int \nablah u \cdot \frac{\vec{Q}_1}{\sigma_t} \ud \x - \int_{\Gamma_0} \jump{u} \avg{\frac{\vec{Q}_1\cdot\n}{\sigma_t}} \ud s - 2\int_{\Gamma_b} u\, \Jin \ud s \,, \quad \forall u \in Y_p \,, 
	\end{multline}
where $\kappa$ is the penalty parameter. 
This discretization is derived by approximating the scalar flux and each component of the current with a degree-$p$ DG finite element space. 
The first moment equation is integrated-by-parts twice in order to derive a discrete elimination of the current. 
The numerical flux is chosen so that the current can be eliminated on each element, leading to a discretization of the second-order form of the VEF equations. 
The symmetric, positive semi-definite penalty bilinear form $\int_{\Gamma_0} \kappa\, \jump{u}\jump{\varphi} \ud s$ is added to stabilize the discretization. 
It is well known that $\kappa \propto \sigma_t^{-1} p^2 / h$ is required to ensure that the penalty bilinear form dominates the negative definite terms in the discretization, resulting in an overall linear system that is positive definite and stable with respect to $h$ and $p$ \cite{Arnold1982,Arnold2002}. 
We note that the penalty parameter's proportionality constant is often problem-dependent. 
While the penalty term is required for stability, the presence of the penalty term degrades the performance of multigrid preconditioners and thus specialized solvers, such as the uniform subspace correction preconditioner developed by \citet{Pazner2021}, must be used to achieve a preconditioned iterative solver that converges independent of the mesh size, polynomial degree, and penalty parameter. 

\citet{dgvef_olivier} also present discretizations of the VEF moment system that are analogs of the second method of Bassi and Rebay (BR2) \cite{Bassi2000} and the Local Discontinuous Galerkin (LDG) method \cite{Cockburn1998}.
These methods use alternate stabilization techniques that avoid the problem-dependent prescription of the penalty parameter. 
However, the BR2 and LDG stabilization terms are more expensive to compute (and more complicated to implement). 
Only minor differences in solution quality and iterative efficiency were seen between IP, BR2, and LDG, and CG. 
Thus, for clarity of presentation we present only the IP and CG methods. 

Evaluating the VEF data when the angular flux is linearly anisotropic gives
	\begin{equation}
		\E(\psi_0) = \frac{1}{3}\I \,, \quad E_b(\psi_0) = E_{b0} \,, 
	\end{equation}
where $E_{b0}$ is computed using the \Sn quadrature rule following Eq.~\ref{eq:eb0_disc}. 
The use of numerical quadrature to evaluate this term is motivated by Remark \ref{rem:disc_Eb}. 
The diffusion operator is then 
	\begin{multline}
		\int_{\Gamma_b} E_{b0}\, u\varphi \ud s + \int_{\Gamma_0} \kappa\, \jump{u}\jump{\varphi} \ud s - \int_{\Gamma_0} \jump{u}\avg{D\nablah\varphi\cdot\n} \ud s - \int_{\Gamma_0} \avg{D\nablah u\cdot\n} \jump{\varphi} \ud s \\
		+ \int \nablah u \cdot D\nablah\varphi \ud \x + \int \sigma_a\, u \varphi \ud \x \\
		= \int u\, Q_0 \ud \x + \int \nablah u \cdot \frac{\vec{Q}_1}{\sigma_t} \ud \x - \int_{\Gamma_0} \jump{u} \avg{\frac{\vec{Q}_1\cdot\n}{\sigma_t}} \ud s - 2\int_{\Gamma_b} u\, \Jin \ud s \,,
	\end{multline}
where $D = \frac{1}{3\sigma_t}$ is the diffusion coefficient. 
The above corresponds to a standard interior penalty discretization of radiation diffusion with Marshak boundary conditions. 
Note that the right hand side includes additional sources depending on $\vec{Q}_1$ since we have not made the assumption that the transport equation's fixed-source, $q$, is at most linearly anisotropic, as is customary for the diffusion approximation. 
Inclusion of the first moment of the source facilitates investigating the accuracy of the method with the method of manufactured solutions on a transport problem as performed in \S\ref{sec:mms}. 

Next, we must determine the correction terms by computing the directional derivative of each term in the VEF discretization (Eq.~\ref{eq:ipvef}) with $\varphi = \varphi_0$ fixed. 
For terms without VEF data, the directional derivative is zero as they do not depend on the transport solution. 
For the boundary factor-dependent term, the correction term is:  
	\begin{equation}
		D\bracket{\int_{\Gamma_b} E_b\, u \varphi_0 \ud s}(\psi_0, \psi) = \int_{\Gamma_b} D[E_b](\psi_0,\psi)\, u \varphi_0 \ud s = \int_{\Gamma_b} u\,\beta(\psi) \ud s \,, 
	\end{equation}
where $\beta(\psi)$ is the correction factor defined in Eq.~\ref{eq:bdr_correction_factor_disc} that uses \Sn quadrature to compute the angular moments of the discrete transport solution. 
We have used $D[E_b](\psi_0,\psi) = \beta(\psi)/\phi_0$ and have made the approximation that $\varphi_0/\phi_0 \approx 1$ (see Remark \ref{rem:independent_phis}). 
Analogously, for terms with the Eddington tensor
	\begin{equation}
		D[\E \varphi_0](\psi_0,\psi) = \frac{\varphi_0}{\phi_0}\T(\psi) = \T(\psi) \,, 
	\end{equation}
where $\T(\psi)$ is the correction tensor from Eq.~\ref{eq:correction_tensor_disc} and we have again made the approximation that $\varphi_0/\phi_0 = 1$. 
Thus, the correction terms can be derived by setting terms without angular flux dependence to zero and replacing 
	\begin{equation}
		E_b\varphi \rightarrow \beta \,, \quad \E\varphi \rightarrow \T \,. 
	\end{equation}
The interior penalty SMM discretization is then: find $\varphi \in Y_p$ such that 
	\begin{multline} \label{eq:ipsmm}
		\int_{\Gamma_b} E_{b0}\,u\varphi \ud s + \int_{\Gamma_0} \kappa\, \jump{u}\jump{\varphi} \ud s - \int_{\Gamma_0} \jump{u}\avg{D\nablah\varphi\cdot\n} \ud s - \int_{\Gamma_0} \avg{D\nablah u\cdot\n} \jump{\varphi} \ud s \\
		+ \int \nablah u \cdot D\nablah\varphi \ud \x + \int \sigma_a\, u \varphi \ud \x \\
		= \int u\, Q_0 \ud \x + \int \nablah u \cdot \frac{\vec{Q}_1}{\sigma_t} \ud \x - \int_{\Gamma_0} \jump{u} \avg{\frac{\vec{Q}_1\cdot\n}{\sigma_t}} \ud s - \int_{\Gamma_b} u\paren{2\Jin + \beta} \ud s \\
		+ \int_{\Gamma_0} \jump{u}\avg{\frac{1}{\sigma_t}\nablah\cdot\T\cdot\n} \ud s + \int_{\Gamma_0} \avg{\frac{\nablah u}{\sigma_t}} \cdot\jump{\T\n} \ud s - \int \nablah u\cdot \frac{1}{\sigma_t}\nablah\cdot\T \ud \x \,, \quad \forall u \in Y_p \,, 
	\end{multline}
where the local divergence of the correction tensor is computed with Eq.~\ref{eq:broken_Tdiv}. 
The above represents a standard interior penalty discretization of diffusion with Marshak boundary conditions that is corrected by transport-dependent volumetric and boundary source terms. 

A continuous finite element discretization can be extracted from Eq.~\ref{eq:ipsmm} by setting the test and trial spaces to $V_p$, the degree-$p$ continuous finite element space. 
For any function $v \in V_p$, 
	\begin{equation}
		\jump{v} = 0 \,, \quad \avg{v} = v \,, \quad \forall \mathcal{F} \in \Gamma_0 \,. 
	\end{equation}
In addition, $v\in V_p$ satisfies $\nablah v = \nabla v$ \cite[Prop.~3.2.1]{quateroni}. 
In other words, the space $V_p$ has the required continuity to allow the broken and global gradients to be equal for any $v \in V_p$. 
Since DG is used for the transport discretization, the correction tensor is still generally discontinuous across interior mesh interfaces such that $\jump{\T\n}$ is non-zero. 
Thus, the CG SMM moment discretization is: find $\varphi \in V_p$ such that 
	\begin{multline}
		\int_{\Gamma_b} E_{b0}\,u\varphi \ud s + \int \nabla u \cdot D\nabla\varphi \ud \x + \int \sigma_a\, u \varphi \ud \x \\
		= \int u\, Q_0 \ud \x + \int \nabla u \cdot \frac{\vec{Q}_1}{\sigma_t} \ud \x - \int_{\Gamma_b} u\paren{2\Jin + \beta} \ud s \\
		+ \int_{\Gamma_0} \avg{\frac{\nabla u}{\sigma_t}} \cdot\jump{\T\n} \ud s - \int \nabla u\cdot \frac{1}{\sigma_t}\nablah\cdot\T \ud \x \,, \quad \forall u \in V_p \,. 
	\end{multline}
The CG SMM moment discretization can also be derived directly from the CG VEF moment discretization in \citet[Eq.~77]{dgvef_olivier} via linearization or algebraically manipulating the closures. 

\subsection{Raviart Thomas and Hybridized Raviart Thomas}
From \citet[Eqs.~61]{rtvef_olivier}, the RT VEF discretization is: find $(\varphi,\vec{J}) \in Y_p \times \RT_p$ such that 
	\begin{subequations} \label{eq:rtvef}
	\begin{equation}
		\int u\, \nabla\cdot\vec{J} \ud \x + \int \sigma_a\, u\varphi \ud \x = \int u\, Q_0 \ud \x \,, \quad \forall u \in Y_p \,,
	\end{equation}
	\begin{multline} \label{eq:rtvef_first}
		\int_{\Gamma_0} \jump{\vec{v}\cdot\avg{\E\n}} \avg{\varphi} \ud s - \int \nablah \vec{v} : \E\varphi \ud \x + \int \sigma_t\, \vec{v}\cdot\vec{J} \ud \x + \int_{\Gamma_b} \frac{1}{E_b}\!\paren{\vec{v}\cdot\E\n}\!\paren{\vec{J}\cdot\n} \ud s \\= \int \vec{v}\cdot\vec{Q}_1 \ud \x + 2\int_{\Gamma_b} \frac{1}{E_b}\vec{v}\cdot\E\n\, \Jin \ud s \,, \quad \forall \vec{v} \in \RT_p \,, 
	\end{multline}
	\end{subequations}
where $\vec{v},\vec{J} \in \RT_p$ are implicitly represented with the contravariant Piola transformation. 
Here, the presence of the Eddington tensor inside the strong form of the first moment's divergence term $\nabla\cdot\paren{\E\varphi}$ prevents the straightforward application of standard mixed finite element techniques. 
To illustrate this, we multiply this divergence term by a test function $\vec{v} \in \RT_p$ and integrate over a single element $K$: 
	\begin{equation}
		\int_K \vec{v}\cdot\nabla\cdot\paren{\E\varphi} \ud \x = \int_{\partial K} \vec{v} \cdot\widehat{\E\varphi}\n \ud s - \int_K \nabla\vec{v} : \E\varphi \ud \x \,,
	\end{equation}
where $\widehat{\E\varphi}\n$ is the numerical flux for the product of the Eddington tensor and VEF scalar flux that arises due to the discontinuous approximations used for the angular flux, and thus the Eddington tensor, and VEF scalar flux. 
When assembled into a global system, the element-boundary term gives rise to an interior face term of the form: 
	\begin{equation}
		\int_{\Gamma_0} \jump{\vec{v}\cdot\widehat{\E\varphi}\n} \ud s \,, 
	\end{equation}
along with a contribution on the boundary of the domain related to the boundary conditions that we omit for clarity of presentation. 
Note that $\vec{v} \in \RT_p$ has a continuous normal component but discontinuous tangential components. 
Thus, when $\E = \frac{1}{3}\I$, the interior face term vanishes since $\jump{\vec{v}\cdot \I\n} = \jump{\vec{v}\cdot\n} = 0$. 
For a general Eddington tensor, however, this term remains, requiring the application of DG-like techniques to treat the discontinuity in $\vec{v}\cdot\E\n$. 
We also note that the volumetric term $\int \nablah\vec{v} : \E\varphi \ud \x$ requires computing the gradient of an $\RT_p$ vector, a procedure significantly complicated by the use of the Piola transformation (see \cite[\S3.2]{rtvef_olivier} for more details). 
Note that \citet{rtvef_olivier} also present a VEF discretization where each component of the current is approximated with continuous finite elements. 
However, such a method could not be scalably preconditioned and solved even for a radiation diffusion problem where $\E = 1/3\I$ and $E_b = 1/2$. 
Thus, we not consider this method here. 

We first undo the choice of numerical flux and the application of the boundary conditions. 
We chose the following numerical flux for the product of the Eddington tensor and VEF scalar flux: 
	\begin{equation}
		\widehat{\E\varphi}\n = \avg{\E\n}\!\avg{\varphi} \,.
	\end{equation}
This choice allows the RT VEF discretization to limit to a standard RT discretization of radiation diffusion in the thick diffusion limit. 
Noting that VEF closes the pressure as $\P = \E \varphi$, we write the numerical flux term in terms of a numerical flux for the pressure denoted by $\widehat{\P}\n$. 
By solving for $\varphi$ in Eq.~\ref{eq:vef_bc}, the boundary scalar flux is given by: 
	\begin{equation}
		\varphi = \frac{1}{E_b}\paren{\vec{J}\cdot\n - 2\Jin} \,, \quad \x \in \partial \D \,, 
	\end{equation}
so that
	\begin{equation}
		\int_{\Gamma_b}\frac{1}{E_b}\!\paren{\vec{v}\cdot\E\n}\!\paren{\vec{J}\cdot\n} \ud s - 2 \int_{\Gamma_b} \frac{1}{E_b}\!\vec{v}\cdot\E\n\,\Jin\ud s = \int_{\Gamma_b} \vec{v}\cdot\E\n\,\varphi \ud s \,. 
	\end{equation}
Thus, we can recast the the RT VEF first moment (Eq.~\ref{eq:rtvef_first}) in terms of the unclosed pressure as: 
	\begin{equation}
		\int_{\Gamma_b} \vec{v}\cdot\P\n \ud s + \int_{\Gamma_0} \jump{\vec{v}\cdot \widehat{\P}\n} \ud s - \int \nablah \vec{v} : \P \ud \x + \int \sigma_t\, \vec{v}\cdot\vec{J} \ud \x = \int \vec{v}\cdot\vec{Q}_1 \ud \x \,, \quad \forall \vec{v} \in \RT_p \,, 
	\end{equation}
where we have substituted the pressure in place of the product $\E\varphi$. 
In this form, the SMM discretization can be derived by defining the pressure on the interior of the domain, the numerical flux of the normal component of the pressure on interior faces, and the boundary conditions. 

On the interior of the domain, we use the SMM closure $\P = \T + 1/3\I\varphi$. 
The volumetric term is then 
	\begin{equation}
		\int \nablah \vec{v} : \P \ud \x = \int \nablah \vec{v} : \T \ud \x + \frac{1}{3}\int \nabla\cdot\vec{v}\, \varphi \ud \x \,. 
	\end{equation}
Here, $\nablah\vec{v} : \I\varphi = \nablah\cdot\vec{v} = \nabla\cdot\vec{v}$ since $\nablah\cdot\vec{v} = \nabla\cdot\vec{v}$ for $\vec{v} \in \RT_p$ \cite[Prop.~3.2.2]{quateroni}.  
For the numerical flux, we simply use the average: 
	\begin{equation}
		\widehat{\P}\n = \avg{\T\n + \frac{1}{3}\varphi \n} = \avg{\T\n} + \frac{1}{3}\avg{\varphi}\n \,. 
	\end{equation}
With this choice of numerical flux, the interior face bilinear form simplifies to 
	\begin{equation}
		\int_{\Gamma_0} \jump{\vec{v}\cdot\widehat{\P}\n} \ud s = \int_{\Gamma_0} \jump{\vec{v}} \cdot\avg{\T\n} \ud s + \frac{1}{3}\int_{\Gamma_0} \jump{\vec{v}\cdot\n} \avg{\varphi} \ud s = \int_{\Gamma_0} \jump{\vec{v}} \cdot\avg{\T\n} \ud s \,, 
	\end{equation}
where $\jump{\vec{v}\cdot\n} = 0$ since the normal component is continuous for $\vec{v}\in\RT_p$. 
Finally, on the boundary of the domain, we solve the SMM boundary conditions in Eq.~\ref{eq:smm_bc_disc} for the scalar flux and insert this relationship into $\P\n = \T\n + \n/3 \varphi$ to yield:  
	\begin{equation}
		\P\n = \T\n + \frac{1}{3}\varphi \n = \T\n + \frac{\n}{3E_{b0}}\paren{\vec{J}\cdot\n - 2\Jin - \beta} \,, \quad \x \in \partial \D \,. 
	\end{equation}
Thus, the RT SMM discretization is: find $(\varphi,\vec{J}) \in Y_p \times \RT_p$ such that 
	\begin{subequations} \label{eq:rtsmm}
	\begin{equation}
		\int u\, \nabla\cdot\vec{J} \ud \x + \int \sigma_a\, u\varphi \ud \x = \int u\, Q_0 \ud \x \,, \quad \forall u \in Y_p \,,
	\end{equation}
	\begin{multline}
		-\frac{1}{3}\int \nabla\cdot\vec{v}\, \varphi \ud \x + \int\sigma_t\, \vec{v}\cdot\vec{J} \ud \x + \int_{\Gamma_b} \frac{1}{3E_{b0}}(\vec{v}\cdot\n)(\vec{J}\cdot\n) \ud s = \int \vec{v}\cdot\Qone \ud \x - \int_{\Gamma_b} \vec{v}\cdot\T\n \ud s \\
		+\int_{\Gamma_b} \frac{1}{3E_{b0}}\vec{v}\cdot\n\paren{2\Jin + \beta} \ud s - \int_{\Gamma_0} \jump{\vec{v}}\cdot\avg{\T\n} \ud s + \int \nablah \vec{v} : \T \ud \x \,, \quad \forall \vec{v} \in \RT_p \,. 
	\end{multline}
	\end{subequations}
Observe that the left hand side of Eqs.~\ref{eq:rtsmm} is equivalent to an RT discretization of radiation diffusion with Marshak boundary conditions. 

The hybridized Raviart Thomas discretization can be derived by directly applying the standard hybridization process (see \citet{doi:10.1137/17M1132562} and the references therein) to the RT SMM moment system. 
That is, the current is approximated in the broken RT space, $\hRT_p$, and continuity of the current in the normal component is enforced with a Lagrange multiplier defined on the trace space of $\RT_p$, $\Lambda_p$. 
The use of a discontinuous approximation for the scalar flux and current allows eliminating these variables locally on each element in favor of the Lagrange multiplier. 
This reduced system is both significantly smaller than the block system and can be effectively preconditioned with AMG. 
Once the Lagrange multiplier has been solved for, the scalar flux and current can be found through element-local back substitution. 
A derivation of a hybridized radiation diffusion method is outlined in \citet[\S 7.1]{rtvef_olivier}. 

The hybridized system is: find $(\vec{J},\varphi,\lambda) \in \hRT_p \times Y_p \times \Lambda_p$ such that 
	\begin{subequations}
	\begin{equation}
		\int u\,\nablah\cdot\vec{J} \ud \x + \int \sigma_a\, u\varphi \ud \x = \int u\, Q_0 \ud \x \,, \quad \forall u \in Y_p \,, 
	\end{equation}
	\begin{equation}
		-\frac{1}{3}\int\nablah\cdot\vec{v}\, \varphi \ud \x + \int\sigma_t\,\vec{v}\cdot\vec{J} \ud \x + \int_{\Gamma_b} \frac{1}{3E_{b0}} (\vec{v}\cdot\n)(\vec{J}\cdot\n) \ud s + \int_{\Gamma_0} \jump{\vec{v}\cdot\n}\, \lambda \ud s = \mathcal{S} \,, \quad \forall \vec{v} \in \hRT_p
	\end{equation}
	\begin{equation} \label{eq:weak_cts}
		\int_{\Gamma_0} \mu\, \jump{\vec{J}\cdot\n} \ud s = 0 \,, \quad \forall \mu \in \Lambda_p \,, 
	\end{equation}
	\end{subequations}
where 
	\begin{equation}
		\mathcal{S} = \int \vec{v}\cdot\Qone \ud \x - \int_{\Gamma_b} \vec{v}\cdot\T\n \ud s 
		+\int_{\Gamma_b} \frac{1}{3E_{b0}}\vec{v}\cdot\n\paren{2\Jin + \beta} \ud s - \int_{\Gamma_0} \jump{\vec{v}}\cdot\avg{\T\n} \ud s + \int \nablah \vec{v} : \T \ud \x
	\end{equation}
is the source term for the RT SMM discretization. 
As noted in \citet{doi:10.1137/17M1132562}, hybridization can be viewed as an algebraic process similar to static condensation. 
Furthermore, it was shown in \citet[\S7.1]{rtvef_olivier} that the weak continuity condition in Eq.~\ref{eq:weak_cts} implies that the normal component of the current is continuous in a pointwise, strong sense. 
Thus, the unhybridized and hybridized RT discretizations are equivalent.
The details of an efficient implementation are discussed in \citet[\S7.3]{rtvef_olivier} and \citet{doi:10.1137/17M1132562}. 

\section{Results}
We now present numerical results demonstrating the iterative efficiency and computational performance of the discrete SMMs. 
The SMMs are compared to the VEF algorithms from \citet{dgvef_olivier} and \citet{rtvef_olivier}. 
The methods were implemented using the MFEM finite element framework \cite{mfem-paper}. 
We use MFEM's conjugate gradient and BiCGStab implementations along with BoomerAMG from the \emph{hypre} sparse linear algebra library \cite{hypre}. 
KINSOL, from the Sundials package \cite{hindmarsh2005sundials}, provided the Anderson-accelerated fixed-point solvers. 
As described in \citet[\S2]{hindmarsh2005sundials}, the fixed-point iteration is terminated when the max norm of the difference between successive iterates is below the iterative tolerance. 
We use the high-order DG \Sn transport solver from \citet{graph_sweeps} to invert the streaming and collision operator needed in the moment algorithms. 
We refer to the methods as IP, CG, RT, and HRT for the interior penalty, continuous finite element, Raviart Thomas, and hybridized Raviart Thomas methods, respectively. 
For IP, we use
	\begin{equation}
		\kappa = \avg{\frac{(p+1)^2}{\sigma_t h}} \,, 
	\end{equation}
following the standard prescription used for model elliptic problems \cite{Bassi1997,Arnold1982} unless otherwise specified. 

The IP methods are preconditioned using the subspace correction method introduced in \citet{Pazner2021} and extended to the non-symmetric VEF system in \citet[\S6]{dgvef_olivier}. 
This method is a two-stage preconditioner that applies AMG to the DG system assembled onto a continuous finite element space and uses a simple Jacobi smoother to attack errors present in the degrees of freedom located on mesh interfaces. 
The CG and HRT methods are preconditioned by AMG. 
Finally, the RT methods are preconditioned with block preconditioners \cite{benzi_golub_liesen_2005}. 
The RT methods admit a block system of the form: 
	\begin{equation} \label{eq:rt_block_sys}
		\begin{bmatrix} 
			\mat{M}_t & \mat{G} \\ \mat{D} & \mat{M}_a 
		\end{bmatrix} 
		\begin{bmatrix} 
			\vec{J} \\ \varphi 
		\end{bmatrix}
		= \begin{bmatrix} 
			g \\ f 
		\end{bmatrix} \,,
	\end{equation}
where $\mat{M}_t$ is the total interaction mass matrix, $\mat{M}_a$ the absorption mass matrix, $\mat{D}$ and $\mat{G}$ represent the discrete divergence and gradient, respectively, and $f$ and $g$ represent the source terms.  
For SMM, $\mat{G} = -1/3 \mat{D}^T$. 
For VEF, $\mat{G}$ depends on the Eddington tensor and thus $\mat{G}$ is in general not proportional to the transpose of the divergence matrix, $\mat{D}$. 
We use block diagonal and lower block triangular preconditioners of the form: 
	\begin{equation}
		\mat{P}_\text{diag} = \begin{bmatrix} 
			\mat{M}_t & \\ & \hmat{S} 
		\end{bmatrix} \,, \quad 
		\mat{P}_\text{tri} = \begin{bmatrix} 
			\mat{M}_t & \\ \mat{D} & \hmat{S}
		\end{bmatrix} \,,
	\end{equation}
where $\hmat{S} = \mat{M}_a - \mat{D}\hmat{M}_t^{-1}\mat{G}$ is an approximate Schur complement formed using the lumped total interaction mass matrix, $\hmat{M}_t$, computed by summing the off-diagonal entries of $\mat{M}_t$ into the diagonal and setting the off-diagonal entries to zero. 
Lumping an $\RT_p$ mass matrix produces a diagonal matrix and thus $\hmat{M}_t^{-1}$ can be formed efficiently without fill-in. 
The approximate inverse of the diagonal preconditioner is applied with approximate inverses of the diagonal blocks while the lower block triangular inverse is applied with block forward substitution. 
We use a single iteration of Gauss-Seidel and AMG to approximately invert the total interaction mass matrix and approximate Schur complement, respectively. 
Note that for MINRES to effectively solve the RT SMM system, the system must be represented in a symmetric indefinite form and the preconditioner must be symmetric positive definite. 
Thus, for RT SMM, we multiply the top row of Eq.~\ref{eq:rt_block_sys} by negative three so that the off-diagonal blocks are symmetric and the scaled $\mat{M}_t$ is negative definite. 
For the block diagonal preconditioner, $\mat{M}_t$ is scaled by positive three so that $\mat{P}_\text{diag}$ is symmetric, positive definite. 
Finally, we note that the lower block triangular preconditioner is a non-symmetric operator and thus cannot be used with MINRES.  
More details on the preconditioning strategies used to solve the moment systems are provided in \citet[\S6]{dgvef_olivier} and \citet[\S6.4]{rtvef_olivier}. 

\subsection{Method of Manufactured Solutions} \label{sec:mms}
The accuracy of the methods are determined with the method of manufactured solutions (MMS). 
The solution is set to 
	\begin{equation} \label{eq:mms_psi}
		\psi = \frac{1}{4\pi}\bracket{\alpha(\x) + \Omegahat\cdot\vec{\beta}(\x) + \Omegahat\otimes\Omegahat : \mat{\Theta}(\x)} \,,
	\end{equation}
where 
	\begin{subequations}
	\begin{equation}
		\alpha(\x) = \sin(\pi x) \sin(\pi y) + \delta \,, 
	\end{equation}
	\begin{equation}
		\vec{\beta}(\x) = \begin{bmatrix} 
			\sin\!\paren{\frac{2\pi(x+\omega)}{1+2\omega}}\sin\!\paren{\frac{2\pi(y+\omega)}{1+2\omega}}\\
			\sin\!\paren{\frac{2\pi(x+\omega)}{1+2\omega}}\sin\!\paren{\frac{2\pi(y+\omega)}{1+2\omega}}
		\end{bmatrix} \,,
	\end{equation}
	\begin{equation} \label{eq:mmsH}
		\mat{\Theta}(\x) = \begin{bmatrix} 
			\frac{1}{2}\sin\!\paren{\frac{3\pi(x+\zeta)}{1 + 2\zeta}}\sin\!\paren{\frac{3\pi(y+\zeta)}{1+2\zeta}}
			& \sin\!\paren{\frac{2\pi(x+\omega)}{1+2\omega}}\sin\!\paren{\frac{2\pi(y+\omega)}{1+2\omega}}\\
			\sin\!\paren{\frac{2\pi(x+\omega)}{1+2\omega}}\sin\!\paren{\frac{2\pi(y+\omega)}{1+2\omega}}
			& \frac{1}{4}\sin\!\paren{\frac{3\pi(x+\zeta)}{1 + 2\zeta}}\sin\!\paren{\frac{3\pi(y+\zeta)}{1+2\zeta}}
		\end{bmatrix} \,. 
	\end{equation}
	\end{subequations}
Here, $\delta = 1.25$ is used to ensure $\psi>0$ and $\zeta = 0.1$ and $\omega = 0.05$ are used to test spatially dependent, non-isotropic inflow boundary conditions. 
The domain is $\D = [0,1]^2$. 
With this definition: 
	\begin{subequations}
	\begin{equation}
		\phi(\x) = \alpha(\x) + \frac{1}{3}\tr\mat{\Theta}(\x) \,,
	\end{equation}
	\begin{equation}
		\vec{J}(\x) = \frac{1}{3}\vec{\beta}(\x) \,,
	\end{equation}
	\begin{equation}
		\P(\x) = \frac{\alpha(\x)}{3}\I + \frac{1}{15}\begin{bmatrix} 
			3 \Theta_{11}(\x) + \Theta_{22}(\x) & \Theta_{12}(\x) \\ \Theta_{21}(\x) & \Theta_{11}(\x) + 3 \Theta_{22}(\x) 
		\end{bmatrix} \,. 
	\end{equation}
	\end{subequations}
This leads to an exact Eddington tensor $\E = \P/\phi$ that is dense and spatially varying. 
The MMS $\psi$ and $\phi$ are substituted into the transport equation to solve for the MMS source $q$ that forces the solution to Eq.~\ref{eq:mms_psi}. 

The accuracy of the SMM discretizations are investigated in isolation by computing the SMM closures from the MMS angular flux and setting the sources $Q_0$ and $\vec{Q}_1$ to the moments of the MMS source. 
This is accomplished by projecting the MMS angular flux onto a degree-$p$ DG finite element space and using Level Symmetric $S_4$ angular quadrature to compute the SMM closures, the moments of the MMS source, and the inflow partial current, $\Jin$. 
The SMM moment systems are then solved as if $\T$, $\beta$, $Q_0$, $\Qone$, and $\Jin$ are given data. 
Errors are calculated with the $L^2(\D)$ norm for scalars and the $[L^2(\D)]^2$ norm for vectors given by
	\begin{equation}
		\| u \| = \sqrt{\int u^2 \ud \x} \,,
	\end{equation}
and
	\begin{equation}
		\| \vec{v} \| = \sqrt{\int \vec{v}\cdot\vec{v} \ud \x}\,,
	\end{equation}
respectively. We also use the $L^2(\D)$ projection operator $\Pi_p : L^2(\D) \rightarrow Y_p$ such that 
	\begin{equation}
		\int u(v - \Pi_p v) \ud \x = 0 \,, \quad \forall u \in Y_p \,, 
	\end{equation}
for some $v \in L^2(\D)$. 
In particular, $\Pi_p$ is used to project the exact MMS scalar flux onto a $Y_p$ finite element grid function in order to investigate a superconvergence property of mixed finite elements. 

\begin{figure}
\centering
\includegraphics[width=.35\textwidth]{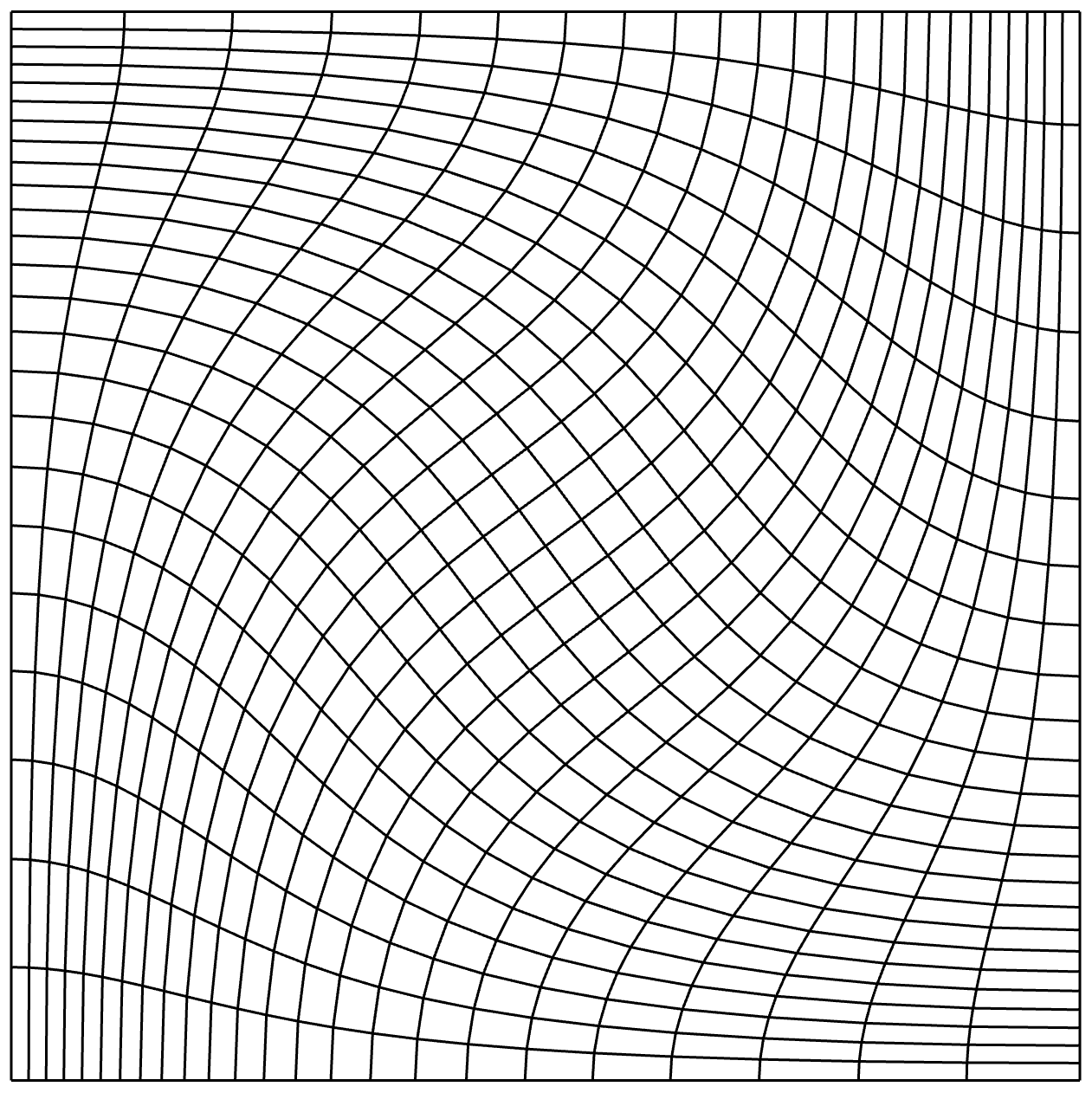}
\caption{A depiction of a third-order mesh generated by distorting an orthogonal mesh according to the Taylor Green vortex. Refinements of this mesh are used in calculating the error with the method of manufactured solutions.}
\label{fig:tgmesh}
\end{figure}
We use refinements of a third-order mesh created by distorting an orthogonal mesh according to the velocity field of the Taylor Green vortex. 
This mesh distortion is generated by advecting the mesh control points with 
	\begin{equation}
		\x = \int_0^T \mat{v} \ud t \,,
	\end{equation}
where the final time $T=0.3\pi$ and 
	\begin{equation}
		\mat{v} = \begin{bmatrix} 
			\sin(x) \cos(y) \\ 
			-\cos(x) \sin(y) 
		\end{bmatrix}
	\end{equation}
is the analytic solution of the Taylor Green vortex. 
300 forward Euler steps were used to advect the mesh. 
An example mesh is shown in Fig.~\ref{fig:tgmesh}. 

\begin{figure}
\centering
\begin{subfigure}{.32\textwidth}
	\centering
	\includegraphics[width=\textwidth]{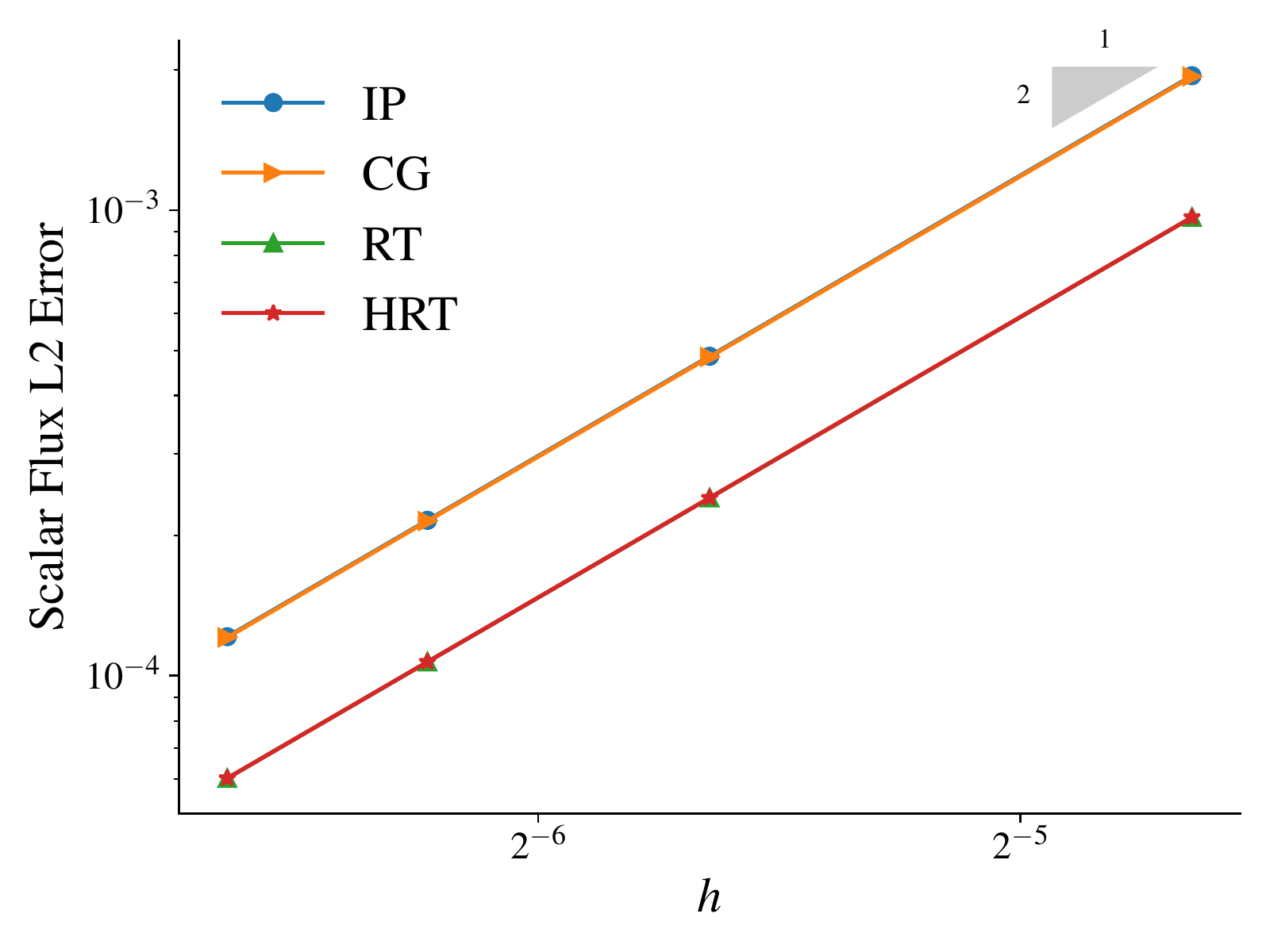} 
	\caption{}
\end{subfigure}
\begin{subfigure}{.32\textwidth}
	\centering
	\includegraphics[width=\textwidth]{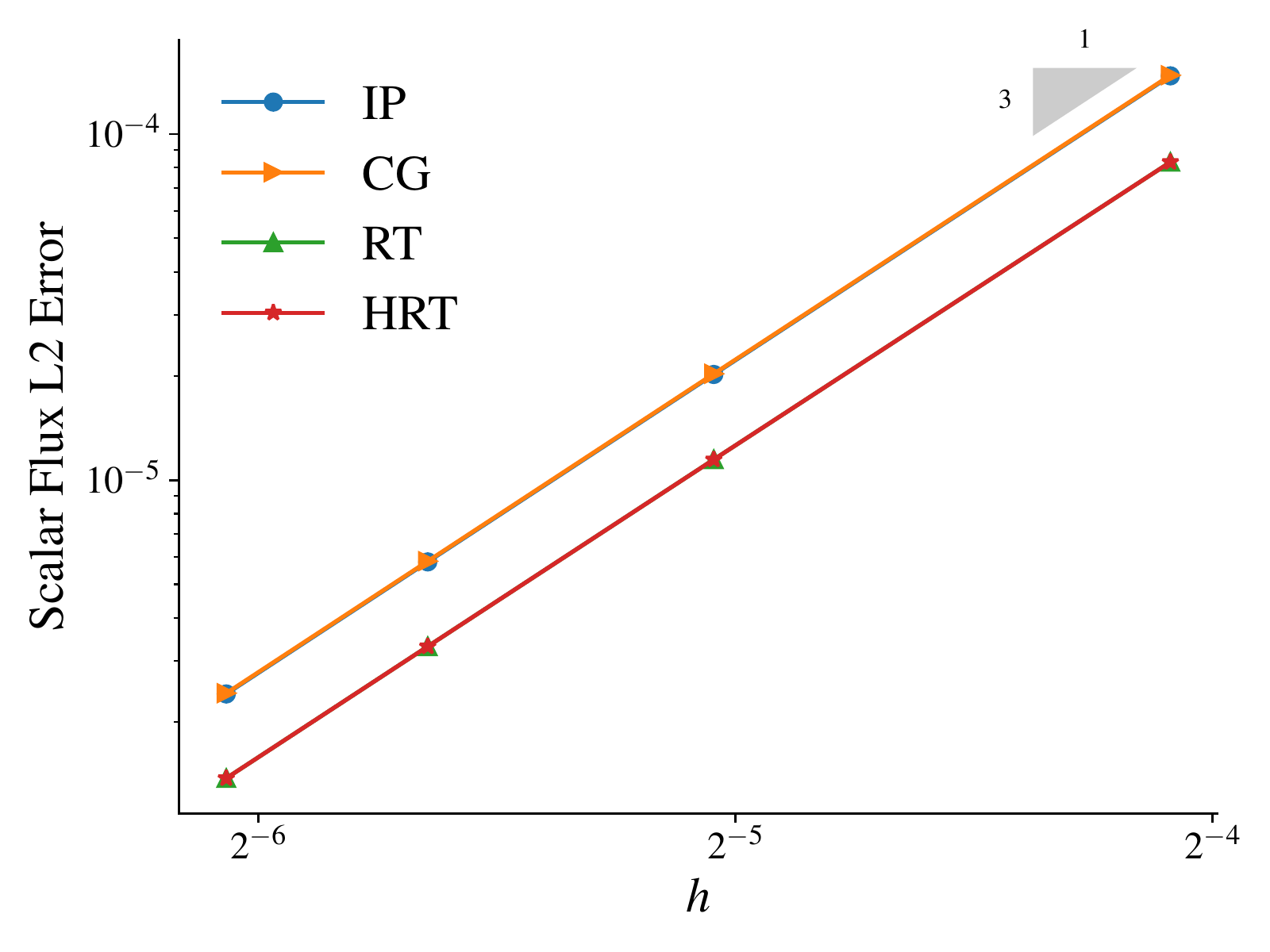} 
	\caption{}
\end{subfigure}
\begin{subfigure}{.32\textwidth}
	\centering
	\includegraphics[width=\textwidth]{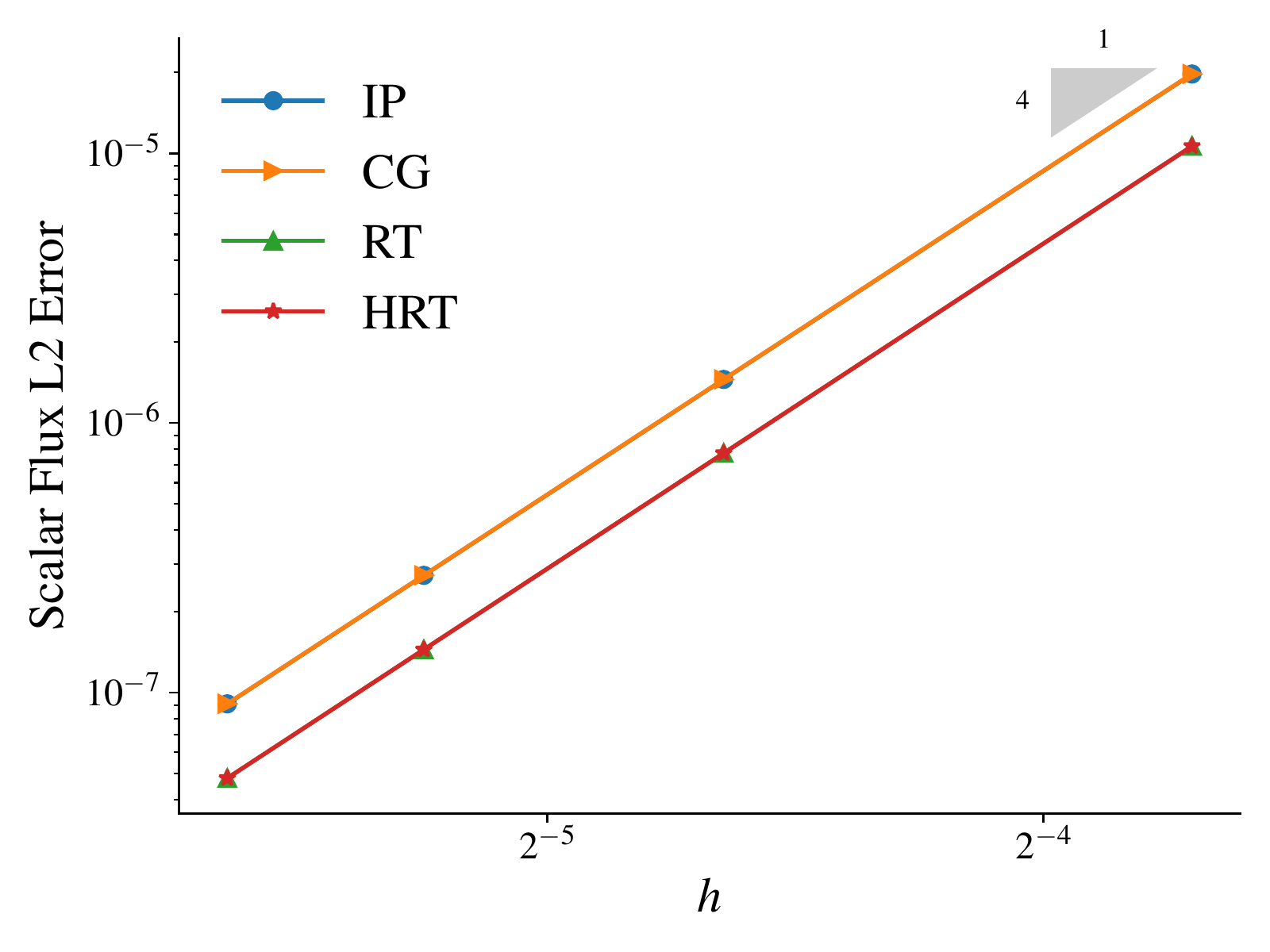} 
	\caption{}
\end{subfigure}
\caption{Plots of the error in the scalar flux as the mesh is refined for (a) linear, (b) quadratic, and (c) cubic finite elements. A quadratically anisotropic MMS transport problem is used. }
\label{fig:mms}
\end{figure}
Figure \ref{fig:mms} shows the error in the scalar flux on this MMS problem for a range of finite element polynomial degrees. 
All four methods achieve optimal $\mathcal{O}(h^{p+1})$ convergence. 
The IP and CG and RT and HRT methods have nearly identical error behavior, respectively, with the RT and HRT methods having a lower constant. 
This is expected since the RT and HRT methods solve for both the scalar flux and current and thus do more computational work than IP and CG. 

\begin{table}
\centering
\caption{MMS order of accuracy and constants. }
\label{tab:vmms}
\begin{tabular}{cccccccccccc}
\toprule
 &  &  & \multicolumn{2}{c}{$\| \varphi - \varphi_\text{ex}\|$}  &  & \multicolumn{2}{c}{$\| \varphi - \Pi \varphi_\text{ex}\|$}  &  & \multicolumn{2}{c}{$\| \vec{J} - \vec{J}_\text{ex}\|$} \\
\cmidrule{4-5}\cmidrule{7-8}\cmidrule{10-11}
$p$ & Value & & RT & HRT & & RT & HRT & & RT & HRT \\
\midrule
\multirow{2}{*}{1} & Order & & 2.002 & 2.002 & & 2.175 & 2.175 & & 0.993 & 0.993 \\
 & Constant & & 0.608 & 0.608 & & 0.145 & 0.145 & & 0.439 & 0.439 \\
\addlinespace
\multirow{2}{*}{2} & Order & & 2.989 & 2.989 & & 2.964 & 2.964 & & 2.521 & 2.521 \\
 & Constant & & 0.396 & 0.396 & & 0.118 & 0.118 & & 0.605 & 0.605 \\
\addlinespace
\multirow{2}{*}{3} & Order & & 4.006 & 4.006 & & 4.254 & 4.254 & & 2.971 & 2.971 \\
 & Constant & & 0.309 & 0.309 & & 0.098 & 0.098 & & 0.286 & 0.286 \\
\bottomrule
\end{tabular}
\end{table}
Table \ref{tab:vmms} investigates a superconvergence property of mixed finite elements and the error behavior of the current. 
The order of accuracy and constant for the RT and HRT methods was experimentally determined by computing the error on four mesh sizes and applying logarithmic regression. 
We compare the scalar flux accuracy, the scalar flux accuracy when the exact solution is first projected onto $Y_p$ using the $L^2(\D)$ projection operator $\Pi_p$, and the accuracy of the current. 
As seen in Fig.~\ref{fig:mms}, the scalar flux converges optimally. 
For standard mixed finite element radiation diffusion, the projected error measure should superconverge at $\mathcal{O}(h^{p+2})$ and the current should be $\mathcal{O}(h^{p+1})$. 
However, on this quadratically anisotropic transport problem, the superconvergence property is lost as the projected error measure converges at only $\mathcal{O}(h^{p+1})$.
Furthermore, the current converges at only $\mathcal{O}(h^p)$ for odd polynomial degrees and $\mathcal{O}(h^{p+1/2})$ for even. 
This behavior was also seen in \citet{rtvef_olivier} for the RT and HRT VEF methods. 
We note that the RT and HRT SMM moment discretizations produce solutions that are equivalent to machine precision. 
This differs from the behavior observed for the VEF methods where the RT and HRT differed on the order of the discretization error \cite{rtvef_olivier}. 
This may be due to the ability to exactly integrate the closure terms whereas with VEF the closures cannot be integrated exactly (see Remark \ref{rem:exact_integration}). 

\subsection{Thick Diffusion Limit} \label{sec:tdl}
\begin{figure}
\centering
\includegraphics[width=.65\textwidth]{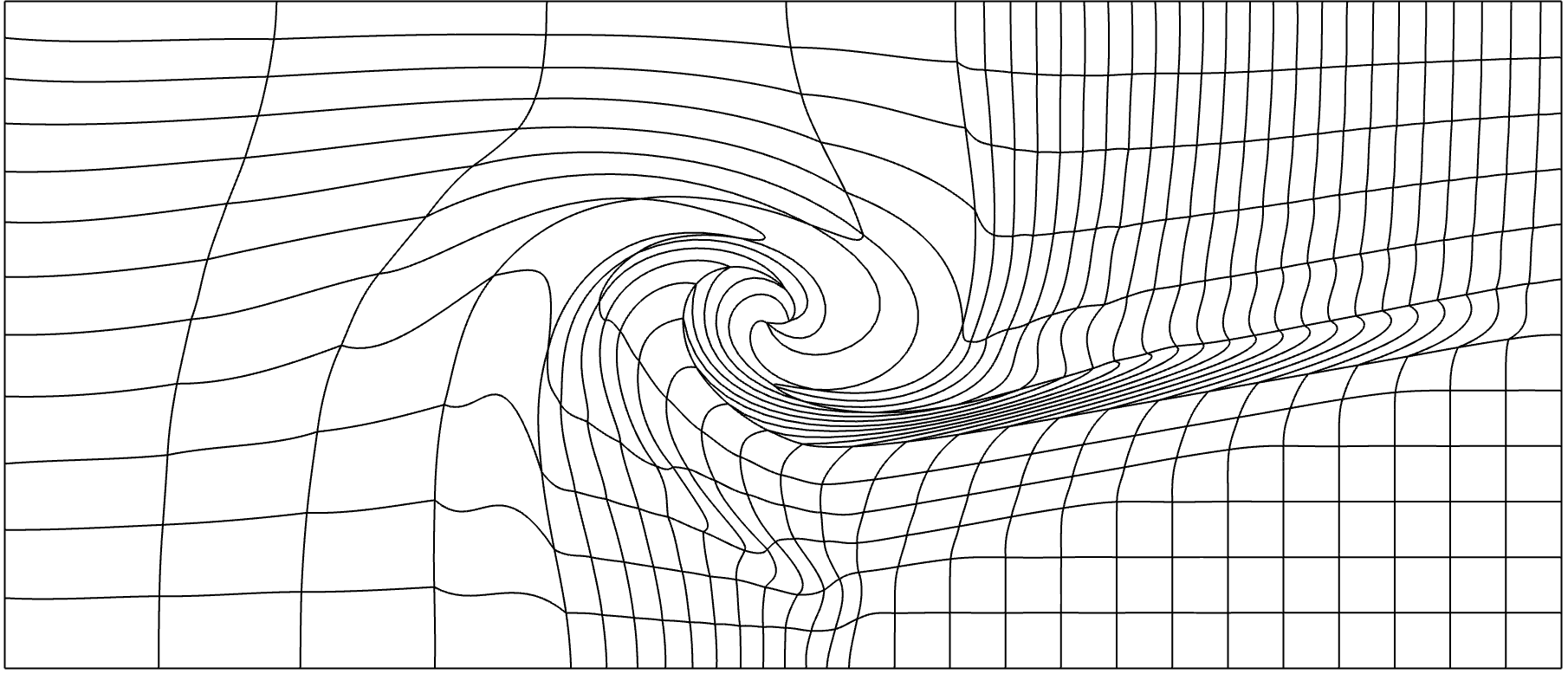}
\caption{A depiction of the triple point mesh used to stress the VEF algorithms on a severely distorted, third-order mesh. This mesh was generated with a Lagrangian hydrodynamics simulation. }
\label{fig:3point_mesh}
\end{figure}
The iterative convergence rates of the SMM algorithms are investigated in the thick diffusion limit. 
The material data are set to 
	\begin{equation}
		\sigma_t = 1/\epsilon \,, \quad \sigma_a = \epsilon \,, \quad \sigma_s = 1/\epsilon - \epsilon \,, \quad q = \epsilon \,, 
	\end{equation}
where $\epsilon \in (0,1]$ is a parameter that characterizes the mean free path and the thick diffusion limit corresponds to the limit $\epsilon \rightarrow 0$. 
We use two coarse meshes that do not resolve the mean free path to stress test the convergence of the SMM algorithm. 
The first is an orthogonal $8\times 8$ mesh with $\D = [0,1]^2$. 
The second is the triple point mesh shown in Fig.~\ref{fig:3point_mesh}. 
On the triple point mesh, the angular flux is only approximately inverted due to the presence of reentrant faces that give rise to upper block triangular entries in the streaming and collision operator.
The upper block triangular portion is iteratively lagged so that the classical transport sweep can be applied. 
Due to this approximation, iterative convergence of the moment algorithms is expected to degrade on the triple point mesh. 
On both meshes, we use Level Symmetric S$_4$ angular quadrature. 
The SMM and VEF algorithms are compared using $p=2$. 
The coupled transport-moment systems are solved with fixed-point iteration to a tolerance of $10^{-6}$. 

\begin{table}
\centering
\caption{The number of fixed-point iterations in the thick diffusion limit. }
\label{tab:tdl_orthog}
\begin{tabular}{ccccccccccccc}
\toprule
 & \multicolumn{2}{c}{IP}  &  & \multicolumn{2}{c}{CG}  &  & \multicolumn{2}{c}{RT}  &  & \multicolumn{2}{c}{HRT} \\
\cmidrule{2-3}\cmidrule{5-6}\cmidrule{8-9}\cmidrule{11-12}
$\epsilon$ & VEF & SMM & & VEF & SMM & & VEF & SMM & & VEF & SMM \\
\midrule
$10^{-1}$ & 8 & 10 & & 8 & 10 & & 8 & 10 & & 8 & 10 \\
$10^{-2}$ & 6 & 8 & & 6 & 8 & & 6 & 8 & & 6 & 8 \\
$10^{-3}$ & 4 & 5 & & 4 & 5 & & 4 & 6 & & 4 & 6 \\
$10^{-4}$ & 3 & 4 & & 3 & 4 & & 3 & 4 & & 3 & 4 \\
\bottomrule
\end{tabular}
\end{table}

\begin{figure}
\centering
\begin{subfigure}{.24\textwidth}
	\centering
	\includegraphics[width=\textwidth]{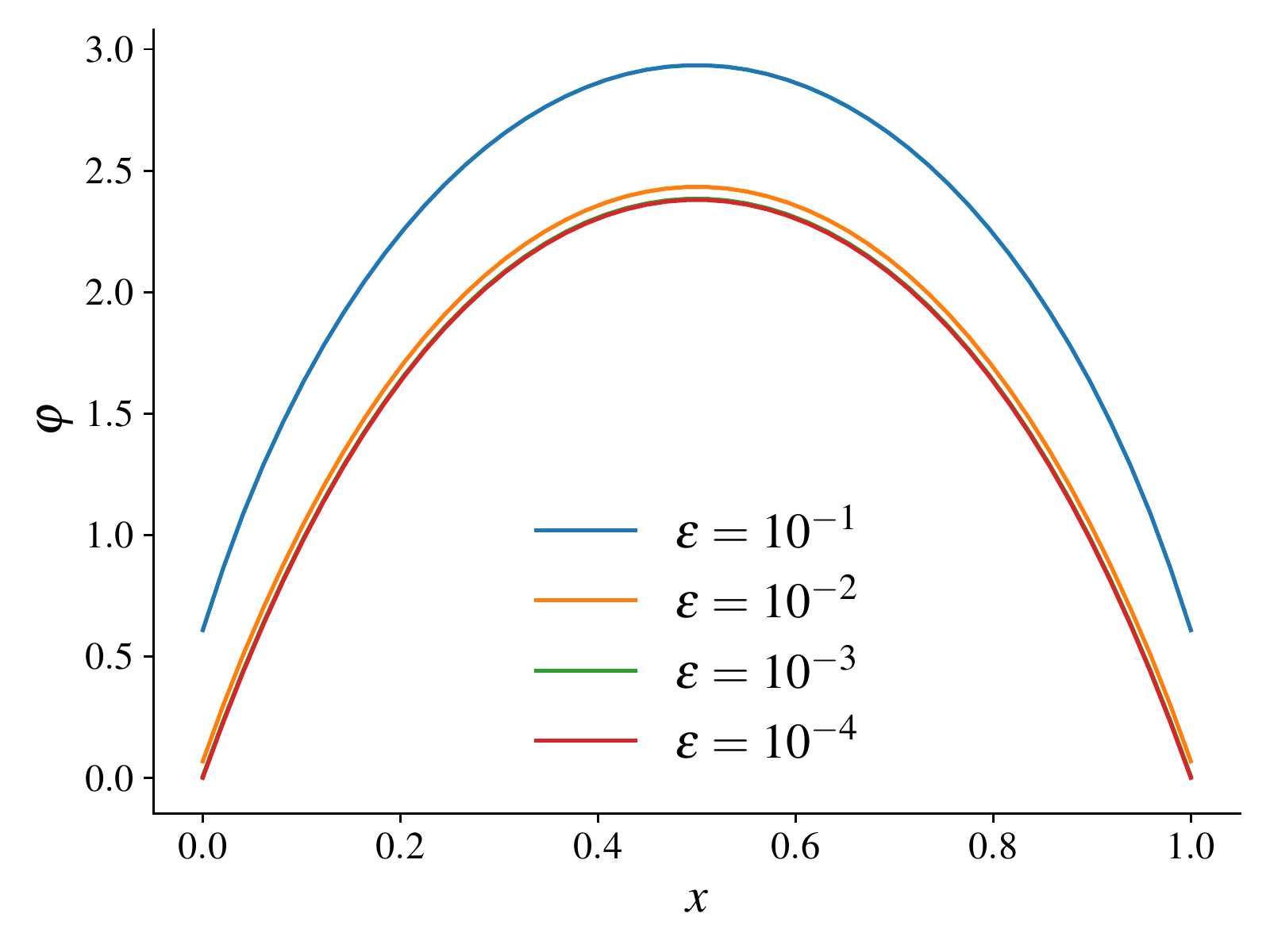}
	\caption{}
\end{subfigure}
\begin{subfigure}{.24\textwidth}
	\centering
	\includegraphics[width=\textwidth]{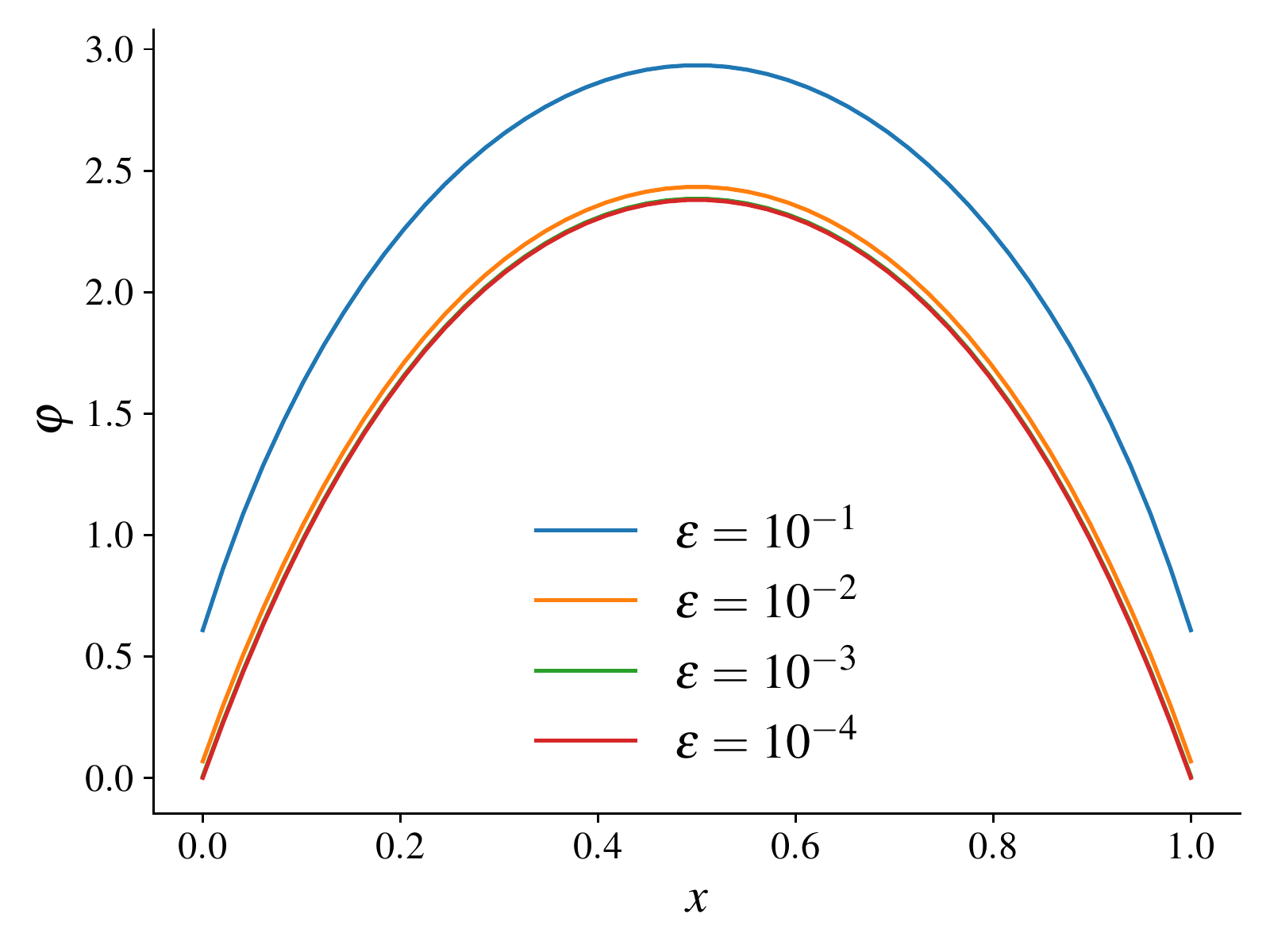}
	\caption{}
\end{subfigure}
\begin{subfigure}{.24\textwidth}
	\centering
	\includegraphics[width=\textwidth]{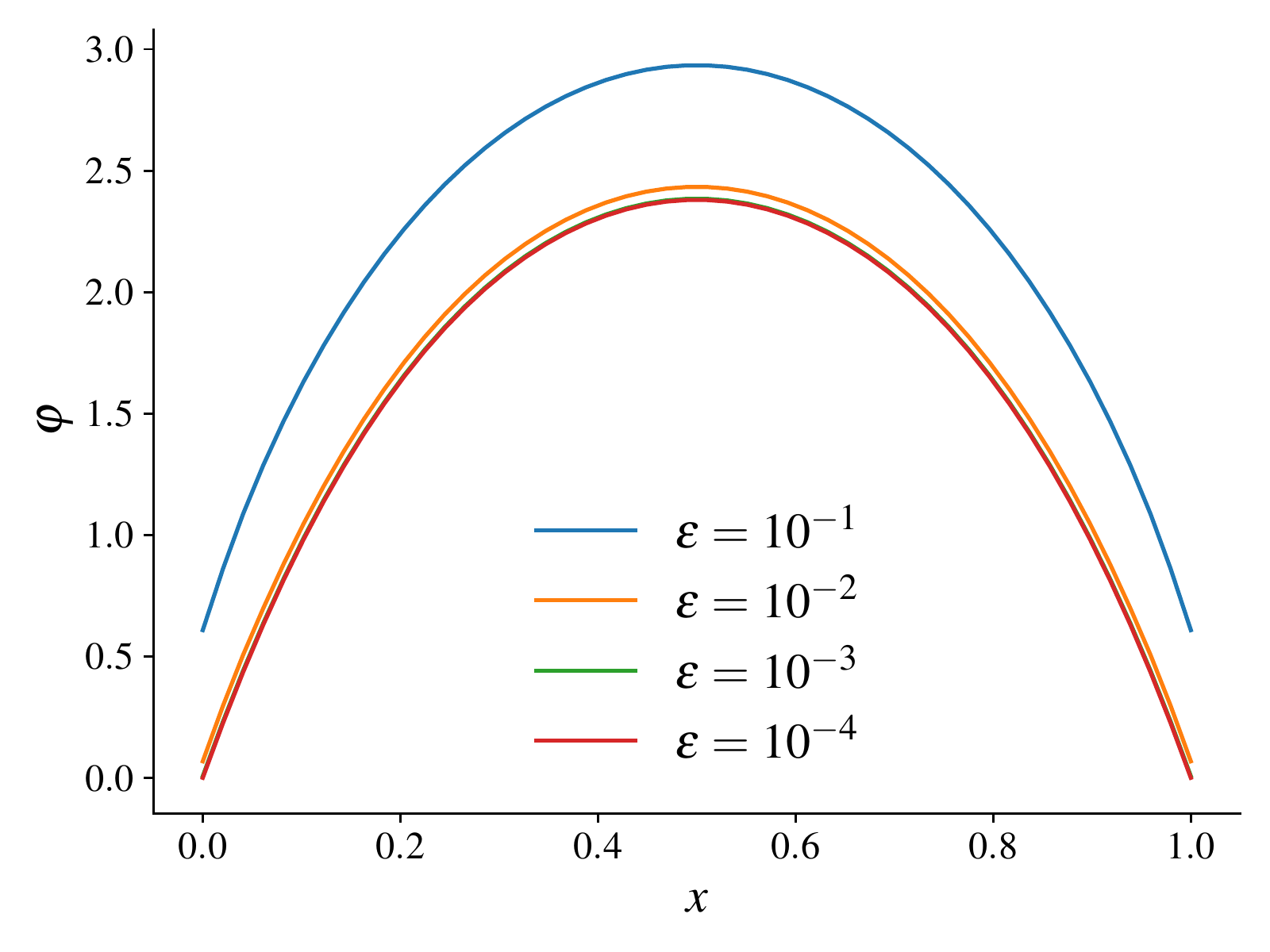}
	\caption{}
\end{subfigure}
\begin{subfigure}{.24\textwidth}
	\centering
	\includegraphics[width=\textwidth]{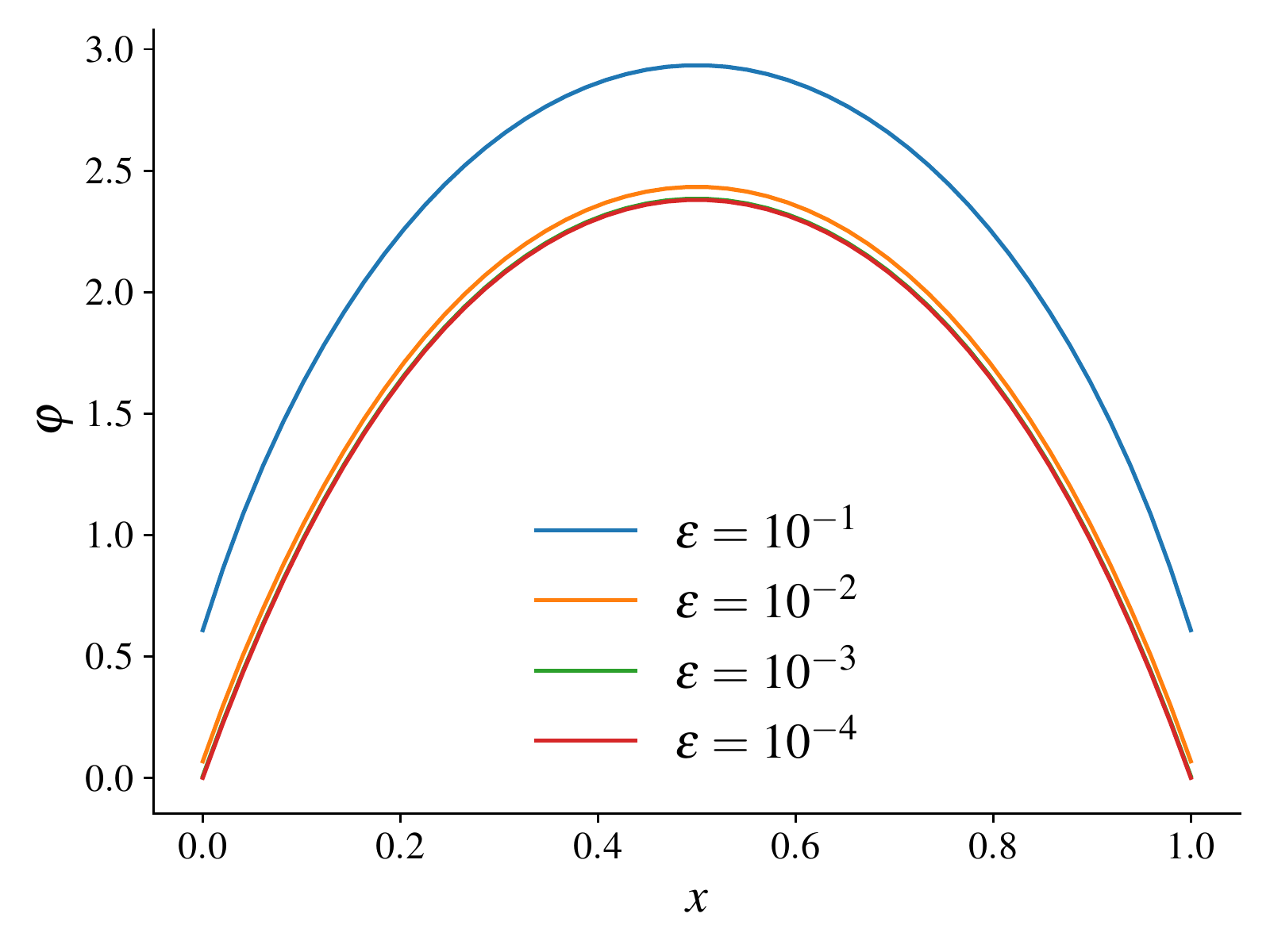}
	\caption{}
\end{subfigure}
\caption{Lineouts of the 2D solution at $y=1/2$ as $\epsilon\rightarrow 0$ for the (a) IP, (b) CG, (c) RT, and (d) HRT methods on an orthogonal $8\times 8$ mesh. The methods all converge to the asymptotic solution indicating they preserve the thick diffusion limit.}
\label{fig:tdl_lineout_orthog}
\end{figure}
Table \ref{tab:tdl_orthog} shows the number of iterations to convergence on the $8\times 8$ orthogonal mesh as $\epsilon \rightarrow 0$. 
Compared to VEF, SMM converged 1-2 iterations slower. 
All of the SMMs required the same number of iterations, indicating the algorithm is insensitive to the choice of moment discretization on this single material problem. 
Lineouts of the solution are shown in Fig.~\ref{fig:tdl_lineout_orthog} demonstrating that the methods are converging to the physically realistic diffusion solution as $\epsilon \rightarrow 0$. 

\begin{table}
\centering
\caption{Iterations to convergence in the thick diffusion limit on the triple point mesh.}
\label{tab:tdl_3point}
\begin{tabular}{ccccccccccccc}
\toprule
 & \multicolumn{2}{c}{IP}  &  & \multicolumn{2}{c}{CG}  &  & \multicolumn{2}{c}{RT}  &  & \multicolumn{2}{c}{HRT} \\
\cmidrule{2-3}\cmidrule{5-6}\cmidrule{8-9}\cmidrule{11-12}
$\epsilon$ & VEF & SMM & & VEF & SMM & & VEF & SMM & & VEF & SMM \\
\midrule
$10^{-1}$ & 18 & 21 & & 18 & 21 & & 21 & 20 & & 21 & 20 \\
$10^{-2}$ & 11 & 12 & & 11 & 12 & & 19 & 17 & & 19 & 17 \\
$10^{-3}$ & 8 & 8 & & 8 & 8 & & 13 & 13 & & 13 & 13 \\
$10^{-4}$ & 6 & 6 & & 6 & 6 & & 8 & 8 & & 8 & 8 \\
\bottomrule
\end{tabular}
\end{table}

\begin{figure}
\centering
\begin{subfigure}{.24\textwidth}
	\centering
	\includegraphics[width=\textwidth]{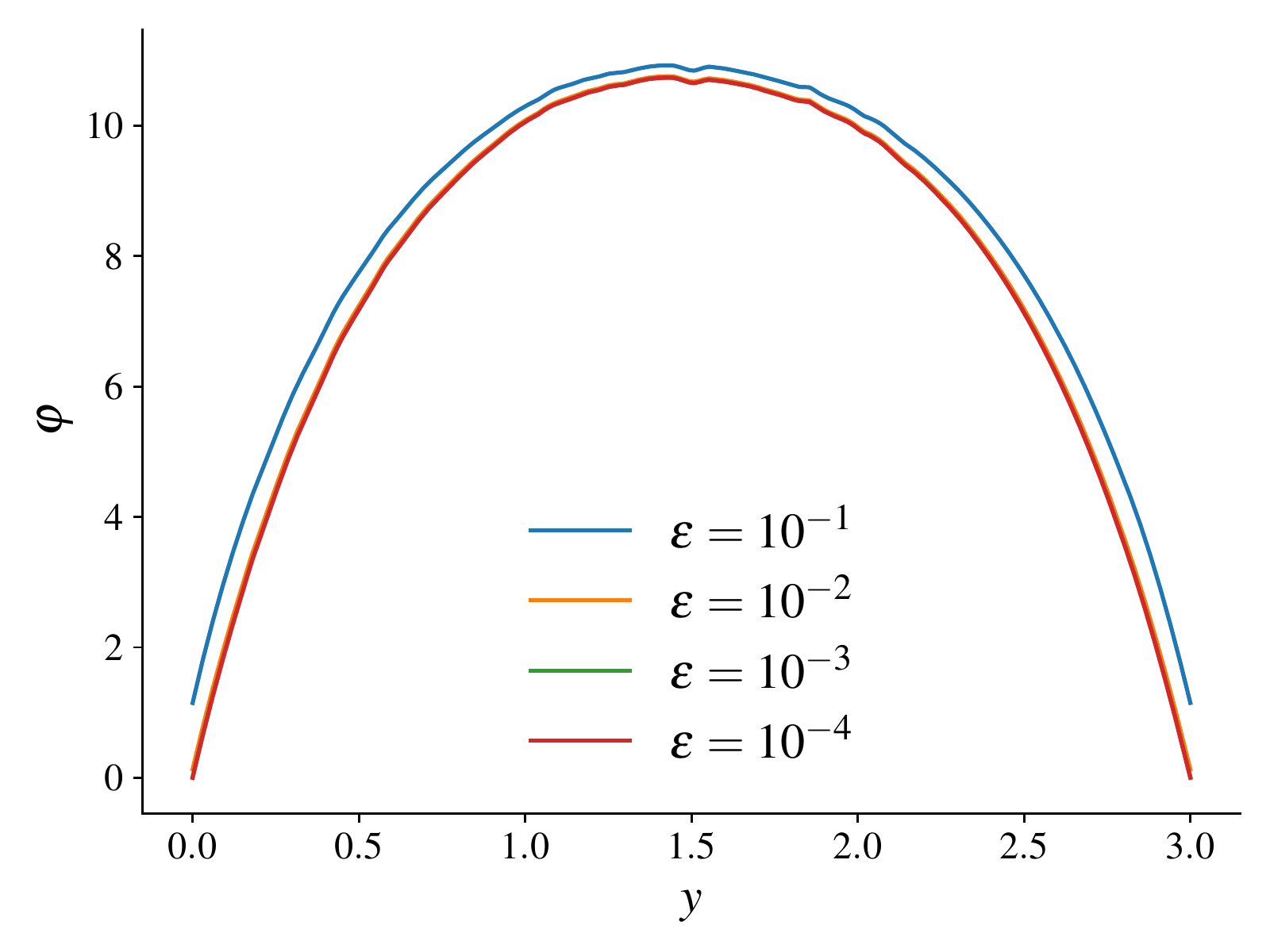}
	\caption{}
\end{subfigure}
\begin{subfigure}{.24\textwidth}
	\centering
	\includegraphics[width=\textwidth]{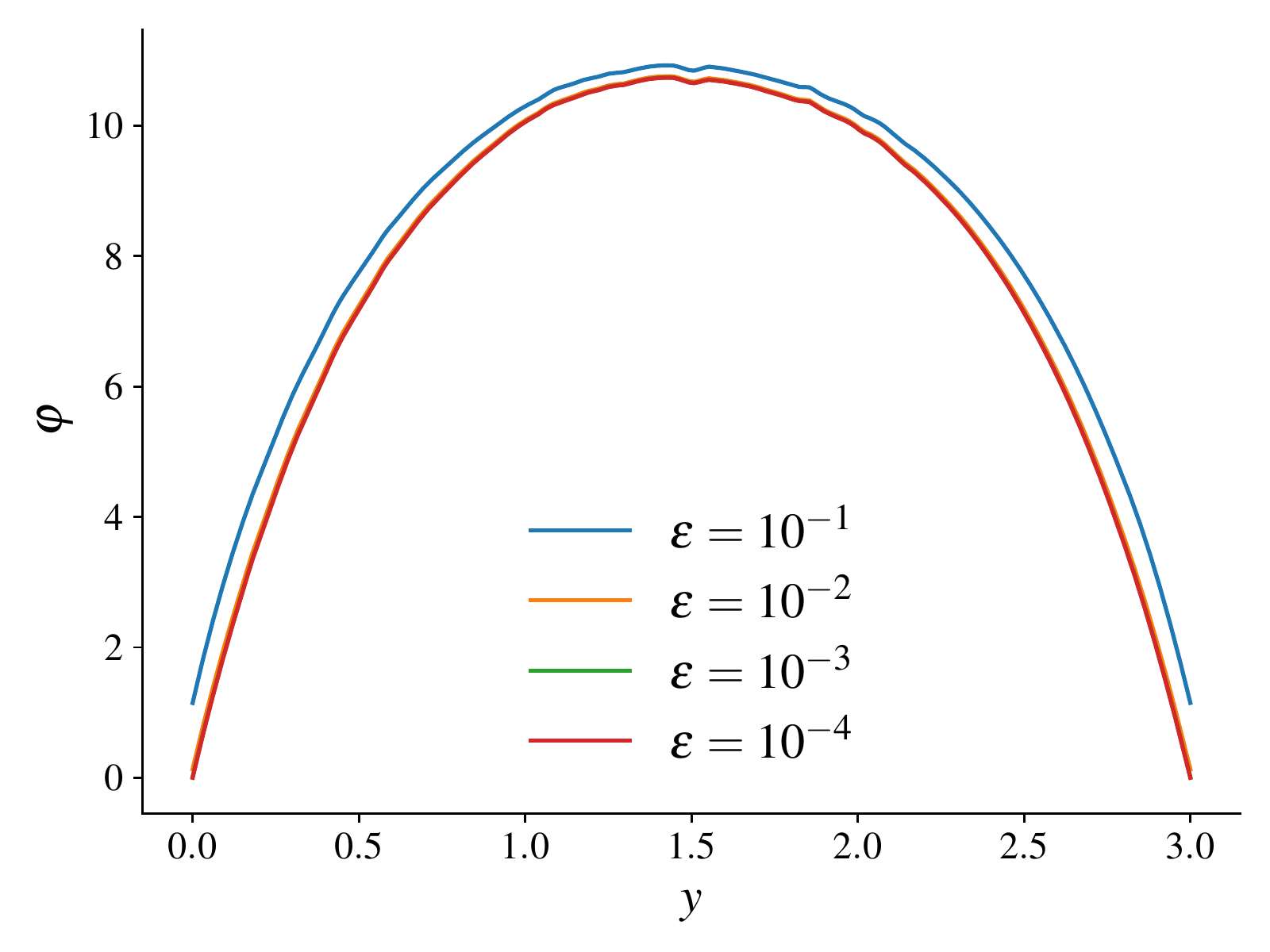}
	\caption{}
\end{subfigure}
\begin{subfigure}{.24\textwidth}
	\centering
	\includegraphics[width=\textwidth]{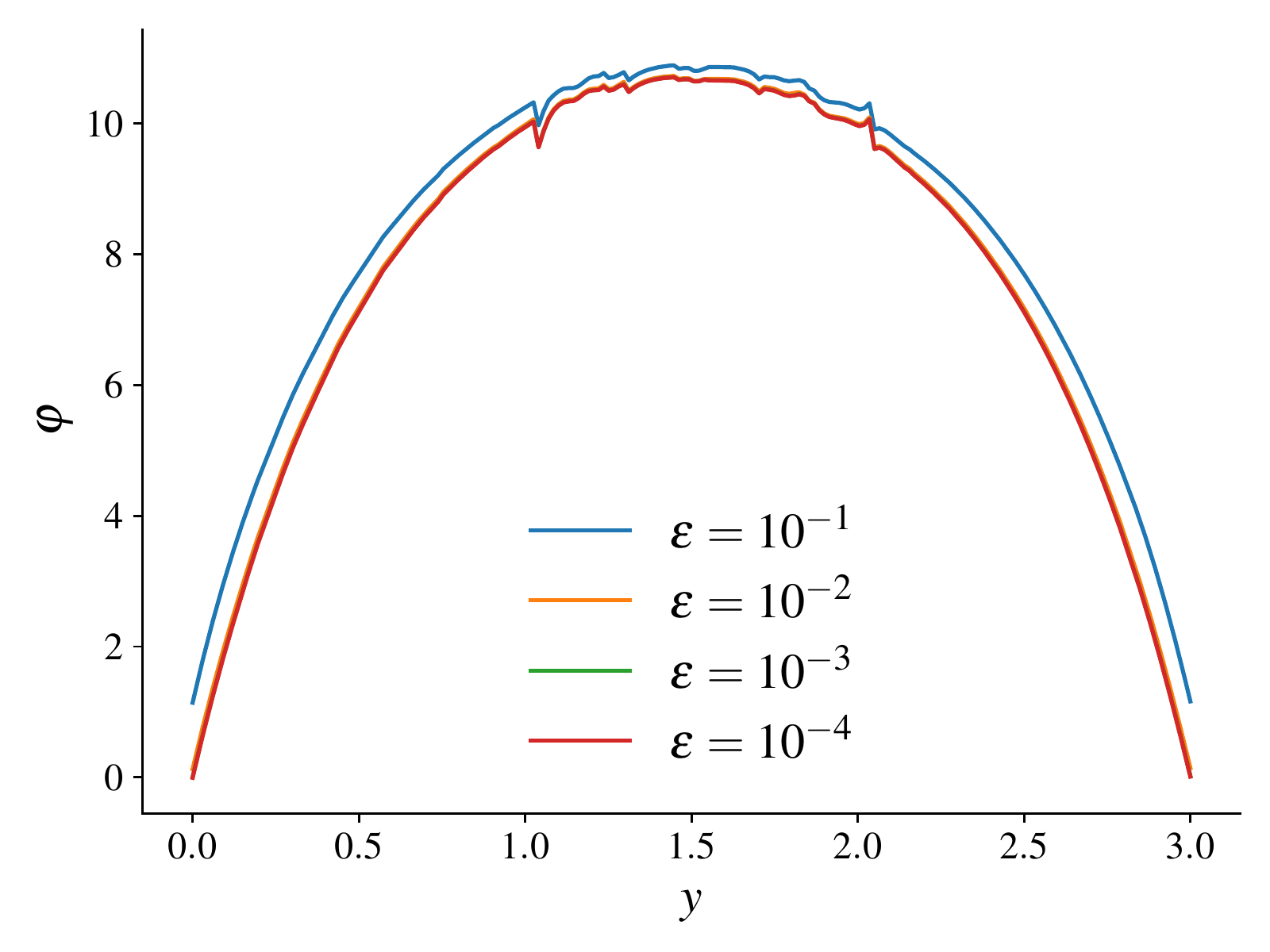}
	\caption{}
\end{subfigure}
\begin{subfigure}{.24\textwidth}
	\centering
	\includegraphics[width=\textwidth]{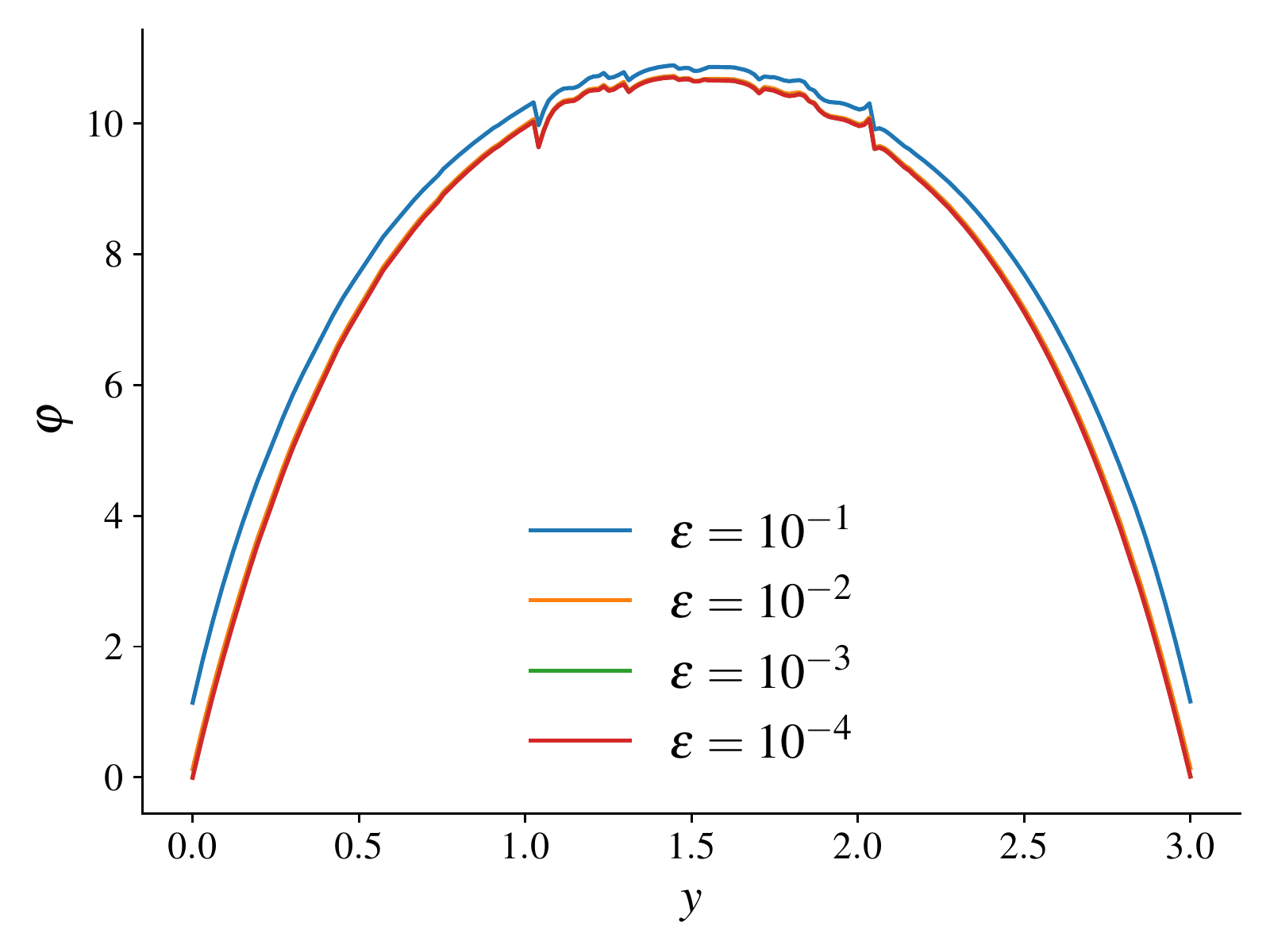}
	\caption{}
\end{subfigure}
\caption{Lineouts of the 2D solution at $x=3.5$ as $\epsilon\rightarrow 0$ for the (a) IP, (b) CG, (c) RT, and (d) HRT methods on the triple point mesh. Non-monotonic oscillations are observed due to imprinting of the highly distorted mesh onto the solution. All methods obtain a non-trivial diffusion solution. }
\label{fig:tdl_lineout_3point}
\end{figure}
Iteration counts and lineouts are presented in Table \ref{tab:tdl_3point} and Fig.~\ref{fig:tdl_lineout_3point}, respectively, for the thick diffusion limit test on the triple point mesh. 
Following the prescription of the penalty parameter used for IP VEF on this problem given in \citet[Eq.~94]{dgvef_olivier}, the IP SMM penalty parameter was scaled by a factor of 169 to ensure stability of the moment system's discretization on the highly distorted triple point mesh. 
Here, the IP and CG methods and the RT and HRT methods converge equivalently but, unlike on the orthogonal mesh problem, these two sets of methods show different performance for the intermediate values of $\epsilon$ for both VEF and SMM. 
In particular, RT and HRT SMM required five more iterations to converge compared to IP and CG for both $\epsilon = 10^{-2}$ and $\epsilon = 10^{-3}$. 
This behavior is also seen for the RT and HRT VEF methods. 
Generally, convergence was slower on the triple point mesh than on the orthogonal mesh. 
However, this discrepancy is reduced as $\epsilon \rightarrow 0$. 
The lineouts indicate that all methods attain the non-trivial, diffusion solution. 
Non-monotonic oscillations are observed at the center of the solutions due to imprinting of the highly distorted mesh. 
Qualitatively, the RT and HRT methods show more oscillations suggesting they are more sensitive to mesh distortion than the IP and CG methods which might also explain their slower convergence compared to IP and CG on intermediate values of $\epsilon$. 

\subsection{Crooked Pipe} \label{sec:cp}
\begin{figure}
\centering
\includegraphics[width=.65\textwidth]{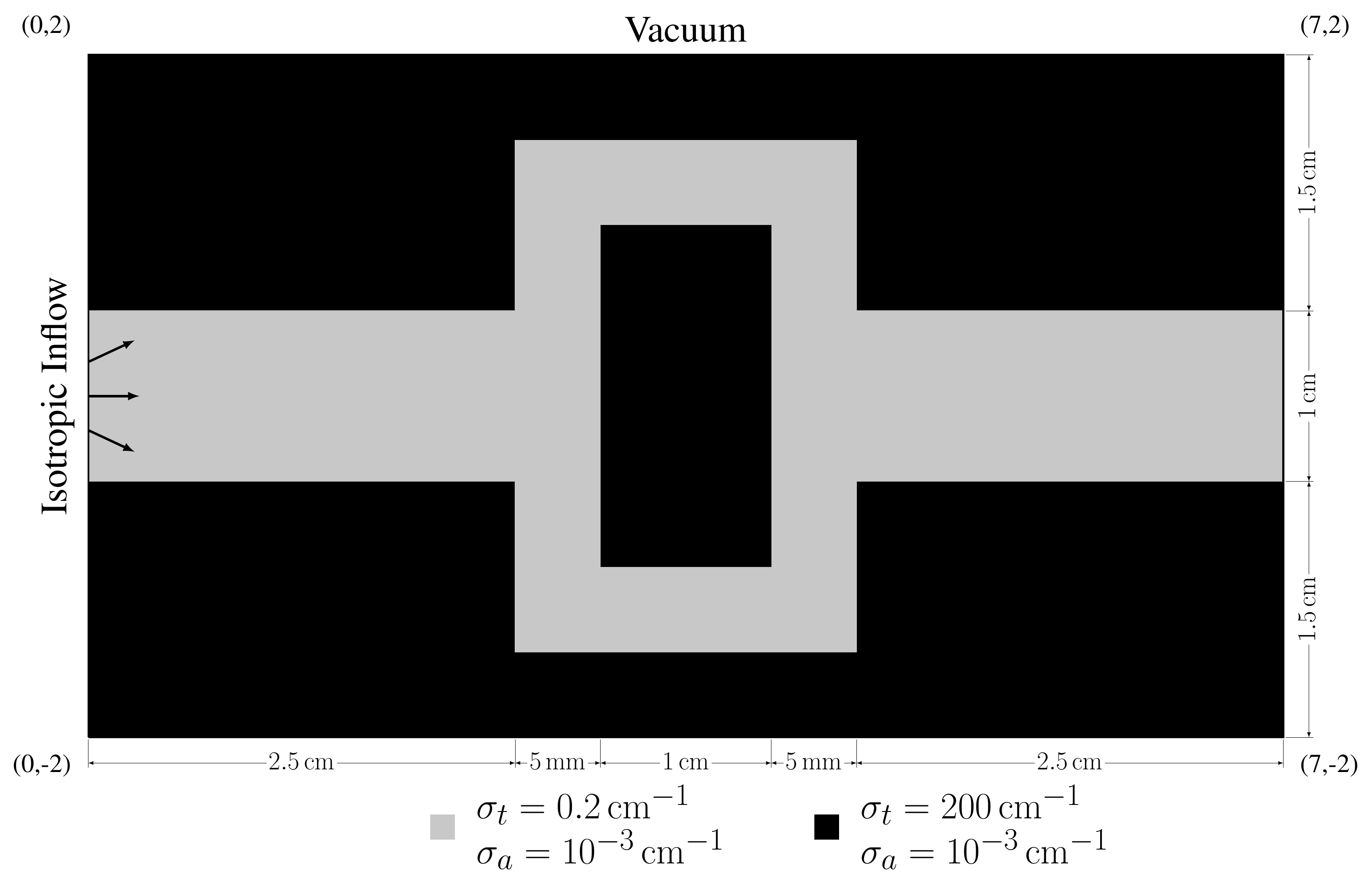}
\caption{The geometry, material data, and boundary conditions for the linearized crooked pipe problem. }
\label{fig:cp_diag}
\end{figure}

We now show convergence in outer fixed-point iterations and inner preconditioned linear solver iterations on a more realistic, multi-material problem. 
The geometry and materials are shown in Fig.~\ref{fig:cp_diag}. 
The problem consists of two materials, the wall and the pipe, which have an 1000x difference in total interaction cross section. 
Time dependence is mocked by including artificial absorption and sources that correspond to backward Euler time integration. 
The time step is set so that $c\Delta t = 10^3$ and the initial condition is $\psi_0 = 10^{-4}$. 
The absorption and source are then $\sigma_a = 1/c\Delta t = 10^{-3} \si{\per\cm}$ and $q = \psi_0/c\Delta t = 10^{-1} \si{\per\cm\cubed\per\s\per\str}$. 
The boundary conditions are set so that isotropic inflow of magnitude $1/2\pi$ enters on the left entrance of the pipe with vacuum on all other surfaces. 
A Level Symmetric S$_{12}$ angular quadrature set is used. 
The zero and scale negative flux fixup \cite{hamilton2009negative} -- a sweep-compatible fixup that zeros out negativities and rescales so that balance is preserved -- is used inside the transport sweep to ensure positivity of the angular flux unless otherwise noted. 
Timing data is presented as the minimum time recorded across five repeated runs. 

The efficiencies of the outer fixed-point and inner preconditioned linear iterations are investigated by refining in $h$ and $p$ on a uniform mesh of quadrilateral elements that are aligned with the materials in the problem. 
The outer solver is Anderson-accelerated fixed-point iteration with a small Anderson space of size two. 
Anderson acceleration is not required for convergence but does provide more uniform convergence with respect to refining in $h$. 
The outer and inner iterative tolerances were $10^{-6}$ and $10^{-8}$, respectively. 
The IP methods were preconditioned by the uniform subspace correction method, CG and HRT with AMG, and RT with block preconditioners.
The non-symmetric lower block triangular and symmetric block diagonal preconditioners were used for the RT VEF and RT SMM moment systems, respectively. 
The VEF moment systems were solved with BiCGStab while the IP, CG, and HRT SMM moment systems were solved with conjugate gradient. 
The RT SMM moment system was solved with MINRES. 
Note that the lower block triangular preconditioner used to precondition the RT VEF moment system is more expensive, and thus more effective, than the block diagonal preconditioner used for RT SMM. 
However, the lower block triangular preconditioner cannot be used with MINRES since it is not symmetric. 
The previous outer iteration's moment solution was used as an initial guess for the inner iteration so that the initial guess becomes progressively more accurate as the outer iteration converges. 

\begin{table}
\centering
\caption{Anderson-accelerated fixed-point iterations on the crooked pipe problem.}
\label{tab:cp_outer}
\begin{tabular}{cccccccccccccc}
\toprule
 &  & \multicolumn{2}{c}{IP}  &  & \multicolumn{2}{c}{CG}  &  & \multicolumn{2}{c}{RT}  &  & \multicolumn{2}{c}{HRT} \\
\cmidrule{3-4}\cmidrule{6-7}\cmidrule{9-10}\cmidrule{12-13}
 & $N_e$ & VEF & SMM & & VEF & SMM & & VEF & SMM & & VEF & SMM \\
\midrule
\multirow{4}{*}{\rotatebox{90}{$p=1$}} & 112 & 10 & 9 & & 10 & 10 & & 12 & 12 & & 12 & 12 \\
 & 448 & 11 & 11 & & 11 & 11 & & 13 & 12 & & 13 & 12 \\
 & 1792 & 13 & 12 & & 13 & 13 & & 16 & 14 & & 16 & 14 \\
 & 7168 & 14 & 14 & & 14 & 14 & & 16 & 15 & & 16 & 15 \\
\addlinespace
\multirow{4}{*}{\rotatebox{90}{$p=2$}} & 112 & 13 & 13 & & 13 & 14 & & 14 & 13 & & 14 & 13 \\
 & 448 & 14 & 14 & & 14 & 14 & & 16 & 15 & & 16 & 15 \\
 & 1792 & 15 & 14 & & 15 & 14 & & 16 & 15 & & 16 & 15 \\
 & 7168 & 15 & 14 & & 15 & 14 & & 16 & 16 & & 16 & 16 \\
\addlinespace
\multirow{4}{*}{\rotatebox{90}{$p=3$}} & 112 & 14 & 14 & & 14 & 14 & & 16 & 14 & & 16 & 14 \\
 & 448 & 15 & 14 & & 16 & 14 & & 16 & 15 & & 16 & 15 \\
 & 1792 & 15 & 16 & & 15 & 16 & & 16 & 16 & & 16 & 16 \\
 & 7168 & 15 & 18 & & 16 & 18 & & 17 & 15 & & 17 & 15 \\
\bottomrule
\end{tabular}
\end{table}
Table \ref{tab:cp_outer} presents the number of Anderson-accelerated fixed-point iterations until convergence for the four SMM algorithms each compared to the corresponding VEF algorithm with an analogous discretization for the moment system. 
Convergence was identical for the IP and CG and RT and HRT SMM methods, respectively. 
Compared to VEF, SMM converged between two iterations faster and three iterations slower. 

\begin{table}
\centering
\caption{The average number of inner, preconditioned linear iterations per outer iteration. 
}
\label{tab:cp_inner}
\begin{tabular}{cccccccccccccc}
\toprule
 &  & \multicolumn{2}{c}{IP}  &  & \multicolumn{2}{c}{CG}  &  & \multicolumn{2}{c}{RT}  &  & \multicolumn{2}{c}{HRT} \\
\cmidrule{3-4}\cmidrule{6-7}\cmidrule{9-10}\cmidrule{12-13}
 & $N_e$ & VEF & SMM & & VEF & SMM & & VEF & SMM & & VEF & SMM \\
\midrule
\multirow{4}{*}{\rotatebox{90}{$p=1$}} & 112 & 11.27 & 19.70 & & 5.18 & 7.82 & & 10.23 & 28.92 & & 4.38 & 8.50 \\
 & 448 & 11.75 & 18.75 & & 5.33 & 8.92 & & 11.36 & 37.23 & & 4.64 & 9.83 \\
 & 1792 & 11.64 & 19.62 & & 5.36 & 9.71 & & 11.24 & 36.67 & & 4.76 & 9.43 \\
 & 7168 & 11.93 & 18.73 & & 5.13 & 9.60 & & 12.35 & 37.12 & & 5.06 & 9.73 \\
\addlinespace
\multirow{4}{*}{\rotatebox{90}{$p=2$}} & 112 & 10.93 & 18.71 & & 6.57 & 10.47 & & 15.13 & 50.71 & & 6.27 & 12.46 \\
 & 448 & 11.73 & 18.67 & & 6.53 & 11.20 & & 16.41 & 54.75 & & 6.47 & 12.33 \\
 & 1792 & 12.06 & 19.00 & & 6.69 & 11.00 & & 17.24 & 57.88 & & 7.24 & 14.47 \\
 & 7168 & 12.06 & 20.40 & & 6.75 & 12.00 & & 18.35 & 60.06 & & 7.82 & 15.19 \\
\addlinespace
\multirow{4}{*}{\rotatebox{90}{$p=3$}} & 112 & 13.73 & 24.60 & & 8.00 & 13.60 & & 15.00 & 50.00 & & 5.71 & 11.93 \\
 & 448 & 14.44 & 25.53 & & 8.65 & 14.40 & & 17.12 & 51.06 & & 5.65 & 11.87 \\
 & 1792 & 15.31 & 26.00 & & 9.38 & 14.65 & & 18.00 & 52.53 & & 5.82 & 11.56 \\
 & 7168 & 15.69 & 25.26 & & 9.88 & 14.63 & & 18.50 & 54.38 & & 6.33 & 13.27 \\
\bottomrule
\end{tabular}
\end{table}
\begin{table}
\centering
\caption{The average amount of time spent solving the moment systems per outer iteration. }
\label{tab:cp_solve_time}
\begin{tabular}{cccccccccccccc}
\toprule
 &  & \multicolumn{2}{c}{IP}  &  & \multicolumn{2}{c}{CG}  &  & \multicolumn{2}{c}{RT}  &  & \multicolumn{2}{c}{HRT} \\
\cmidrule{3-4}\cmidrule{6-7}\cmidrule{9-10}\cmidrule{12-13}
 & $N_e$ & VEF & SMM & & VEF & SMM & & VEF & SMM & & VEF & SMM \\
\midrule
\multirow{4}{*}{\rotatebox{90}{$p=1$}} & 112 & 2.15 & 2.09 & & 1.62 & 1.59 & & 5.20 & 4.88 & & 2.30 & 2.24 \\
 & 448 & 7.89 & 7.48 & & 6.01 & 5.91 & & 22.33 & 24.73 & & 8.97 & 9.03 \\
 & 1792 & 30.75 & 29.62 & & 23.57 & 23.10 & & 89.39 & 101.30 & & 36.57 & 36.76 \\
 & 7168 & 125.88 & 114.97 & & 93.18 & 93.30 & & 398.94 & 402.85 & & 149.14 & 146.81 \\
\addlinespace
\multirow{4}{*}{\rotatebox{90}{$p=2$}} & 112 & 4.41 & 4.20 & & 2.79 & 2.64 & & 22.57 & 27.52 & & 3.57 & 3.49 \\
 & 448 & 17.33 & 15.58 & & 10.73 & 10.52 & & 99.97 & 126.22 & & 13.82 & 13.76 \\
 & 1792 & 72.03 & 64.75 & & 43.52 & 41.12 & & 440.46 & 525.82 & & 60.56 & 59.62 \\
 & 7168 & 286.13 & 268.64 & & 173.53 & 167.71 & & 2053.12 & 2377.09 & & 246.79 & 247.66 \\
\addlinespace
\multirow{4}{*}{\rotatebox{90}{$p=3$}} & 112 & 11.65 & 11.13 & & 6.54 & 6.16 & & 65.32 & 80.92 & & 4.92 & 5.15 \\
 & 448 & 48.72 & 46.28 & & 26.39 & 25.11 & & 300.24 & 340.78 & & 21.28 & 21.19 \\
 & 1792 & 202.41 & 181.64 & & 112.72 & 102.71 & & 1415.54 & 1518.80 & & 85.99 & 85.48 \\
 & 7168 & 903.56 & 772.37 & & 471.34 & 400.64 & & 6312.46 & 6829.70 & & 364.36 & 361.74 \\
\bottomrule
\end{tabular}
\end{table}
The average number of inner iterations per outer iteration is presented in Table \ref{tab:cp_inner}. 
All methods were scalable in $h$ and $p$. 
Recall that different solvers were used to solve the VEF and SMM moment systems for each discretization type. 
In particular, BiCGStab performs roughly twice as much work per iteration as conjugate gradient and thus the IP, CG, and HRT SMM moment solves need only take fewer than twice as many iterations as the corresponding VEF moment system to be cheaper in cost. 
This is corroborated by Table \ref{tab:cp_solve_time} where the average amount of time spent solving the moment systems is shown.
Here, the IP and CG SMM moment solves are shown to be cheaper than the corresponding VEF moment solve, especially for $p=2$ and $p=3$. 
The VEF and SMM HRT moment solves were roughly equal in cost. 
The RT SMM moment solve was more expensive than the RT VEF moment solve due to the use of differing preconditioners. 
Block diagonal-preconditioned MINRES applied to the symmetric SMM moment operator did not outperform lower block triangular-preconditioned BiCGStab applied to the non-symmetric VEF moment system. 
This indicates that it might be advantageous to use lower block triangular preconditioning and BiCGStab to solve the RT SMM moment system even though the system is symmetric. 

\begin{table}
\centering
\caption{The average amount of time spent assembling the moment systems per outer iteration. 
}
\label{tab:cp_assembly}
\begin{tabular}{cccccccccccccc}
\toprule
 &  & \multicolumn{2}{c}{IP}  &  & \multicolumn{2}{c}{CG}  &  & \multicolumn{2}{c}{RT}  &  & \multicolumn{2}{c}{HRT} \\
\cmidrule{3-4}\cmidrule{6-7}\cmidrule{9-10}\cmidrule{12-13}
 & $N_e$ & VEF & SMM & & VEF & SMM & & VEF & SMM & & VEF & SMM \\
\midrule
\multirow{4}{*}{\rotatebox{90}{$p=1$}} & 112 & 13.08 & 6.98 & & 13.07 & 5.64 & & 22.51 & 15.31 & & 22.88 & 21.71 \\
 & 448 & 49.89 & 25.08 & & 49.78 & 17.52 & & 70.83 & 50.81 & & 72.37 & 63.78 \\
 & 1792 & 195.30 & 96.38 & & 198.01 & 59.61 & & 242.67 & 182.26 & & 256.36 & 207.78 \\
 & 7168 & 784.84 & 364.82 & & 788.22 & 220.13 & & 932.15 & 678.96 & & 963.29 & 730.83 \\
\addlinespace
\multirow{4}{*}{\rotatebox{90}{$p=2$}} & 112 & 24.91 & 11.32 & & 24.81 & 8.46 & & 39.69 & 31.76 & & 43.39 & 40.83 \\
 & 448 & 95.98 & 42.90 & & 97.21 & 28.86 & & 133.31 & 101.62 & & 145.04 & 119.96 \\
 & 1792 & 384.43 & 169.98 & & 395.10 & 107.31 & & 496.14 & 366.28 & & 558.86 & 395.79 \\
 & 7168 & 1549.05 & 678.84 & & 1586.76 & 412.36 & & 1947.54 & 1344.27 & & 2157.27 & 1383.17 \\
\addlinespace
\multirow{4}{*}{\rotatebox{90}{$p=3$}} & 112 & 46.98 & 19.69 & & 46.85 & 15.15 & & 75.94 & 50.64 & & 94.05 & 74.56 \\
 & 448 & 189.96 & 74.65 & & 185.81 & 52.21 & & 274.23 & 178.82 & & 353.88 & 223.77 \\
 & 1792 & 747.61 & 289.70 & & 765.16 & 193.40 & & 1078.80 & 654.50 & & 1393.50 & 743.02 \\
 & 7168 & 3110.29 & 1131.54 & & 3152.60 & 732.44 & & 4465.49 & 2543.70 & & 5585.02 & 2630.79 \\
\bottomrule
\end{tabular}
\end{table}
Table \ref{tab:cp_assembly} shows the average cost per outer iteration of assembling the VEF and SMM moment systems. 
Recalling Remark \ref{rem:smm_advantage}, SMM is expected to have lower assembly costs due to the reduced cost associated with assembling a linear form versus a bilinear form. 
For the most refined problem with $p=3$, assembling the VEF moment system was 2.7x, 4.3x, 1.8x, and 2.2x as expensive as assembling the SMM moment system for the IP, CG, RT, and HRT discretizations, respectively. 
Assembly is a much larger cost than solving the resulting linear system iteratively. 
For example, the assembly cost for IP VEF ranges between 3.4x and 6.3x more expensive than the solve cost. 
For IP SMM, this ratio ranged only between 1.5x and 3.4x, suggesting that assembly occupies a relatively smaller portion of the IP SMM cost at each iteration compared to IP VEF. 

\begin{table}
\centering
\caption{The cost of the full algorithm to solve the linearized crooked pipe problem. 
}
\label{tab:cp_total}
\begin{tabular}{cccccccccccccc}
\toprule
 &  & \multicolumn{2}{c}{IP}  &  & \multicolumn{2}{c}{CG}  &  & \multicolumn{2}{c}{RT}  &  & \multicolumn{2}{c}{HRT} \\
\cmidrule{3-4}\cmidrule{6-7}\cmidrule{9-10}\cmidrule{12-13}
 & $N_e$ & VEF & SMM & & VEF & SMM & & VEF & SMM & & VEF & SMM \\
\midrule
\multirow{4}{*}{\rotatebox{90}{$p=1$}} & 112 & 2.08 & 1.82 & & 2.07 & 1.96 & & 2.64 & 2.49 & & 2.55 & 2.55 \\
 & 448 & 7.76 & 7.45 & & 7.74 & 7.34 & & 9.63 & 8.57 & & 9.28 & 8.56 \\
 & 1792 & 34.20 & 30.89 & & 34.71 & 32.31 & & 43.49 & 37.60 & & 42.96 & 37.39 \\
 & 7168 & 148.10 & 139.61 & & 146.48 & 137.63 & & 176.14 & 157.74 & & 170.21 & 155.24 \\
\addlinespace
\multirow{4}{*}{\rotatebox{90}{$p=2$}} & 112 & 4.30 & 4.11 & & 4.27 & 4.35 & & 5.20 & 4.67 & & 4.89 & 4.51 \\
 & 448 & 16.60 & 15.88 & & 16.53 & 15.65 & & 21.27 & 19.54 & & 19.65 & 18.43 \\
 & 1792 & 71.46 & 63.11 & & 71.05 & 61.38 & & 85.62 & 76.62 & & 78.42 & 70.60 \\
 & 7168 & 284.71 & 251.91 & & 283.51 & 247.63 & & 343.61 & 328.58 & & 310.90 & 294.44 \\
\addlinespace
\multirow{4}{*}{\rotatebox{90}{$p=3$}} & 112 & 9.19 & 8.85 & & 9.14 & 8.66 & & 12.21 & 10.35 & & 11.07 & 9.61 \\
 & 448 & 39.06 & 34.57 & & 40.05 & 33.67 & & 47.79 & 42.81 & & 42.90 & 38.48 \\
 & 1792 & 153.92 & 153.98 & & 152.63 & 151.51 & & 194.26 & 181.54 & & 172.90 & 161.17 \\
 & 7168 & 643.23 & 711.21 & & 668.00 & 700.95 & & 858.34 & 719.65 & & 749.42 & 620.24 \\
\bottomrule
\end{tabular}
\end{table}
Next, we discuss the total cost of solving the crooked pipe problem with VEF and SMM. 
Table \ref{tab:cp_total} shows the total cost of the full VEF and SMM algorithms for each discretization type on refinements of the crooked pipe. 
For a given discretization type, the algorithm that required the fewest sweeps was cheapest. 
On the most refined problem with $p=3$, the discretization choice and closure choice led to a relative standard deviation of only 10\% in total runtime. 
In other words, due to the high cost of the transport sweep, the algorithmic choices investigated here only lead to a small change in total runtime. 

\begin{table}
\centering
\caption{SMM iterations to convergence on the crooked pipe with and without a fixup.}
\label{tab:cp_fixup}
\begin{tabular}{cccccccccccccc}
\toprule
 &  & \multicolumn{2}{c}{IP}  &  & \multicolumn{2}{c}{CG}  &  & \multicolumn{2}{c}{RT}  &  & \multicolumn{2}{c}{HRT} \\
\cmidrule{3-4}\cmidrule{6-7}\cmidrule{9-10}\cmidrule{12-13}
 & $N_e$ & Fixup & No Fixup & & Fixup & No Fixup & & Fixup & No Fixup & & Fixup & No Fixup \\
\midrule
\multirow{4}{*}{\rotatebox{90}{$p=1$}} & 112 & 9 & 10 & & 10 & 10 & & 12 & 12 & & 12 & 12 \\
 & 448 & 11 & 11 & & 11 & 11 & & 12 & 12 & & 12 & 12 \\
 & 1792 & 12 & 13 & & 13 & 13 & & 14 & 14 & & 14 & 14 \\
 & 7168 & 14 & 14 & & 14 & 14 & & 15 & 15 & & 15 & 15 \\
\addlinespace
\multirow{4}{*}{\rotatebox{90}{$p=2$}} & 112 & 13 & 13 & & 14 & 12 & & 13 & 13 & & 13 & 13 \\
 & 448 & 14 & 14 & & 14 & 14 & & 15 & 15 & & 15 & 15 \\
 & 1792 & 14 & 14 & & 14 & 14 & & 15 & 16 & & 15 & 16 \\
 & 7168 & 14 & 15 & & 14 & 15 & & 16 & 16 & & 16 & 16 \\
\addlinespace
\multirow{4}{*}{\rotatebox{90}{$p=3$}} & 112 & 14 & 14 & & 14 & 14 & & 14 & 15 & & 14 & 15 \\
 & 448 & 14 & 14 & & 14 & 14 & & 15 & 15 & & 15 & 15 \\
 & 1792 & 16 & 15 & & 16 & 15 & & 16 & 16 & & 16 & 16 \\
 & 7168 & 18 & 17 & & 18 & 17 & & 15 & 16 & & 15 & 16 \\
\bottomrule
\end{tabular}
\end{table}
In light of Remark \ref{rem:nofixup}, we investigate the sensitivity of the SMMs to the use of a negative flux fixup. 
Table \ref{tab:cp_fixup} presents the number of Anderson-accelerated fixed-point iterations to convergence for the four SMMs with and without the negative flux fixup. 
Use of the fixup resulted in equivalent convergence to within $\pm 1$ iterations indicating that SMM can converge robustly without a negative flux fixup. 

\subsection{Spatial Convergence to \Sn Solution}
\begin{figure}
\centering
\begin{subfigure}{.49\textwidth}
	\centering
	\includegraphics[width=.75\textwidth]{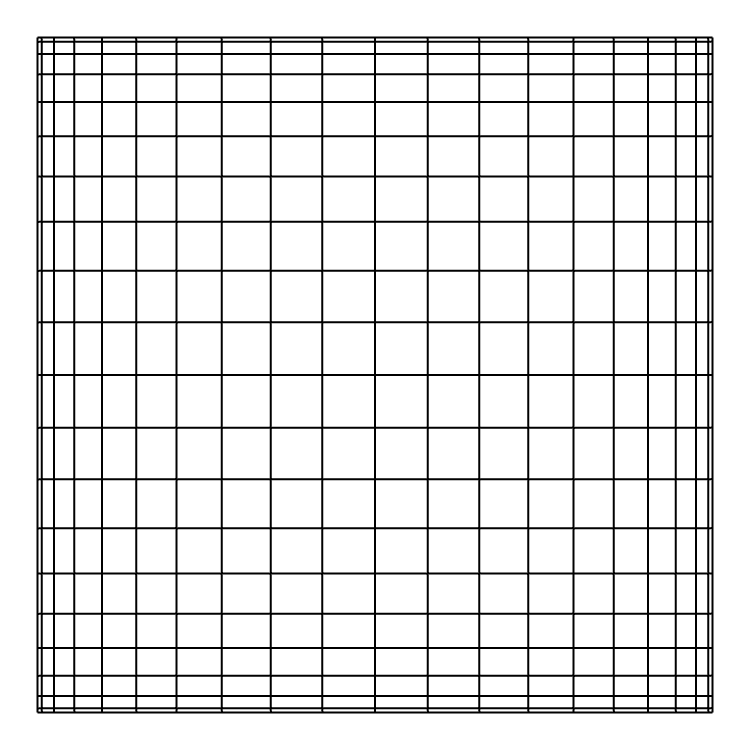}
	\caption{}
	\label{fig:chebmesh}
\end{subfigure}
\begin{subfigure}{.49\textwidth}
	\centering
	\includegraphics[width=\textwidth]{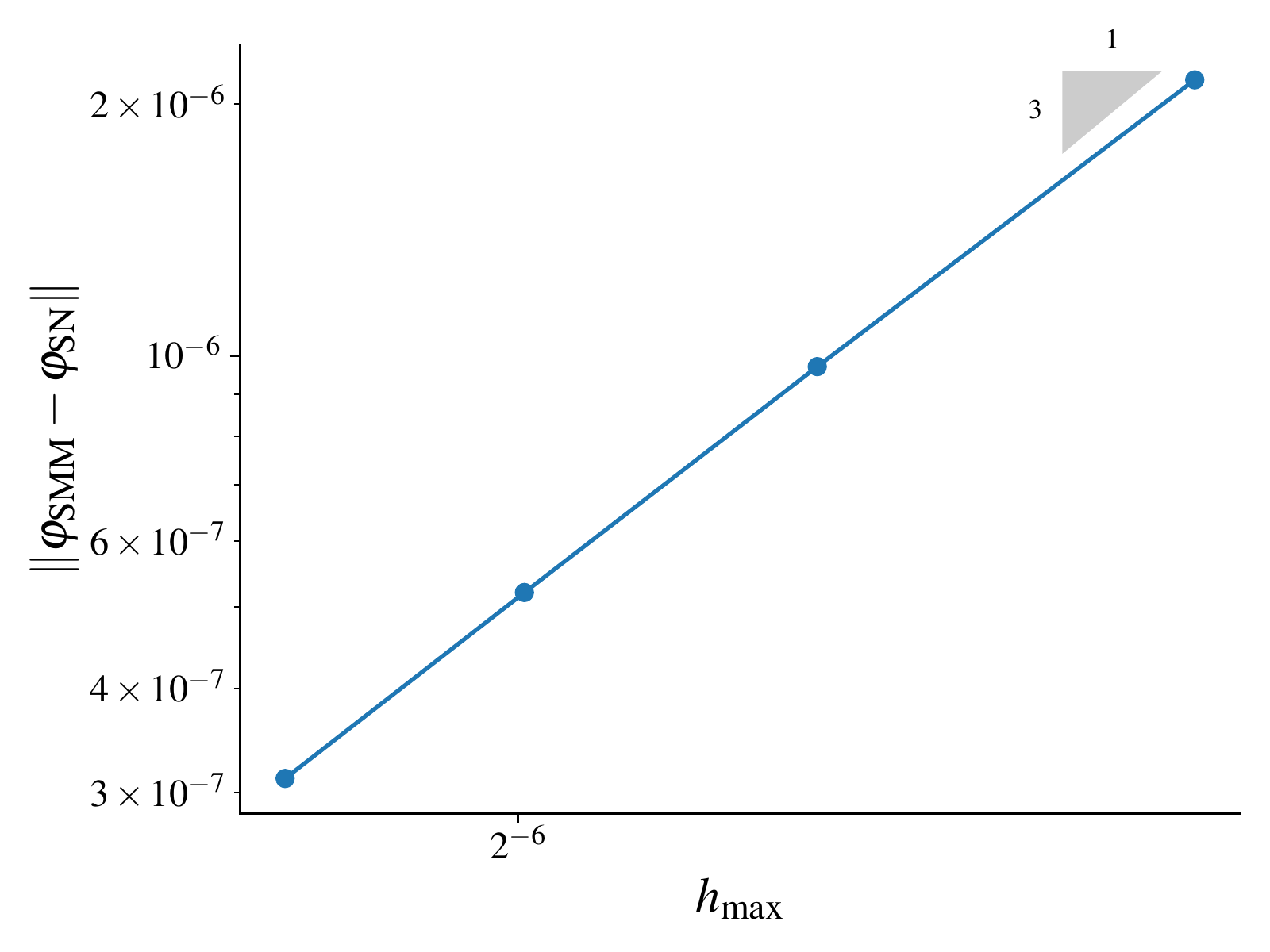}
	\caption{}
	\label{fig:dsacomp}
\end{subfigure}
\caption{(a) An example of a mesh built from a tensor product of Chebyshev points in the interval $[0,1]$ used to resolve the steep gradients in the solution at the boundary of the domain. (b) A plot of the $L^2(\D)$ norm difference between the solutions generated by the IP SMM method and a DG \Sn transport method preconditioned with DSA on the thick diffusion limit problem with $\epsilon = 10^{-2}$. Both methods used $p=2$. The solutions are compared on four meshes generated with a tensor product of 61, 81, 101, and 121 Chebyshev points in each direction. The errors are presented as a function of the maximum characteristic mesh length, $h_\text{max}$. }
\end{figure}

In this section, we compare the solutions generated by IP SMM and a reference transport solution taken to be the high-order DG \Sn discretization and DSA preconditioner of \citet{ldrd_dsa} as the mesh is refined. 
We use the thick diffusion limit problem from \S\ref{sec:tdl} with $\epsilon = 10^{-2}$ on an orthogonal mesh. 
This comparison is facilitated by the following bound on the error. 
Let the asymptotic spatial error for the SMM and \Sn solutions in isolation be written: 
	\begin{subequations}
	\begin{equation}
		\| \varphi_\text{SMM} - \varphi_\text{ex} \| = C_\text{SMM} h^{p+1} \,, 
	\end{equation}
	\begin{equation}
		\| \varphi_\text{SN} - \varphi_\text{ex} \| = C_\text{SN} h^{p+1} \,,
	\end{equation}
	\end{subequations}
where $\varphi_\text{ex}$ is the exact solution, $\varphi_\text{SMM}$ and $\varphi_\text{SN}$ the solutions generated by SMM and \Sn, respectively, $C_i$ the error constants, $h$ the mesh size, and $p$ the finite element polynomial degree. 
Using the triangle inequality, observe that 
	\begin{equation} \label{eq:snconv}
	\begin{aligned}
		\| \varphi_\text{SMM} - \varphi_\text{SN} \| &= \| (\varphi_\text{SMM} - \varphi_\text{ex}) + (\varphi_\text{ex} - \varphi_\text{SN})\| \\
		&\leq \| \varphi_\text{SMM} - \varphi_\text{ex} \| + \| \varphi_\text{SN} - \varphi_\text{ex} \| \\
		&= (C_\text{SMM} + C_\text{SN})h^{p+1} \,. 
	\end{aligned}
	\end{equation}
Thus, we expect that the difference between the solutions generated by SMM and \Sn to converge with order $p+1$. 
Note that this behavior is dependent on both solutions being in the asymptotic error regime such that their errors are well characterized by $C_ih^{p+1}$. 
Due to this, we use refinements of a non-uniform mesh built from a tensor product of the one-dimensional Chebyshev points. 
The clustering of the Chebyshev points at the endpoints allows such a mesh to capture the boundary layer in the thick diffusion limit problem. 
An example mesh is shown in Fig.~\ref{fig:chebmesh}. 

Figure \ref{fig:dsacomp} shows the $L^2(\D)$ norm between the SMM and \Sn solutions as a function of $h_\text{max}$, the maximum characteristic mesh length in the mesh. 
We use Level Symmetric S$_4$ angular quadrature and $p=2$. 
The solutions are computed on four meshes built from 61, 81, 101, and 121 Chebyshev points in each direction. 
Using $p=2$, Eq.~\ref{eq:snconv} predicts that the two solutions will converge at $\mathcal{O}(h^{3})$. 
Using logarithmic regression, the experimentally observed order of convergence was $h_\text{max}^{2.77}$. 
This result indicates that SMM does in fact converge to the \Sn solution and can converge with high-order accuracy on a problem that is smooth in space and angle, corroborating the claims made in Remark \ref{rem:linearization}. 

\subsection{Weak Scaling}
We present a weak scaling study of the IP SMM and IP VEF algorithms with $p=2$ on the crooked pipe problem from \S\ref{sec:cp}. 
Uniform refinements are used in tandem with increasing the parallel partitioning by a factor of four so that the number of degrees of freedom per processor remains fixed. 
The following results were generated on 29 nodes of the \texttt{rztopaz} machine at LLNL which has two 18-core Intel Xeon E5-2695 CPUs and 128GB of memory per node. 
Timing data is presented as the minimum time measured across three repeated runs. 

We use a parallel block Jacobi transport sweep to approximately invert the streaming and collision operator at each iteration. 
That is, each processor performs a local sweep on its processor-local domain using incoming angular flux information from the previous iteration. 
This allows each processor to sweep independently from each other at the expense of no longer exactly inverting the streaming and collision operator. 
\resp[pbj]{In many problems of interest, a parallel block Jacobi-based algorithm can outperform a full parallel sweep-based algorithm since parallel block Jacobi converges quickly when the mean free path is small, the previous time step's solution often provides a good initial guess, and parallel block Jacobi has lower communication costs compared to the full parallel sweep.}

Results are presented for SMM without a negative flux fixup applied to the angular flux since convergence was not affected by the lack of a negative flux fixup in the serial crooked pipe problem from \S\ref{sec:cp}.
However, the zero and scale fixup is still used for comparisons to VEF since a fixup is required for the VEF closures to be well defined.  

\begin{table}
\centering
\caption{Weak scaling the iterations required by IP VEF and IP SMM algorithms. 
}
\label{tab:weak_itr}
\begin{tabular}{cccccccccccccc}
\toprule
 &  & \multicolumn{2}{c}{Outer}  &  & \multicolumn{2}{c}{Min Inner}  &  & \multicolumn{2}{c}{Max Inner}  &  & \multicolumn{2}{c}{Avg.~Inner} \\
\cmidrule{3-4}\cmidrule{6-7}\cmidrule{9-10}\cmidrule{12-13}
Processors & DOF & VEF & SMM & & VEF & SMM & & VEF & SMM & & VEF & SMM \\
\midrule
1 & \num{4032} & 30 & 32 & & 2 & 3 & & 18 & 27 & & 9.27 & 13.44 \\
4 & \num{16128} & 51 & 55 & & 2 & 3 & & 22 & 32 & & 10.45 & 16.13 \\
16 & \num{64512} & 76 & 80 & & 2 & 2 & & 30 & 34 & & 10.72 & 15.62 \\
64 & \num{258048} & 135 & 135 & & 1 & 2 & & 35 & 37 & & 10.88 & 15.64 \\
256 & \num{1032192} & 261 & 277 & & 1 & 1 & & 50 & 37 & & 10.59 & 14.78 \\
1024 & \num{4128768} & 540 & 592 & & 1 & 1 & & 95 & 37 & & 10.53 & 13.63 \\
\bottomrule
\end{tabular}
\end{table}
Table \ref{tab:weak_itr} shows the inner and outer iteration counts in the weak scaling study. 
Here, the outer iteration scales with the problem size due to the use of the approximate parallel block Jacobi sweep. 
Compared to VEF, SMM required more iterations to converge, suggesting SMM is more sensitive to the parallel block Jacobi sweep than VEF. 
This sensitivity could arise from SMM's use of unnormalized closures that depend on both the transport solution's magnitude and angular shape whereas VEF uses more stable, normalized closures that depend only on the angular shape of the solution. 

The VEF and SMM moment systems were solved with uniform subspace correction-preconditioned BiCGStab and conjugate gradient, respectively. 
The previous outer iteration's solution was used as an initial guess for the inner solver. 
Both methods were scalable in terms of the average number of iterations to converge at each outer iteration. 
However, VEF shows a clear dependence on the weak scaled problem size in terms of the maximum number of inner iterations required, rising from 18 BiCGStab iterations at the smallest problem size to 95 at the largest. 
On the other hand, the maximum number of conjugate gradient iterations in the SMM algorithm varied only between 27 and 37. 
Due to the lagging of angular flux data on parallel boundaries, the closures have non-physical discontinuities in the early stages of the iteration. 
For VEF, these non-physical discontinuities affect the left hand side of the moment system, degrading the effectiveness of the preconditioned iterative solver. 
The SMM algorithm avoids this behavior because the non-physical closure terms are present in the source terms only. 

\begin{table}
\centering
\caption{Timing data for the weak scaling study of IP VEF and IP SMM. }
\label{tab:weak_time}
\begin{tabular}{ccccccccccc}
\toprule
 &  & \multicolumn{2}{c}{Avg.~Sweep}  &  & \multicolumn{2}{c}{Avg.~Moment}  &  & \multicolumn{2}{c}{Total} \\
\cmidrule{3-4}\cmidrule{6-7}\cmidrule{9-10}
Processors & DOF & VEF & SMM & & VEF & SMM & & VEF & SMM \\
\midrule
1 & \num{4032} & 0.95 & 0.95 & & 0.11 & 0.06 & & 32.37 & 32.75 \\
4 & \num{16128} & 1.12 & 1.06 & & 0.12 & 0.07 & & 64.27 & 62.74 \\
16 & \num{64512} & 1.44 & 1.31 & & 0.16 & 0.09 & & 121.74 & 112.39 \\
64 & \num{258048} & 1.62 & 1.48 & & 0.20 & 0.11 & & 247.57 & 215.81 \\
256 & \num{1032192} & 1.64 & 1.48 & & 0.21 & 0.12 & & 484.85 & 445.19 \\
1024 & \num{4128768} & 1.70 & 1.50 & & 0.23 & 0.13 & & 1048.52 & 972.07 \\
\bottomrule
\end{tabular}
\end{table}
Timing data is provided in Table \ref{tab:weak_time} for the average parallel block Jacobi sweep cost, average cost of forming and solving the moment system, and the total cost. 
Observe that the average sweep cost is higher for VEF due to the additional costs associated with applying the negative flux fixup within the sweep. 
As expected, the moment solve cost is lower for SMM compared to VEF with VEF approximately 1.7x more expensive per iteration. 
On the largest problem, the total cost of the algorithm was 8\% lower for SMM compared to VEF. 
Note that this is despite SMM requiring 52 more sweeps than VEF on this problem size. 
This discrepancy is caused by both the faster sweep times from not using a fixup and the faster SMM moment solve. 

\begin{response}[fixup]
We stress that we have used a fixup only to ensure that the moment system's closures are computed using a positive solution as discussed in Remark \ref{rem:nofixup}. 
For moment methods in a multiphysics setting, the intent is to use the moment system's solution to couple to other physics components. 
Thus, applying the fixup to the angular flux within the sweep increases the cost of the sweep but does not alter the overall simulation's robustness to negativity. 
Due to this, SMM's ability to converge without a fixup provides a computational advantage over VEF.  
\end{response} 

\begin{figure}
\centering
\includegraphics[width=.5\textwidth]{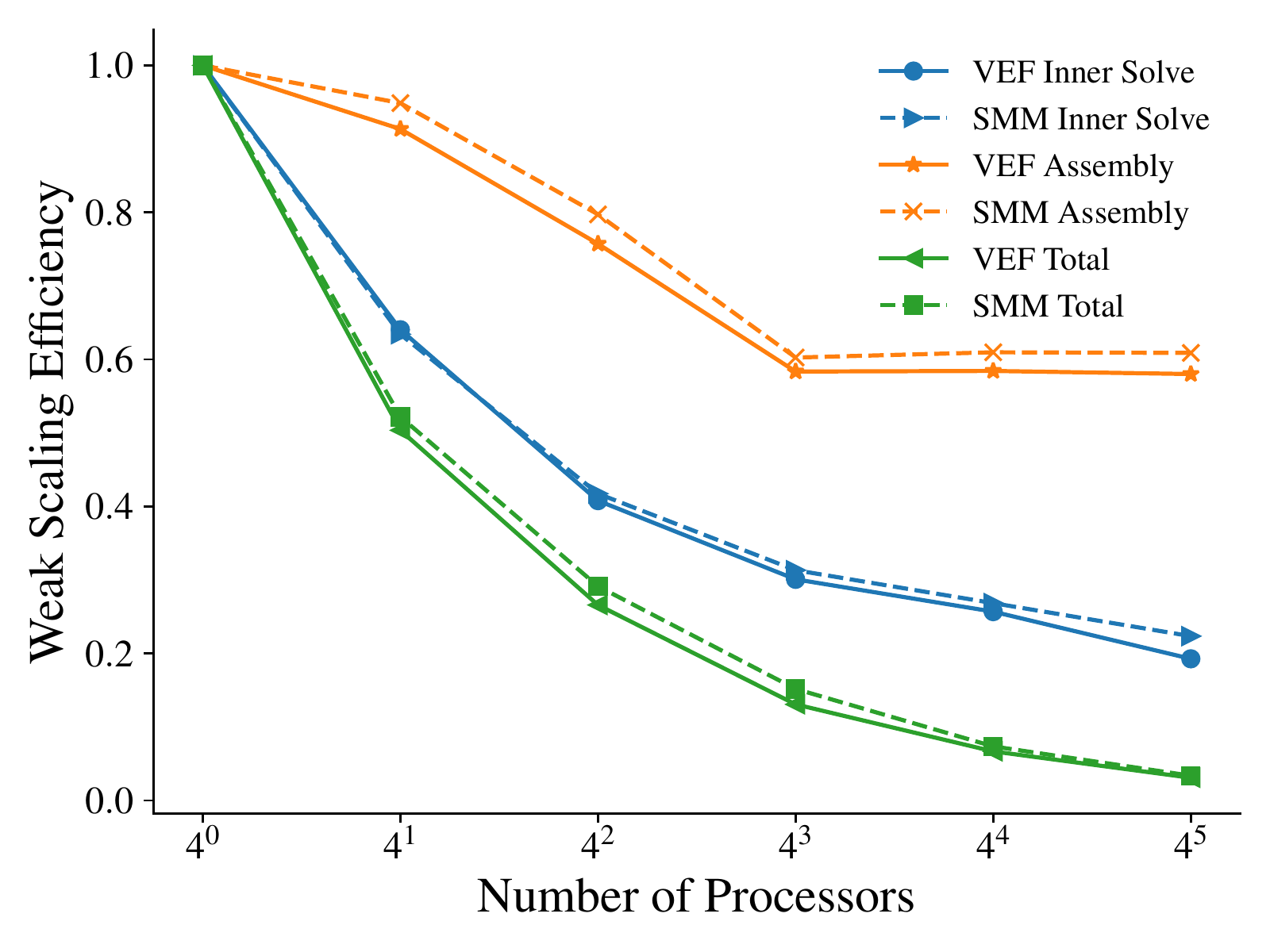}
\caption{The weak scaling efficiency for the average inner solve cost, average moment assembly cost, and total cost of the algorithm for the IP VEF and IP SMM algorithms. Scaling of the total algorithm is limited by the parallel block Jacobi transport sweep. }
\label{fig:weak_eff}
\end{figure}
Finally, we discuss the weak scaling efficiency of the algorithms. 
Let 
	\begin{equation} \label{eq:scaling}
		\varepsilon_n = \frac{\text{solve time with one processor}}{\text{solve time with $n$ processors}}
	\end{equation}
denote the weak scaling efficiency where ideal weak scaling is characterized by $\varepsilon_n = 1$. 
Due to the unavoidable communication costs associated with solving elliptic equations, ideal weak scaling is not expected. 
\begin{response}[parpart]
In addition, the problem size per processor was chosen based on the computation and storage costs associated with the high-dimensional transport sweep not the lower-dimensional moment solve.  
This leads to lower parallel efficiency in the moment solve due to the small number of moment unknowns per processor inducing higher relative communication overhead. 
\end{response}%
Figure \ref{fig:weak_eff} shows the weak scaling efficiency of the average inner solve cost, average assembly moment system assembly cost, and the total cost of the algorithm for the IP VEF and IP SMM algorithms. 
The total solve weak scales poorly due to the dependence on the number of outer iterations on the parallel partitioning from the use of the parallel block Jacobi sweep. 
Both the VEF and SMM moment solves saturate near 20\% efficient while assembly scales at 60\% for both methods. 

\subsection{Strong Scaling}
Finally, we investigate strong scaling for the IP VEF and IP SMM algorithms with $p=2$ on the crooked pipe problem from \S\ref{sec:cp}. 
The problem size was fixed at \num{28672} equally sized elements with S$_{12}$ angular quadrature. 
This led to \num{285048} scalar flux unknowns and \num{21676032} angular flux unknowns. 
Fixed-point iteration without Anderson acceleration was used. 
The outer and inner tolerances were $10^{-6}$ and $10^{-8}$, respectively. 
As with the weak scaling from the previous section, a parallel block Jacobi sweep was used to invert the streaming and collision operator and a negative flux fixup is used for VEF only. 
Timing data is reported as the minimum time measured across three repeated runs. 
Performance is characterized on a single node of the \texttt{rztopaz} machine. 

\begin{table}
\centering
\caption{Iterations to convergence for a strong scaling study of IP VEF and IP SMM.}
\label{tab:strong_itr}
\begin{tabular}{ccccccccccccc}
\toprule
 & \multicolumn{2}{c}{Outer}  &  & \multicolumn{2}{c}{Min Inner}  &  & \multicolumn{2}{c}{Max Inner}  &  & \multicolumn{2}{c}{Avg.~Inner} \\
\cmidrule{2-3}\cmidrule{5-6}\cmidrule{8-9}\cmidrule{11-12}
Processors & VEF & SMM & & VEF & SMM & & VEF & SMM & & VEF & SMM \\
\midrule
1 & 47 & 53 & & 1 & 1 & & 18 & 31 & & 8.26 & 12.60 \\
2 & 57 & 67 & & 1 & 1 & & 25 & 31 & & 9.37 & 13.96 \\
4 & 68 & 77 & & 1 & 2 & & 33 & 33 & & 10.96 & 14.99 \\
8 & 76 & 84 & & 1 & 2 & & 44 & 35 & & 11.70 & 15.64 \\
16 & 83 & 86 & & 1 & 2 & & 42 & 36 & & 11.07 & 15.98 \\
32 & 112 & 113 & & 1 & 2 & & 40 & 36 & & 11.51 & 16.28 \\
\bottomrule
\end{tabular}
\end{table}
Table \ref{tab:strong_itr} presents the outer and inner iteration data for the strong scaling study. 
SMM took between 1-10 outer, fixed-point iterations more to converge than VEF. 
As with weak scaling, VEF shows an increase in the maximum number of inner iterations as the parallel partitioning increases that is not present for the SMM moment solve. 
Both methods were uniform in average number of inner iterations. 

\begin{figure}
\centering
\includegraphics[width=.5\textwidth]{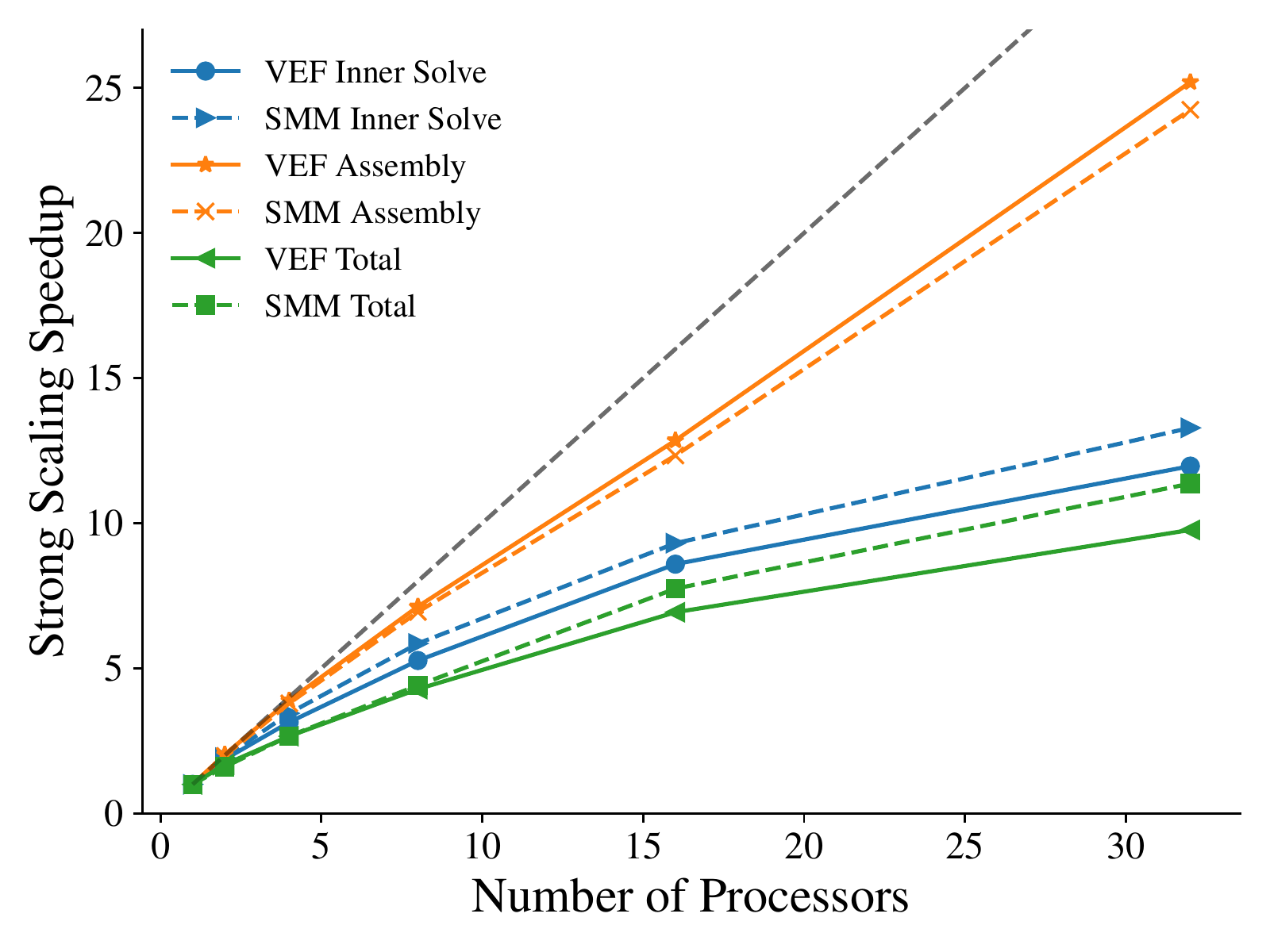}
\caption{Strong scaling speedup as a function of the number of processors on the crooked pipe problem with \num{28672} elements, S$_{12}$ angular quadrature, and $p=2$. The speedup of the average inner solve cost, average moment system assembly cost, and total cost of the algorithm are compared for the IP VEF and IP SMM algorithms. }
\label{fig:strong_eff}
\end{figure}
The strong scaling speedup, defined equivalently to $\varepsilon_n$ in Eq.~\ref{eq:scaling}, is plotted in Fig.~\ref{fig:strong_eff}. 
Here, ideal speedup is $\varepsilon_n = n$. 
The average moment system assembly cost strong scales well for both VEF and SMM due to its low communication overhead with VEF and SMM achieving speedups using 32 processors of 25.2x and 24.2x, respectively.  
The inner solve requires more communication and thus the speedup on 32 processors for the average inner solve cost was reduced to 12x and 13.2x for VEF and SMM, respectively. 
The total cost is primarily hindered by the scaling of outer iterations with processors. 
The total speedup on 32 processors was 9.8x for VEF and 11.4x for SMM. 

\section{Conclusions}
We have presented a framework for moment methods that encompasses the Variable Eddington Factor (VEF) method and the Second Moment Method (SMM). 
Both methods are iterative schemes to solve the Boltzmann transport equation centered around the use of the first two angular moments of the transport equation with suitable closures to accelerate slow-to-converge physics such as scattering. 
The VEF and SMM algorithms are differentiated only by their choice of closure: VEF uses nonlinear, multiplicative closures where SMM uses linear, additive closures. 
We demonstrated the close connection between these two algorithms both by way of algebraically reformulating the closures that define each moment method and through linearization of the VEF algorithm about a linearly anisotropic solution. 
We stress that this linearization process is actually exact since the transport equation itself is linear and thus the SMM algorithm can obtain the transport solution even in transport regimes. 

The close connection between the VEF and SMM algorithms was used to convert existing independent VEF methods -- where the closed moment system is not discretized to be algebraically equivalent to the moments of the discrete transport equation -- into corresponding discrete SMMs. 
In particular, we use algebraic manipulations and linearization to systematically convert the interior penalty (IP), continuous finite element (CG), Raviart Thomas (RT) mixed finite element, and hybridized Raviart Thomas (HRT) mixed finite element VEF moment system discretizations from \citet{dgvef_olivier} and \citet{rtvef_olivier} into discretizations of the SMM moment system. 
Due to iteratively lagging the additive SMM closures in the SMM algorithm, the discrete SMM moment systems represent standard IP, CG, RT, and HRT discretizations of radiation diffusion with transport-dependent correction sources and can thus directly leverage existing preconditioned iterative solver technology. 
We use the subspace correction preconditioner from \citet{Pazner2021} for the IP discretization, algebraic multigrid (AMG) for CG and HRT, and AMG-based block preconditioners for RT. 
The resulting discretizations and solvers are coupled to a high-order Discontinuous Galerkin discretization of the \Sn transport equations to form robust and efficient independent SMMs for linear, steady-state transport problems. 

The SMMs were verified to converge the scalar flux with optimal $\mathcal{O}(h^{p+1})$ accuracy on refinements of a third-order mesh using a manufactured, quadratically anisotropic transport solution. 
As with the mixed finite element VEF methods from \citet{rtvef_olivier}, the RT and HRT SMMs exhibited suboptimal convergence for the current. 
The iterative efficiency of the algorithms was tested in the thick diffusion limit on an orthogonal mesh and a severely distorted, third-order mesh generated with a Lagrangian hydrodynamics code. 
In both cases, the SMMs converged robustly and performed similarly to the analogous VEF methods. 

The methods were also tested on a multi-material problem designed to emulate the first time step of a thermal radiative transfer calculation. 
This problem had an 1000x difference in total cross section and sharp material discontinuities that were aligned with the mesh. 
Using solvers designed for symmetric problems, the SMM moment systems were efficiently solved uniformly with respect to the mesh size and polynomial degree. 
Using a small Anderson space of size two, the outer fixed-point iteration converged robustly as well. 
Compared to VEF, the SMMs converged nearly identically with some discretizations and problem sizes converging up to two iterations faster and others no more than three iterations slower. 
While the SMM moment system was cheaper to assemble and solve than the corresponding VEF moment system, the overall time-to-solution varied by only a relative standard deviation of 10\% across the four discretizations and two closures. 
This invariance in cost is due to the high cost of the transport sweep, which dominates the cost of the operations associated with the moment system. 
We also demonstrated that the SMMs were able to robustly converge even without a negative flux fixup, an advantage over VEF which requires a negative flux fixup to guarantee the VEF closures are well defined. 

We demonstrated that the IP SMM solution converged to the \Sn transport solution on a thick diffusion limit problem. 
Using a fixed angular quadrature rule, the SMM and \Sn solutions were compared as the mesh was refined and seen to converge at the expected, optimal order of convergence. 
This result indicates that the SMM algorithm does produce high-fidelity, transport solutions. 

Finally, we conducted weak and strong scaling studies to demonstrate the scalability of the IP SMM algorithm. 
We used a parallel block Jacobi transport sweep to invert the streaming and collision operator in parallel. 
Parallel block Jacobi allows each processor to sweep its local domain independently at the expense of no longer being an exact inversion of the streaming and collision operator. 
The SMMs generally required more outer iterations than the corresponding VEF method. 
The IP SMM inner solve was more robust to the parallel partitioning, avoiding the scaling of the maximum number of inner iterations required to converge at each outer iteration exhibited by the IP VEF method. 
These discrepancies may be due to differences in how the VEF and SMM closures interact with the parallel block Jacobi sweep. 
Due to SMM's faster moment assembly and moment solve and SMM's ability to run without a negative flux fixup, SMM was able to outperform VEF in total time-to-solution by 8\% on the largest weak scaling problem. 
On 32 processors, SMM attained a strong scaling speedup of 11.4x. 

Overall, the SMM iteration was insensitive to the choice of discretization with fixed-point iteration counts roughly independent of the choice of the discretization for the SMM moment system. 
Furthermore, the convergence rates of the SMM and VEF algorithms were similar with only minor differences in iterations to convergence. 
Thus, this study indicates that the independent moment algorithm will be robust and efficient regardless of the chosen discretization technique for the moment system or the choice of closure. 
The discretization and closure then have the flexibility to be chosen to satisfy the requirements of the larger algorithm such as to maximize computational efficiency, improve multiphysics compatibility, or simplify software design and maintenance. 
In these regards, SMM is preferred since the SMM moment systems were cheaper to assemble and solve and simpler to implement than the corresponding VEF method. 
For computational efficiency and software design purposes, we recommend the CG-based method since the CG SMM moment system had the fewest unknowns, was the cheapest to assemble and solve, and required the simplest discretization and preconditioning techniques. 
\begin{response}[ip]
The IP and RT-based SMMs are recommended for our intended multiphysics application in the hydrodynamics framework of \cite{blast} since both methods produce solutions in a compatible, discontinuous finite element space and are cheaper to solve and simpler to implement than the corresponding VEF-based algorithms. 
Furthermore, RT SMM has the distinct advantage that the left hand side diffusion operator is equivalent to the RT-based diffusion operator already implemented in \cite{blast}, allowing the design of a transport algorithm that can reuse existing software. 
We stress that suboptimal convergence of the current was observed for both RT VEF and RT SMM and thus it is not clear whether an RT-based moment method would truly improve physics compatibility with other mixed finite element-based multiphysics components. 
\end{response}%
In the future, we plan to investigate whether the algorithmic flexibility seen here for the steady-state, linear transport problem extends to more complex and numerically demanding applications such as thermal radiative transfer and radiation-hydrodynamics.  

\section*{Acknowledgments}
S.O. was supported by the U.S. Department of Energy, Office of Science, Office of Advanced Scientific Computing Research, and the Department of Energy Computational Science Graduate Fellowship under Award Number DE-SC0019323 and by the U.S. Department of Energy through the Los Alamos National Laboratory. 
Los Alamos National Laboratory is operated by Triad National Security, LLC, for the National Nuclear Security Administration of U.S. Department of Energy (Contract No.~89233218CNA000001). 

\bibliographystyle{IEEEtranN}
\bibliography{references}
\end{document}